\documentclass[final,leqno,onefignum,onetabnum]{siamltex}
\usepackage{enumerate,amsmath}
\usepackage{multirow}
\usepackage{algorithm,algorithmic}

\usepackage{psfrag}
\usepackage{fullpage}
\usepackage[crop=pdfcrop]{pstool}
\usepackage[margin=8pt,font=small,labelfont=bf,format=plain,indention=.5cm,labelsep=quad,skip=6pt]{caption}
\usepackage{float}
\usepackage{caption}
\usepackage{xmpmulti}
\usepackage{multicol}
\usepackage{epstopdf}
\usepackage{subfigure}
\usepackage{amsmath,amssymb}
\usepackage{amsopn}
\newcommand{\reviewerA}[1]{{#1}}
\newcommand{\reviewerB}[1]{{#1}}
\newcommand{\ourReReading}[1]{{#1}}

\newcommand{\restrictT}{\ensuremath{\reviewerB{\bds{R}}}}	
\newcommand{\restrict}{\ensuremath{\reviewerB{\bds{R}^T}}}	
\newcommand{\restrictArg}[1]{\ensuremath{\restrict #1 }}	
\newcommand{\restrictOrth}{\ensuremath{\reviewerB{\bds{R}^T_\perp}}}	
\newcommand{\restrictOrthT}{\ensuremath{\reviewerB{\bds{R}_\perp}}}	
\newcommand{\restrictOrthArg}[1]{\ensuremath{\restrictOrth #1 }}	
\newcommand{\restrictEntry}[1]{\ensuremath{\reviewerB{\bds{r}^T_{#1}}}}	
\newcommand{\restrictTEntry}[1]{\ensuremath{\reviewerB{\bds{r}_{#1}}}}	
\newcommand{\restrictEntryArg}[2]{\restrictEntry{#1} #2 }	
\newcommand{\restrictjArg}[1]{\restrictEntryArg{j}{#1}}	
\newcommand{\identity}{\ensuremath{\bds{I}}}	
\newcommand{\identityArg}[1]{\identity_{#1}}	
\newcommand{\prolongate}{\ensuremath{\reviewerB{\bds{P}}}}	
\newcommand{\prolongateArg}[1]{\prolongate #1}	
\newcommand{\prolongateOrth}{\ensuremath{\reviewerB{\bds{P}_\perp}}}	
\newcommand{\prolongateOrthArg}[1]{\prolongateOrth #1}	
\newcommand{\nrestrict}{\ensuremath{{\bar\ndof}}}	
\newcommand{\nrestrictPOD}{\ensuremath{q}}	
\newcommand{\totalTS}{\nptfntot}	
\newcommand{\timebasisj}{\tevbasfcssp}	
\newcommand{\timebasisjArg}[1]{\timebasis_{#1}}	
\newcommand{\timebasisjm}{\timebasisj^\timestepdummy}	
\newcommand{\timebasismArg}[1]{\timebasis^{#1}}	
\newcommand{\unrollfuncNo}{\bds{h}}
\newcommand{\unrollmfuncNo}{\unrollfuncNo^\timestepdummy}
\newcommand{\unrollfunc}[1]{{\unrollfuncNo}(#1)}
\newcommand{\unrollmfunc}[1]{{\unrollmfuncNo}(#1)}
\newcommand{\unknown}{\state}	
\newcommand{\unknownTime}{\ensuremath{\mathbf w}}	
\newcommand{\unknownTimeInit}[1]{\unknownTime^{#1,(0)}}	
\newcommand{\parametric}{\parv}	
\newcommand{\param}{\ensuremath{\parametric}}	
\newcommand{\paramTraini}[1]{\ensuremath{\bar\parametric_{#1}}}	
\newcommand{\paramDomain}{\ensuremath{\mathcal D}}	
\newcommand{\paramTrain}{\ensuremath{\{\paramTraini{i}\}_{i=1}^{\ntrain}}}	
\newcommand{\paramTrainNfunctions}{\ensuremath{\{\paramTraini{i}\}_{i=1}^{\nfunctions}}}	
\newcommand{\paramOnlinei}[1]{\ensuremath{\parametric^\star_{#1}}}	
\newcommand{\paramOnline}{\ensuremath{\{\paramOnlinei{i}\}_{i=1}^{\nonline}}}	
\newcommand{\nparam}{\pardim}	
\newcommand{\timestepdummy}{{\reviewerA{n}}}	
\newcommand{\pararealit}{\ensuremath{\reviewerA{k}}}	
\newcommand{\nstate}{\ensuremath{{\hat N}}}	
\newcommand{\memory}{\ensuremath{\alpha}}	
\newcommand{\forecastFunctionSymb}{\reviewerB{f}}	
\newcommand{\forecastFunction}[3]{\forecastFunctionSymb(#3;#2,#1)}	
\newcommand{\forecastFunctionFour}[4]{\forecastFunctionSymb_{#4}(#3;#2,#1)}	
\newcommand{\forecastFunctionjNo}{\forecastFunctionSymb_j}	
\newcommand{\forecastFunctionj}[3]{\forecastFunctionSymb_j(#3;#2,#1)}	
\newcommand{\forecastFunctionjArgs}[3]{(#3;#2,#1)}	
\newcommand{\forecastFunctionjmNo}{\forecastFunctionSymb_j^\timestepdummy}	
\newcommand{\forecastFunctionjm}[3]{\forecastFunctionSymb_j^\timestepdummy(#3;#2,#1)}	
\newcommand{\forecastFunctionjmArgs}[3]{(#3;#2,#1)}	
\newcommand{\forecastFunctionmFour}[4]{\forecastFunctionSymb_{#4}^\timestepdummy(#3;#2,#1)}	
\newcommand{\ntrain}{\ensuremath{N_\text{train}}}	
\newcommand{\nonline}{\ensuremath{N_\text{online}}}	

\usepackage{mathtools}
\newcommand{\defeq}{\vcentcolon=}

\newtheorem{remark}{Remark}[section]
\newenvironment{proofofthing}[2]{\begin{trivlist}
\item {\bf {Proof of #1 \ref{#2}}.\:}}
{$\blacksquare $\end{trivlist}}
\definecolor{Blue}{rgb}{0,0,1}
\definecolor{Red}{rgb}{1,0,0}

\newcommand{\bmat}[1]{\begin{bmatrix}#1\end{bmatrix}}

\newcommand{\vectomat}[2]{\begin{bmatrix}#1_{1}\ \cdots\ #1_{#2}\end{bmatrix}}
\newcommand{\entrytodiag}[2]{\diag(#1_{1},\ldots,#1_{#2})}
\newcommand{\umvec}{\uvec^k}
\newcommand{\vmvec}{\vvec^k}

\newcommand{\R}{\mathbb{R}}
\newcommand{\RR}[1]{\R^{#1}}
\newcommand{\RRstar}[1]{\mathbb{R}_\star^{#1}}
\newcommand{\RRplus}{\R_{+}}
\newcommand{\orthogonalMatrix}[2]{\mathbb{V}_{#2}(\RR{#1})}

\newcommand{\bds}{\mathbf}

\newcommand{\ndof}{N}

\newcommand{\state}{\bds{x}}
\newcommand{\stateExact}{\bds{x}^\star}
\newcommand{\stateNP}{\reviewerB{\bds{x}}}
\newcommand{\stateNPscalar}{x}
\newcommand{\stateNPExact}{\stateNP^\star}
\newcommand{\stateNPExactscalar}{\stateNPscalar^\star}
\newcommand{\initstate}{\state^0}
\newcommand{\initstateNP}{\stateNP^0}
\newcommand{\initstateNPscalar}{\stateNPscalar^0}

\newcommand{\stateNPEntryNo}{\reviewerB{x}}
\newcommand{\stateNPEntry}[1]{\stateNPEntryNo_{#1}}

\newcommand{\stateDummyEntryNo}{{\xi}}

\newcommand{\stateDummy}[1]{\bds{\stateDummyEntryNo}}
\newcommand{\stateContinuousDummyEntryNo}{\reviewerB{y}}

\newcommand{\stateContinuousDummy}[1]{\bds{\stateContinuousDummyEntryNo}}
\newcommand{\tintend}{T_{\mathrm{final}}}
\newcommand{\finalT}{\tintend}	
\newcommand{\tvar}{t}
\newcommand{\pardim}{p}
\newcommand{\parsp}{\mathcal{D}}
\newcommand{\parv}{\boldsymbol{\mu}}
\newcommand{\parvEntry}[1]{{p_{#1}}}
\newcommand{\ivpfunc}{\reviewerB{\bds{g}}}
\newcommand{\ivpfuncNP}{\bds{g}}

\newcommand{\tptcrs}[1]{{T}^{#1}}
\newcommand{\tixcrs}{\reviewerB{\timestepdummy}}
\newcommand{\tsscrs}{H}
\newcommand{\tsscrsk}{\reviewerB{\tsscrs}}
\newcommand{\nptcrs}{\reviewerB{M}}
\newcommand{\nCoarseIntervals}{\nptcrs}
\newcommand{\nIntervals}{\nCoarseIntervals}
\newcommand{\nTimeIntervals}{\nCoarseIntervals}
\newcommand{\tptfn}[1]{t^{#1}}

\newcommand{\tssfn}{h}
\newcommand{\nfinepercoarse}{\reviewerB{\bar m}}
\newcommand{\tssfnArg}[1]{\reviewerB{\tssfn}}
\newcommand{\tssfnk}{\tssfnArg{\timestepdummy}}

        \usepackage[utf8]{inputenc}
        \usepackage{BOONDOX-calo}
\newcommand{\timeInstanceSet}{\mathcal{t}}
\DeclareMathAlphabet{\mathpzc}{OT1}{pzc}{m}{it}
\newcommand{\timeInstanceSetInterval}[1]{\mathcal t^{#1}}
\newcommand{\timeInstanceSetCoarse}{\mathcal T}

\newcommand{\toset}[4]{\{#1\}_{#2=#3}^{#4}}

\newcommand{\crstofnNo}{{\reviewerB{\nfinepercoarse}}}
\newcommand{\crstofn}[1]{\crstofnNo #1 }
\newcommand{\nptfn}{m} 
\newcommand{\nptfnArg}[1]{\nfinepercoarse} 
\newcommand{\nptfnk}{\nptfnArg{\timestepdummy}} 
\newcommand{\nTimestepsIntervalArg}[1]{\nfinepercoarse} 
\newcommand{\nptfntot}{m}
\newcommand{\nTimestepsFine}{\nptfntot}

\newcommand{\red}[1]{\hat{#1}}
\newcommand{\reddim}{{\red{\ndof}}}
\newcommand{\redstate}{\red{\state}}
\newcommand{\redstateEntry}[1]{\red x_{#1}}
\newcommand{\redstateExact}{{\redstate}^\star}

\newcommand{\tribasmtx}{\bds{\Phi}}
\newcommand{\trialbasis}{\tribasmtx}
\newcommand{\podStateBasis}{\bds{U}}
\newcommand{\podTemporalBasisNo}{\bds{V}}
\newcommand{\podTemporalBasisj}{\podTemporalBasisNo_j}
\newcommand{\podTemporalBasisQj}{\bds{Q}_j}
\newcommand{\podTemporalBasisRj}{\bds{R}_j}
\newcommand{\podTemporalBasisVecNo}{\bds{v}}
\newcommand{\podTemporalBasisVec}[2]{\podTemporalBasisVecNo_{#1}^{#2}}

\newcommand{\snpsht}{\bds{X}}
\newcommand{\parsnpsht}[1]{\snpsht_{#1}}

\newcommand{\testmtx}{\bds{\Psi}}
\newcommand{\testbasmtx}{\testmtx}
\newcommand{\testbasis}{\testmtx}
\newcommand{\projmtx}{\testmtx}
\newcommand{\projerrorj}{\varepsilon_j}
\newcommand{\relmagj}{m_j}
\newcommand{\apx}[1]{\tilde{#1}}
\newcommand{\apxstate}{\apx{\state}}
\newcommand{\apxstateExact}{{\apxstate}^\star}

\newcommand{\tevix}{j}
\newcommand{\dimBasis}{a}
\newcommand{\dimBasism}{\reviewerB{\dimBasis}}
\newcommand{\dimfcssp}{\reviewerB{\dimBasis}}
\newcommand{\dimBasisj}{\dimfcssp}
\newcommand{\dimBasisjm}{\reviewerB{\dimBasis}}
\newcommand{\timebasis}{\bds{\Xi}}	
\newcommand{\tevbasfcssp}{\timebasis_{\tevix}}

\newcommand{\redCoordRestrict}{\hat{\bds{r}}}
\newcommand{\redCoordRestrictTwo}{\hat{\hat{\bds{r}}}}
\newcommand{\redCoordRestrictArgs}[1]{\redCoordRestrict_{#1}}
\newcommand{\redCoordRestrictTwoArgs}[2]{\redCoordRestrictTwo_{#1}^{#2}}
\newcommand{\prevtpts}{\alpha}

\newcommand{\fcarg}{\reviewerB{\bds{w}}}

\newcommand{\fccoeff}{\fcarg_{\tevix}}
\newcommand{\fccoeffm}{\fccoeff^\timestepdummy}
\newcommand{\fccoeffArg}[2]{\fccoeff(#2;#1)}
\newcommand{\fccoeffmArg}[2]{\fccoeffm(#2;#1)}
\newcommand{\smplmtxNo}{\reviewerB{\bds{Z}}}
\newcommand{\smplmtx}[2]{\reviewerB{\smplmtxNo_{#1}}}
\newcommand{\sampleMat}[1]{\smplmtx{#1}}

\newcommand{\sampleMatLocal}[1]{\smplmtxNo^{#1}}
\newcommand{\unitvec}[1]{\bds{e}_{#1}}
\newcommand{\ones}[1]{\bds{1}_{#1}}

\newcommand{\crspropsym}{\mathcal{G}}
\newcommand{\crsprop}[3]{\crspropsym ( #3;#2,#1 )}
\newcommand{\crspropforesymAll}{{\crspropsym_\mathrm{LF}}}
\newcommand{\crspropforesymAllGlobal}{\crspropsym_\mathrm{GF}}
\newcommand{\crspropBE}{\crspropsym_{\mathrm{BE}}}
\newcommand{\crspropCN}{\crspropsym_{\mathrm{CN}}}
\newcommand{\crspropfore}[3]{\crspropforesymAll( #3;#2,#1 )}
\newcommand{\crspropforeGlobal}[2]{\crspropforesymAllGlobal( #2;#1 )}
\newcommand{\crspropforesym}[2]{\crspropforesymAll_{#1}^{#2}}
\newcommand{\crspropforesymjm}{\crspropforesym{j}{\timestepdummy}}
\newcommand{\forecastQuantityArg}[3]{\delta_{#1}^{#2}(#3)}
\newcommand{\forecastQuantityVecArg}[2]{\boldsymbol\delta^{#1}(#2)}
\newcommand{\fnpropsym}{\mathcal{F}}
\newcommand{\fnpropsymBE}{\mathcal{F}_{\mathrm{BE}}}
\newcommand{\stabilityConstant}{C_\fnpropsym}

\newcommand{\overallCoarse}{A}

\newcommand{\overallFineCoarse}{B}

\newcommand{\constCoarse}{\alpha_A}
\newcommand{\constCoarseLFnj}{\alpha_{j}^n}
\newcommand{\constCoarseGFnj}{\bar\alpha_{j}^n}
\newcommand{\constCoarseComplexNoPar}{\sqrt{\nrestrict}
\lipschitzProlong\max_{j\innat{\nrestrict},n\innatZero{\nptcrs-1}}\lipschitzRestrictEntry{j}\constCoarseLFnj}

\newcommand{\constFineCoarse}{\alpha_B}
\newcommand{\constFineCoarseComplexNoPar}{(\lipschitzOrth+
\sqrt{\nrestrict}\lipschitzProlong
\max_{j\innat{\nrestrict}}\lipschitzRestrictEntry{j}(\constCoarseLFnj+1))}

\newcommand{\stabilityCoarse}{C_A}
\newcommand{\stabilityCoarseLFnj}{C_{j}^n}
\newcommand{\stabilityCoarseLFnjMod}{D_{j}^n}
\newcommand{\stabilityCoarseGFnj}{\bar C_{j}^n}
\newcommand{\stabilityCoarseComplexNoPar}{\stabilityConstant\max_{j\innat{\nrestrict},n\innatZero{\nptcrs-1}}\stabilityCoarseLFnj}

\newcommand{\stabilityFineCoarse}{C_B}
\newcommand{\stabilityFineCoarseComplexNoPar}{\stabilityConstant}

\newcommand{\normedQuantityGappyLarger}[2]{\beta_{#1}^{#2}}
\newcommand{\normedQuantityGappy}[2]{\kappa_{#1}^{#2}}
\newcommand{\normedQuantityGappyGlobal}[2]{\bar\kappa_{#1}^{#2}}
\newcommand{\normedQuantityGappySmall}[2]{\lambda_{#1}^{#2}}
\newcommand{\normedQuantityGappySmallGlobal}[2]{\bar\lambda_{#1}^{#2}}
\newcommand{\costsolve}{\tau_\fnpropsym}
\newcommand{\fnprop}[3]{\fnpropsym ( #3;#2,#1 )}

\newcommand{\apxsolFOM}[2]{\stateNP^{#1}_{#2}}
\newcommand{\Mmat}{\bds{M}}
\newcommand{\MmatEntry}[2]{m_{#1,#2}}
\newcommand{\apxsolFOMscalar}[2]{\stateNPscalar^{#1}_{#2}}

\newcommand{\finesolFOM}[2]{{\reviewerB{\bds{f}^{#1}_{#2}}}}
\newcommand{\coarsesolFOM}[2]{\bds{g}^{#1}_{#2}}
\newcommand{\itvar}{\reviewerA{k}}

\newcommand{\range}[1]{\mathrm{Ran}(#1)}

\newcommand{\rescost}{\tau_r}

\newcommand{\linearF}{\bds{A}}
\newcommand{\paramFunction}{\Theta}
\newcommand{\paramFunctionArg}[1]{\paramFunction^{#1}}
\newcommand{\paramFunctioni}{\paramFunctionArg{i}}
\newcommand{\paramFunctionMod}{\bar\Theta}
\newcommand{\paramFunctionModArg}[1]{\paramFunctionMod^{#1}}
\newcommand{\paramFunctionModi}{\paramFunctionModArg{i}}
\newcommand{\initialStateContrib}{\bar\state}
\newcommand{\initialStateContribArg}[1]{\initialStateContrib^{#1}}
\newcommand{\initialStateContribi}{\initialStateContribArg{i}}
\newcommand{\linearFnprop}{\bar\linearF_\fnpropsym}
\newcommand{\linearFnpropBE}{\bar\linearF_\text{BE}}
\newcommand{\mappingTrainingToBasis}{\bds{D}}
\newcommand{\mappingTrainingToBasisEntry}[2]{d_{#1#2}}
\newcommand{\nfunctions}{r}
\newcommand{\f}{\ivpfunc}
\newcommand{\fscalar}{g}
\newcommand{\g}{\ivpfuncNP}
\newcommand{\fArg}[3]{\f\left(#1;#2,#3\right)}
\newcommand{\fArgNoParam}[2]{\g\left(#1;#2\right)}
\newcommand{\linearFscalar}{a}
\newcommand{\linearFnpropscalar}{\bar\linearFscalar_{\mathcal F}}
\newcommand{\linearFnpropBEscalar}{\bar\linearFscalar_\text{BE}}
\newcommand{\linearCrspropscalar}{\bar\linearFscalar_{\mathcal G}}
\newcommand{\linearCrspropscalarLF}{\bar\linearFscalar_\text{LF}}
\newcommand{\linearCrspropscalarLFError}{\rho}

\newcommand{\redf}{\hat\f}
\newcommand{\redfArg}[3]{\redf\left(#1;#2,#3\right)}
\newcommand{\initredstate}{\hat\state^0}

\newcommand{\bm}[1]{{\mathbf #1}}
\newcommand{\nat}[1]{\mathbb{N}(#1)}
\newcommand{\natNo}{\mathbb{N}}
\newcommand{\natZero}[1]{\mathbb{N}_0(#1)}
\newcommand{\innat}[1]{\in\nat{#1}}

\newcommand{\innatZero}[1]{\in\natZero{#1}}
\newcommand{\innatZeroseq}[1]{=0,\ldots,#1}

\newcommand{\zero}{{\bm 0}}

\newcommand{\timestep}{\tssfn}
\newcommand{\stateVar}{{\boldsymbol \xi}}
\newcommand{\stateVarScalar}{{\xi}}
\newcommand{\stateVarEquilibrium}{{\boldsymbol \xi_e}}
\newcommand{\redstateVar}{\hat{\boldsymbol \xi}}
\newcommand{\stateContinuousVar}{\reviewerB{\boldsymbol y}}
\newcommand{\timeVar}{{ t}}
\newcommand{\timeVarArg}[1]{{\timeVar^{#1}}}
\newcommand{\paramVar}{{\boldsymbol \nu}}

\newcommand{\testbasisArg}[3]{\testbasis\left(#1;#2,#3\right)}
\newcommand{\testbasisTypical}{\testbasisArg{\redstate}{\tvar}{\parv}}

\newcommand{\localalgNo}{{\texttt{local\_basis}}}
\newcommand{\localalg}[3]{{\localalgNo}(#1,#2,#3)}
\newcommand{\podalgNo}{{\texttt{pod}}}

\newcommand{\globalForeAlg}{{\texttt{global\_forecast}}}
\newcommand{\localForeAlg}{{\texttt{local\_forecast}}}
\newcommand{\U}{\bds{U}}

\newcommand{\uvec}{\bds{u}}
\newcommand{\vvec}{\bds{v}}
\newcommand{\Sig}{\bds{\Sigma}}

\newcommand{\singularValue}{\sigma}

\newcommand{\singularValueArg}[1]{\ourReReading{\singularValue_{#1}}}
\newcommand{\V}{\bds{V}}

\newcommand{\energyCrit}{\upsilon}
\newcommand{\energyCritSet}{\Upsilon}

\newcommand{\Tm}{\tptcrs{\timestepdummy}}
\newcommand{\Tmp}{\tptcrs{\timestepdummy+1}}
\newcommand{\linearOpForecast}[1]{\reviewerB{{\bds{\Gamma}}_{#1}^\timestepdummy}}
\newcommand{\linearOpForecastGlobal}[1]{\reviewerB{\bar{\bds{\Gamma}}_{#1}^{\timestepdummy}}}
\newcommand{\forecastScalar}[2]{\reviewerB{\gamma}_{#1#2}^\timestepdummy}
\newcommand{\forecastScalarGlobal}[2]{\reviewerB{\bar\gamma}_{#1#2}^{\timestepdummy}}
\newcommand{\forecastScalarScalar}[1]{\reviewerB{\gamma}_{#1}}

\newcommand{\pararealAlgorithmModName}{\reviewerA{\texttt{parareal\_with\_initialization}}}
\newcommand{\initializeAlgorithmName}{\texttt{initialize}}

\newcommand{\jumpTolerance}{\ensuremath{\epsilon}}

\newcommand{\setOfTimeDepFunctions}{\mathcal H}

\newcommand{\setOfTimeDepFunctionsArg}[1]{\setOfTimeDepFunctions^{#1}}

\newcommand{\currentTimeInstance}{\reviewerB{i}}
\newcommand{\resLMno}{\bds{r}}
\newcommand{\resLMone}[1]{\resLMno^{#1}}
\newcommand{\resLM}[2]{\resLMone{#1}(#2)}
\newcommand{\fineFillIn}{\mathcal f}
\newcommand{\fineFillInOne}[1]{\fineFillIn(#1)}
\newcommand{\fineFillInLocalOne}[1]{\fineFillIn^{#1}}
\newcommand{\fineFillInLocal}[2]{\fineFillInLocalOne{#1}(#2)}
\newcommand{\pararealItConverge}{\ensuremath{\reviewerA{K}}}	
\newcommand{\speedupLocal}[1]{S_\mathrm{\localForeLabel\text{-}\localForeLabel}(#1)}	
\newcommand{\speedupGlobal}[1]{S_\mathrm{\globalForeLabel\text{-}\localForeLabel}(#1)}	
\newcommand{\speedupFine}[1]{S_\textit{fine}(#1)}	
\newcommand{\indicatorNo}{\mathbf 1}	
\newcommand{\indicator}[1]{\indicatorNo_{#1}}	

\usepackage{enumitem}

\newlist{Assumption}{enumerate}{1}
\setlist[Assumption]{label=A\arabic*}
\newcommand{\lipschitzProlong}{\ensuremath{M_\prolongate}}	
\newcommand{\lipschitzOrth}{\ensuremath{M_\perp}}	
\newcommand{\lipschitzProlongOrth}{\ensuremath{M_\prolongateOrth}}	
\newcommand{\minSingular}[1]{\ensuremath{\sigma_{\mathrm{min}}}\left(#1\right)}	
\newcommand{\lipschitzRestrictEntry}[1]{\ensuremath{M_{\restrictTEntry{#1}}}}	
\newcommand{\localForeLabel}{LF}
\newcommand{\globalForeLabel}{GF}
\newcommand{\BELabel}{BE}
\newcommand{\CNLabel}{CN}
\newcommand{\timeParallelError}{{e}}
\newcommand{\stabilityCoarseProp}{{\upsilon}}
\newcommand{\stabilityCoarsePropGlobal}{\bar{\upsilon}}

\mathtoolsset{showonlyrefs=true}
\title{Data-driven time parallelism via forecasting}
\author{Kevin Carlberg\thanks{Sandia National Laboratories
(\href{mailto:ktcarlb@sandia.gov}{ktcarlb@sandia.gov}).}
\and Lukas Brencher\thanks{University of Stuttgart
(\href{mailto:lukas.brencher@web.de}{lukas.brencher@web.de},
\href{mailto:haasdonk@mathematik.uni-stuttgart.de}{haasdonk@mathematik.uni-stuttgart.de},
\href{mailto:andrea.barth@mathematik.uni-stuttgart.de}{andrea.barth@mathematik.uni-stuttgart.de}).}
\and Bernard Haasdonk\footnotemark[2] \and Andrea Barth\footnotemark[2]}

\begin{document}
\setlength{\abovedisplayskip}{3pt}
\setlength{\belowdisplayskip}{3pt}
\setlength{\abovedisplayshortskip}{3pt}
\setlength{\belowdisplayshortskip}{3pt}

\maketitle
 
\begin{abstract}
This work proposes a data-driven method for enabling the efficient, stable
time-parallel numerical solution of systems of ordinary differential equations
(ODEs).  The method assumes that low-dimensional bases that accurately capture
the \textit{time evolution} of the dynamical-system state are available; these
bases can be computed from snapshot data by proper orthogonal decomposition (POD) in the case of
parameterized ODEs, for example. The method adopts the parareal framework for
time parallelism, which is defined by an initialization method, a coarse
propagator that advances solutions on a coarse time grid, and a fine
propagator that operates on an underlying fine time grid. Rather than 
employing
usual approaches for initialization and coarse propagation (e.g., a typical
time integrator applied with a large time step), we propose novel
\textit{data-driven} techniques that leverage the available time-evolution
bases. The coarse propagator is defined by a forecast (proposed in Ref.\
\cite{Carlberg_carlberg2015decreasing}) applied locally within each coarse
time interval, which comprises the following steps: (1) apply the fine
propagator for a small number of time steps, (2) approximate the state over
the entire coarse time interval using gappy POD with the local time-evolution
bases, and (3) select the approximation at the end of the time interval as the
propagated state. We also propose both local-forecast initialization (i.e.,
the local-forecast coarse propagator applied sequentially) and global-forecast
initialization (i.e., the local-forecast coarse propagator applied over the
entire time interval with global time-evolution bases).  The method is
particularly well suited for POD-based reduced-order models (ROMs). In this
case, spatial parallelism quickly saturates, as the ROM dynamical system is low
dimensional; thus, time parallelism is needed to enable lower wall times.
Further, the time-evolution bases can be extracted from readily available
data, i.e., the right singular vectors arising during POD computation. In
addition to performing analyses related to the method's accuracy, speedup, 
stability, \reviewerA{and convergence,} we also numerically demonstrate the method's performance.  Here,
numerical experiments on ROMs for a nonlinear convection--reaction problem
demonstrate the method's ability to realize \textit{near-ideal} speedups;
global-forecast initialization with a local-forecast coarse propagator
leads to the best performance.
\end{abstract}

\section{Introduction}

Two emerging trends introduce both challenges and opportunities in
computational science: (1) in future extreme-scale architectures, improved
wall-time performance must be achieved primarily by exposing additional concurrency,
and (2) the rapid increase in the volume of available physical and
computational data presents an opportunity to extract useful 
insights from these data.  
The first of these trends
		can be attributed to the stagnation of clock speeds and attendant increase
		in core counts; further, the execution time and
		energy-consumption costs of communication tend to dominate those of
		computation at extreme scale, thus creating an additional incentive for
		(communication-avoiding) concurrent computation. The second of these trends arises from an
		increase in both the number of sensors and in the quantity of generated data (e.g.,
		particle-image-velocimetry measurement systems generate full
		spatio-temporal datasets), as well as
		the increasing fidelity of physics-based simulations, which generate
		large-scale computational datasets.
Further, these trends expose a unique opportunity: integrating extreme-scale
simulation with data analytics can positively impact both data-intensive
science and extreme-scale computing \cite{doeReport}. 

This is what this work aims to accomplish: we aim to leverage
available \textit{computational data} to improve \textit{concurrency and
parallel performance}
when simulating parameterized dynamical systems.
More precisely, this work considers numerically solving large-scale systems of
parameterized ordinary differential equations (ODEs), which arise in applications ranging
from computational fluid dynamics to molecular dynamics. The above
trends have particular implications in this context.

		\subsection{Numerically solving ODEs: exposing
		concurrency}\label{sec:concurrency}
		First, the sequential nature of numerically solving ODEs (i.e.,
		numerical time integration) typically poses
		the dominant computational bottleneck, both in strong and weak scaling.
\textit{Strong scaling} refers to increasing the number of computing cores used
 to solve a problem of fixed (total) size. In the context of numerically solving
	 ODEs, strong scaling
	 is typically achieved through parallelizing `across the system' by increasing the number of processors over which the
	 problem is decomposed \textit{spatially}; this usually associates with
	 parallelizing the linear-system solve occurring 
	 within each time step for implicit time integration.\footnote{If the system
	 of ODEs is nonlinear and Newton's method is applied to solve each system of algebraic
	 equations, the
	 linear-system solve occurs at each Newton iteration within each time
	 step.} However, spatial parallelism saturates:
	 there exists a number of cores beyond which the speedup decreases due to the
	 dominance of latency and bandwidth costs over savings in sequential computation.  This maximum number of (useful) cores is proportional to
	 the problem size and defines the minimum
	 wall-time achievable by spatial parallelism alone, even in the presence of
	 unlimited computational resources. This wall-time floor can preclude computational models from
	 being employed in time-critical applications (e.g., model predictive
	 control, in-the-field analysis) that demand low simulation times.
\textit{Weak scaling} refers to simultaneously increasing both the number of computing cores and
total problem size such that the problem size per core remains fixed. In the
context of numerically solving ODEs, weak scaling is typically achieved by
refining the spatial discretization (when the ODE associates
with a spatially discretized partial differential equation) as the number of
cores used for spatial parallelism increases. However, in order to prevent 
time-discretization errors from dominating spatial-discretization errors
(and to preserve stability in the case of explicit time integration), spatial
refinement typically requires attendant temporal refinement, which leads to an increase
in the total number of time steps. This implies poor weak scaling, as the wall
time is proportional to the problem size in this case.

To this end, researchers have developed a number of
\textit{time-parallel} methods that `widen the computational front' by
exposing parallelism in the temporal
dimension.\footnote{We note that some specialized Runge--Kutta schemes
achieve parallelism `across the method' \cite{cortialThesis}; however, such
approaches are typically only useful for high-order schemes and can suffer from
dense communication patterns.} In principle, such approaches can mitigate this bottleneck, as they
can decrease the minimum realizable wall time in the strong-scaling case, and
can remove the dependence of the runtime on the total number of time steps in
the weak-scaling case. Broadly, these techniques can be categorized
\cite{gander201550} as iterative methods based on multiple shooting
\cite{nievergelt1964parallel,bellen1989parallel,saha1997parallel,lions2001parareal,Carlberg_farhat2003time},
domain decomposition and waveform relaxation
\cite{gander1998overlapping,vandewalle2013parallel}, and
multigrid
\cite{hackbuschParabolic,lubich1987multi,horton1995space,minion2011hybrid,emmett2012toward,mgrit,Neumuller},
as well as direct methods
\cite{miranker1967parallel,axelsson1985boundary,womble,worley1991parallelizing,sheen2003parallel,maday2008parallelization}.

Perhaps the most well-studied and widely adopted time-parallel method is the
\textit{parareal} technique \cite{lions2001parareal}, which can be
interpreted \cite{gander2007analysis,mgrit} as a deferred/residual-correction
scheme, a multiple-shooting method with a finite-difference Jacobian
approximation, or as a two-level multigrid method.  The parareal method
alternates between (1) time integration using a \textit{fine propagator}
executed in \textit{parallel} 
on a non-overlapping decomposition of the time domain, and (2) 
time integration using a \textit{coarse propagator} executed in
\textit{serial}
on a coarse time discretization defined by boundaries of the temporal
subdomains. 
The update formula associated with sequential coarse time integration aims to
set the discontinuities in the fine solution 
(occurring at temporal-subdomain boundaries) to zero.

The parareal method converges to the solution computed by the fine propagator;
thus the fine propagator is usually chosen to be a typical
\reviewerA{single-step} time integrator
(e.g., Runge--Kutta scheme). On the other hand, the
coarse propagator can be chosen somewhat freely; it determines the parallel
performance of the parareal method.  Desired properties in the coarse
propagator include \textit{accuracy} (i.e., it should incur small error with
respect to the fine propagator \reviewerB{to ensure fast convergence}), \textit{low cost} (i.e., its computational
complexity should not scale with the underlying fine time discretization), and
\textit{stability} (i.e., it should ensure a stable parareal recurrence).  A
primary research area within time-parallel methods aims to develop coarse
propagators that satisfy these properties. 

The most commonly used coarse propagator is simply a typical time integrator (which can
have a lower-order accuracy than the fine propagator
\cite{blouza2011parallel}) applied with coarse time steps
\cite{lions2001parareal,bal2002parareal} or an explicit time
integrator \cite{nielsen2012feasibility} (where stability may preclude use for
large coarse time steps). While straightforward to
implement, the coarse time step is typically outside the asymptotic range of
convergence for the chosen time integrator, which can hamper accuracy
and lead to slow parareal convergence. This approach can be
accelerated by additionally coarsening the spatial discretization
\cite{fischer2005parareal,FarhatCortial2006,cortial2009time}, employed
simplified physics models
\cite{baffico2002parallel,maday2003parallel,blouza2011parallel,engblom2009parallel,maday2007parareal},
or relaxed solver tolerances \cite{guibert2007adaptive}. Some authors have
also employed reduced-order models constructed `on the fly' (i.e., during the
parareal recurrence without any `offline' pre-processing step)
\cite{FarhatCortial2006,cortial2009time,ruprecht2012explicit,chen2014use}.
Instead, this work proposes employing \textit{time-evolution data} that may be
available to devise an accurate, low-cost,
stable coarse propagator. We now describe the source of these data.

		\subsection{Numerically solving ODEs: availability of
		data}\label{sec:introData}
It is often the case that \textit{data} are available about the dynamical
system of interest. These data can arise (1) from experimental analyses, (2) from
numerically solving the system of ODEs over a small time interval, or (3) from
simulating the dynamical system for different parameter instances (if the
dynamical system is parameterized), for example. 

In this work, we assume that data are available related to the \textit{time
evolution} of the dynamical-system \textit{state}. Such data
may be extracted from any of the above sources. For example, these data could
be provided from (1) experimental time traces of state variables at different
spatial coordinates, (2) a time-domain Fourier transform of the short-time ODE
numerical solution, or (3) the singular value decomposition (SVD) of the
numerical spatio-temporal solution to the dynamical system at different parameter
instances. While we focus primarily on the third data source (see Section
\ref{sec:POD}), this is not strictly required for the method to be
employed.

\subsection{Proposed methodology}

The proposed methodology adopts the data-driven forecasting method introduced
in Ref.\
\cite{Carlberg_carlberg2015decreasing} to define both the coarse
propagator and the initial solution used to `seed'
the parareal recurrence.
Given 
bases for the time-evolution of the dynamical-system state\footnote{In
practice, we apply forecasting to a \textit{restriction} of the
state.} (as discussed in
Section \ref{sec:introData} above),
the coarse propagator is
defined on a given coarse time interval by a `local forecast' as follows: (1) apply the fine propagator
for a small number of time steps, (2) apply gappy POD
	\cite{sirovichOrigGappy} with \textit{local} time-evolution bases (with support over
	the coarse time interval) to generate an
	approximation of the state over the entire coarse time interval, and (3)
	select the value of the approximated state at the end of the coarse time
	interval as the propagated state. For initialization, this `local
	forecast' can be applied sequentially; alternatively, a `global forecast' can be applied as follows: 
(1) apply the fine propagator
for a small number of time steps at the beginning of the time interval, (2) apply gappy POD
	 with \textit{global} time-evolution bases (with support over
	the entire time interval) to generate an
	approximation of the state over the entire time interval, and (3)
	select the value of the approximated state at the temporal-subdomain boundaries
	as the initial solution.

The methodology is particular well-suited for projection-based reduced-order
models (ROMs) for two reasons. \textit{First}, dynamical-system ROMs associate
with small-scale ODEs that typically must be integrated over long time
intervals. This occurs because ROMs reduce the \textit{spatial complexity}
(i.e., the cost of each linear-system solve) of large-scale dynamical systems
by reducing the number of degrees of freedom (via projection) and complexity
of evaluating nonlinear terms (e.g., via empirical interpolation
\cite{barrault2004eim,chaturantabut2010nonlinear}, empirical operator
interpolation \cite{drohmann2012reduced}, or gappy POD \cite{CarlbergGappy});
however, ROMs generally do not significantly reduce the associated
\textit{temporal complexity} (i.e., the number of linear-system solves), which
is typically proportional to the spatial dimension of the original large-scale
dynamical system. Thus, ROMs suffer from early spatial-parallelism saturation
associated with strong scaling as discussed in Section \ref{sec:concurrency}.
For example, on a compressible flow problem, the Gauss--Newton with
approximated tensors (GNAT) ROM yielded a 438 factor improvement as measured
in core--hours, but only a 6.86 wall-time speedup \cite{carlberg2011model};
spatial parallelism was saturated with only 12 cores. \textit{Second}, ROMs
already require computational data for their construction. In fact, ROMs based
on proper orthogonal decomposition (POD) already employ the third data set
described in Section \ref{sec:introData}; thus, the proposed coarse propagator
can be computed `for free' in such contexts (see Section \ref{sec:ROM}). Here,
the required time-evolution bases are easily obtained from the right singular
vectors of corresponding snapshot matrices.
Finally, we note that while we present the proposed coarse propagator and
initialization methods in the parareal context, these techniques could also be
applied to alternative time-parallel methods, e.g., PITA
\cite{Carlberg_farhat2003time}, MGRIT \cite{mgrit}.

\subsection{Outline and notation}
To summarize, contributions of this work include:
 \begin{itemize} 
 \item A novel coarse propagator based on local forecasting (Section
 \ref{sec:coarseprop}),
 \item Novel initialization methods based on both local and global
 forecasting (Section \ref{sec:initialSeed}),
\item Error analysis for the local-forecast coarse propagator (Section
\ref{sec:error}) in the general (Theorem \ref{lem:coarseError})
and ideal (Theorem \ref{thm:exactCoarse}) cases,
\item Speedup analysis for all proposed methods (Section \ref{sec:speedup}) in
the general (Theorems \ref{thm:localForecastGen}--\ref{thm:globalForecastGen})
and ideal (Theorems
\ref{thm:idealSpeedup}--\ref{thm:idealSpeedupGlobal}) 
cases,
\item Stability analysis (Section \ref{sec:stability}) of the local-forecast coarse
propagator (Lemma \ref{lem:coarseStable}) and the resulting parareal
recurrence (Theorem \ref{thm:finalStability}),
\item \reviewerA{Convergence analysis (Section \ref{sec:convergence}) of the
local-forecast coarse propagator within the parareal recurrence (Corollary \ref{cor:convergenceProposed}),}
\item Descriptions of how the required method ingredients can be computed via
POD (Section \ref{sec:POD}) for parameterized ODEs (Section \ref{sec:FOM}) and
reduced-order models (Section \ref{sec:ROM}), and
\item Numerical experiments (Section \ref{sec:experiments}) that both highlight the
practical benefits of the proposed methodology and illustrate the theoretical
results.
 \end{itemize}

The paper is structured as follows.
Section \ref{sec:timeparallel} 
introduces the parareal method, 
Section \ref{sec:proposed} describes the proposed methodology, including
algebraic techniques for data-driven global (Section \ref{sec:forecasting})
and local (Section \ref{sec:localforecast})
forecasting, and their application as coarse propagators (Section
\ref{sec:coarseprop}) and initialization
methods (Section \ref{sec:initialSeed}). Section
\ref{sec:analysis} analyzes the proposed method in terms of 
accuracy (Section \ref{sec:error}), cost/speedup (Section
\ref{sec:speedup}), stability (Section \ref{sec:stability})\reviewerA{, and
convergence (Section \ref{sec:convergence})}. Section
\ref{sec:POD} describes how the ingredients of the proposed methodology can be
computed for parameterized ODEs (Section \ref{sec:FOM}) and reduced-order models
(Section \ref{sec:ROM}) using proper orthogonal decomposition (POD),
which is closely related to the singular value decomposition (SVD). Section
\ref{sec:experiments} provides numerical experiments that assess the
performance of the proposed technique in practice. Finally, Section
\ref{sec:conclusions} concludes the manuscript, Appendix \ref{sec:proofs} contains all
proofs, \ourReReading{Appendix \ref{sec:idealSpeedupsNewtonSolver} provides some additional aspects on using 
forecasting for Newton-initialization.}

In the remainder of this paper,  matrices are denoted by capitalized bold
letters, vectors by lowercase bold letters, scalars by unbolded letters. The
columns of a matrix $\bm{A}\in\RR{m\times {n}}$ are denoted by
$\bm{a}_i\in\RR{m}$, $i\innat{{n}}$ with $\nat{n}\defeq\{1,\ldots,
n\}$ such that $\bm{A}\defeq\left[\bm{a}_1\ \cdots\
\bm{a}_{n}\right]$.  The scalar-valued matrix elements are denoted
by $a_{i,j}\in\RR{}$ such that  $\bm{a}_j\defeq\left[a_{1,j}\ \cdots\
a_{m,j}\right]^T$, $j\innat{{n}}$\ourReReading{; we similarly denote
the elements of a vector as $\bm{a}\defeq\left[a_{1}\ \cdots\
a_{m}\right]^T$.} We also define
$\natZero{n}\defeq\{0,\ldots, n\}$.

\section{Time parallelism and parareal}\label{sec:timeparallel}
We consider initial value problems for systems of (possibly nonlinear)
ordinary differential equations (ODEs) of the form
\begin{align}
\label{eq:ODE} 
\frac{d}{dt}
{\stateNPExact} (\tvar) &= \fArgNoParam{\stateNPExact}{\tvar},\quad
\stateNPExact (0) = \initstateNP \ ,
\end{align}
where
$\tvar \in [0, \tintend]$ denotes time with $\tintend \in \RRplus{}$ the
final time,
$\stateNPExact:\RRplus \rightarrow \R^{\ndof}$ denotes the state implicitly
defined as the (exact) solution to problem~\eqref{eq:ODE}, 
$\initstateNP\in \R^{\ndof}$ denotes the initial state, and
$\ivpfuncNP: \R^{\ndof} \times \RRplus \rightarrow
\R^{\ndof}$ with $(\stateVar;\timeVar)\mapsto \g(\stateVar;\timeVar)$ denotes the
velocity\reviewerA{, which may be linear or nonlinear in its first argument}.
Time-parallel methods constitute one approach to improve wall-time performance
when numerically solving such problems.  We now introduce the
parareal method, which we consider in this work.

First\reviewerB{, and without loss of generality,} we introduce a \reviewerB{uniform} \textit{fine} time discretization
characterized by time step 
$\tssfnk\reviewerB{=\finalT/\nptfn}$ and time instances
\reviewerB{$\tptfn{\timestepdummy } = \timestepdummy \tssfnk$, $\timestepdummy \innatZero{\nTimestepsFine}$},
where
$\nTimestepsFine\in\natNo$ denotes the number of total time instances beyond
the initial time $\tptfn{0} = 0$ such that the final time corresponds to $\tptfn{\nTimestepsFine}
= \finalT$. We denote the set of time instances
associated with this discretization as
$\timeInstanceSet\defeq\toset{\tptfn{\timestepdummy }}{\timestepdummy }{0}{\nTimestepsFine}$.
We introduce a `fine propagator'
\ourReReading{ $\fnpropsym:\RR{\ndof}\times \timeInstanceSet\times \timeInstanceSet\rightarrow\RR{\ndof}$  
with $(\stateVar;\timeVarArg{i},\timeVarArg{j})\mapsto\fnprop{\timeVarArg{j}}{\timeVarArg{i}}{\stateVar}$ }
that acts on this discretization and propagates a
state $\stateVar\in\RR{\ndof}$ defined at time
$\timeVarArg{i}$ to 
time $\timeVarArg{j}$ with $j\geq i$. This propagator 
satisfies 
 \begin{equation} \label{eq:fineIntermediate}
\reviewerB{
\fnprop{\timeVarArg{k}}{\timeVarArg{i}}{\stateVar} =
\fnprop{\timeVarArg{k}}{\timeVarArg{j}}{\cdot}\circ
\fnprop{\timeVarArg{j}}{\timeVarArg{i}}{\stateVar}
,\quad 0\leq i\leq j  \leq
k\leq\nTimestepsFine}
 \end{equation}
\noindent and typically corresponds to
the application of a \reviewerA{single-step} time integrator (e.g., 
Runge--Kutta scheme) to numerically solve problem~\eqref{eq:ODE}. \reviewerB{For example,
the backward-Euler fine propagator $\fnpropsymBE$
implicitly satisfies
$
\fnpropsymBE(\stateVar;\timeVarArg{i},\timeVarArg{i+1}) - \stateVar -
\tssfnk\fArgNoParam{\fnpropsymBE(\stateVar;\timeVarArg{i},\timeVarArg{i+1})}{\timeVarArg{i+1}}=0,\
i\in\natZero{\nptfn-1}.$} We define 
\begin{align} \label{eq:propagatorExact}
\begin{split}
\stateNP&:\timeVar\mapsto\fnprop{\timeVar}{0}{\initstateNP}, \quad \ourReReading{\timeVar \in \timeInstanceSet}
\end{split}
\end{align} 
as the associated numerical solution with 
$\stateNP\in(\setOfTimeDepFunctions)^\ndof$, where
$\setOfTimeDepFunctions$ denotes the set of functions from 
$\timeInstanceSet$
to $\RR{}$. Note that Eqs.~\eqref{eq:fineIntermediate} and \eqref{eq:propagatorExact} 
imply $\stateNP(\tptfn{j})
=\fnprop{\tptfn{j}}{\tptfn{i}}{\stateNP(\tptfn{i})},\ 
0\leq i\leq j\leq\nTimestepsFine
$. 
\ourReReading{It is this time-discrete solution $\stateNP$, which we want to approximate with the time-parallel procedure.}

\begin{figure}[bt] 
\centering 
\includegraphics[width=0.55\textwidth]{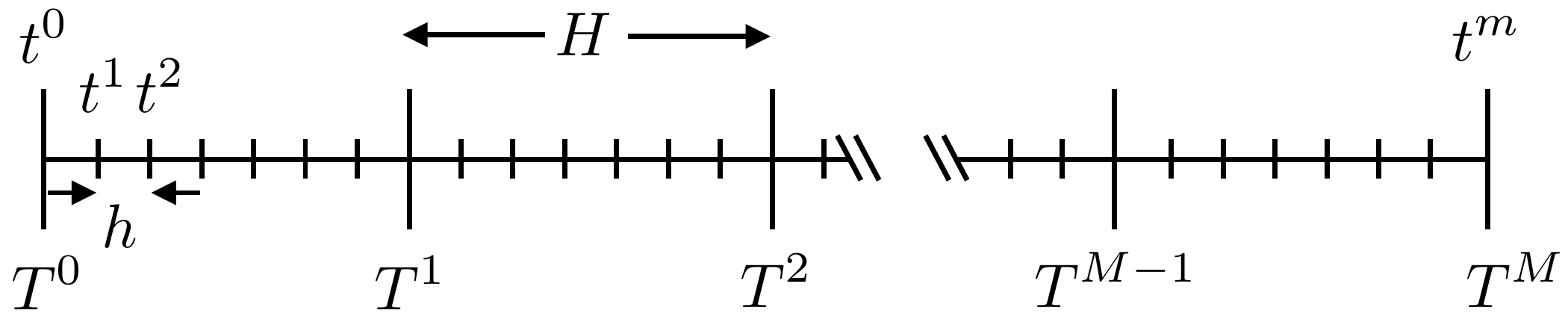} 
\caption{\reviewerB{Coarse and fine time discretizations.}} 
\label{fig:timediscretizations} 
\end{figure}

Analogously, we consider a \textit{coarse} time discretization characterized by
\reviewerB{(uniform) time step} $\tsscrsk\reviewerB{=\finalT/\nptcrs}$ and time instances
\reviewerB{$\tptcrs{\timestepdummy } = \timestepdummy \tsscrsk$, $\timestepdummy \innatZero{\nptcrs}$}, where
$\nptcrs\in\natNo$ denotes the number of coarse time instances beyond the
initial time $\tptcrs{0}=\tptfn{0}=0$ such that the final time corresponds to
$\tptcrs{\nptcrs}=\tptfn{\nTimestepsFine}=\finalT$
(see Figure
\ref{fig:timediscretizations}). 
We denote the set of time instances
associated with the coarse discretization as
$\timeInstanceSetCoarse\defeq\toset{\tptcrs{\timestepdummy }}{\timestepdummy }{0}{\nIntervals}$.
\reviewerB{Further, we assume that the coarse time step is an integral multiple of the
fine time step, i.e., $\tsscrsk =
\nfinepercoarse \tssfn$ with $\nfinepercoarse\in\natNo$. This implies that
the coarse discretization is nested within the fine
discretization $\timeInstanceSetCoarse\subseteq\timeInstanceSet$ 
\reviewerB{such that $\tptcrs{\timestepdummy } = \tptfn{\crstofn{\timestepdummy }}$,
$\timestepdummy \innatZero{\nIntervals}$ and $\nptfn=\nfinepercoarse\nptcrs$.}}
We define the set of fine time instances associated with the
$\timestepdummy$th coarse time interval
as 
\reviewerB{$\timeInstanceSetInterval{\timestepdummy }\defeq\timeInstanceSet\cap
\left[\tptcrs{\timestepdummy },\tptcrs{\timestepdummy +1}\right]
=\toset{\tptfn{i}}{i}{\nfinepercoarse
\timestepdummy }{\nfinepercoarse
(\timestepdummy +1)}
$, $\timestepdummy \innatZero{\nIntervals-1}$.}
%

	Denoting by $\apxsolFOM{\tixcrs}{\itvar}$ the approximation to
$\stateNP(\tptcrs{\tixcrs})$ 
at parareal iteration $\itvar$,
the parareal method 
first computes an initial guess
$\apxsolFOM{\timestepdummy }{0}$, $\timestepdummy \innatZero{\nptcrs-1} $ with
$\apxsolFOM{0}{0}=\initstateNP$ (typically via 
$\apxsolFOM{\timestepdummy +1}{0} = \crsprop{\tptcrs{\tixcrs +
1}}{\tptcrs{\tixcrs}}{\apxsolFOM{\tixcrs}{0}}$), and subsequently executes the following
iterations 
\begin{align}\label{eq:parareal}
\quad \apxsolFOM{\tixcrs +1}{\itvar +1} = \crsprop{\tptcrs{\tixcrs
+ 1}}{\tptcrs{\tixcrs}}{\apxsolFOM{\tixcrs}{\itvar +1}} +  \fnprop{\tptcrs{\tixcrs +
1}}{\tptcrs{\tixcrs}}{\apxsolFOM{\tixcrs}{\itvar}} -
\crsprop{\tptcrs{\tixcrs +
1}}{\tptcrs{\tixcrs}}{\apxsolFOM{\tixcrs}{\itvar}},\ \pararealit \innatZeroseq{\pararealItConverge},\
\timestepdummy =\pararealit ,\ldots,\nptcrs-1,
\end{align}
 where 
 \ourReReading{$ \apxsolFOM{\pararealit}{\pararealit +1}= \apxsolFOM{\pararealit }{\pararealit}$
and} 
 $\pararealItConverge$ is determined by a termination criterion that
	is satisfied when 
	the solution discontinuities at coarse time instances become sufficiently small.
	Here, $\crspropsym:\RR{\ndof}\times \timeInstanceSetCoarse\times
\timeInstanceSetCoarse\rightarrow\RR{\ndof}$ with 
$(\stateVar;\tptcrs{i},\tptcrs{j})\mapsto\crsprop{\tptcrs{j}}{\tptcrs{i}}{\stateVar
}$ denotes a `coarse propagator' that propagates a
state $\stateVar$ defined at (coarse) time instance $\tptcrs{i}$ to 
time instance $\tptcrs{j}$ with 
$j > i$. 
In essence, the parareal method alternates
serial (inexpensive) \textit{coarse} propagation with parallel (expensive) \textit{fine} propagation;
the expectation is that parallelizing the fine propagation can realize
wall-time performance improvements. Algorithm \ref{alg:parareal}---which
enables alternative initializations---reports the particular parareal
algorithm we consider in this work.

Critically, this method exhibits the `finite-termination property', which is
the result 
\begin{equation}\label{eq:finiteTermination}
\apxsolFOM{\timestepdummy }{\pararealit } = \stateNP(\tptcrs{\timestepdummy
}),\quad  \timestepdummy \leq \pararealit +1.
\end{equation}	
This states that the method will terminate in at most
$\pararealItConverge=\nTimeIntervals-1$ parareal iterations; realizing this
`worst-case scenario' implies that the parallelization over time provided no
gain over numerically solving Eq.~\eqref{eq:ODE} using the fine propagator in
serial. 

\begin{algorithm}[htbp]
\caption{ \pararealAlgorithmModName }
\begin{algorithmic}[1]\label{alg:parareal}
\REQUIRE Fine propagator $\fnpropsym$, coarse propagator $\crspropsym$,
initialization algorithm \initializeAlgorithmName,
initial condition $\apxsolFOM{0}{0}$, termination tolerance $\jumpTolerance$
\ENSURE Number of parareal iterations $\pararealItConverge\leftarrow \pararealit $, converged solution
$(\apxsolFOM{1}{0},\ldots,\apxsolFOM{\pararealItConverge+1}{\pararealItConverge},\ldots,\apxsolFOM{\nTimeIntervals}{\pararealItConverge})$
\STATE\label{step:initializeOne} $\pararealit \leftarrow 0$
\STATE\label{step:initialize}
$(\apxsolFOM{1}{0},\ldots,\apxsolFOM{\nTimeIntervals}{0})=
\initializeAlgorithmName(\apxsolFOM{0}{0})$
\FOR[parallel fine propagation]{$\timestepdummy =0
,\ldots,\nIntervals-1$}\label{step:beginFinePropInit}
\STATE\label{step:finePropInit} $
\finesolFOM{\timestepdummy+1}{0}  =
\fnprop{\tptcrs{\timestepdummy +1}}{\tptcrs{\timestepdummy }}{\apxsolFOM{\timestepdummy }{0}}$
\ENDFOR\label{step:initializeLast}
\WHILE{$\max_{\timestepdummy \in\{\pararealit  + 1,\ldots,\nIntervals-1\}}\|
	\finesolFOM{\timestepdummy}{\pararealit}
	-\apxsolFOM{\timestepdummy }{\pararealit  }\|/\|\finesolFOM{\timestepdummy}{\pararealit}
\|>
\jumpTolerance$}\label{step:termination}
\IF[initial-seed coarse propagation]{$\pararealit =0$}
\FOR[parallel coarse propagation]{$\timestepdummy =1 ,\ldots,\nIntervals-1$}
\STATE\label{step:initial_seed_prop} $\coarsesolFOM{\timestepdummy +1}{0} =
\crsprop{\tptcrs{\timestepdummy +1}}{\tptcrs{\timestepdummy }}{\apxsolFOM{\timestepdummy }{0}}$
\ENDFOR
\ENDIF
\STATE \ourReReading{
$ \apxsolFOM{\pararealit}{\pararealit +1}= \apxsolFOM{\pararealit }{\pararealit}$
}
\FOR[serial coarse propagation and correction]{$\timestepdummy =\pararealit ,\ldots,\nIntervals-1$}
\STATE\label{step:serialcoarse} $\coarsesolFOM{\timestepdummy +1}{\pararealit
+1} =
\crsprop{\tptcrs{\timestepdummy +1}}{\tptcrs{\timestepdummy
}}{\apxsolFOM{\timestepdummy }{\pararealit +1}}$
\STATE\label{step:overwrite} $\apxsolFOM{\timestepdummy +1}{\pararealit +1 } =
\coarsesolFOM{\timestepdummy +1}{\pararealit +1}+
\finesolFOM{\timestepdummy +1}{\pararealit  }
 -
\coarsesolFOM{\timestepdummy +1}{\pararealit  }$
\ENDFOR
\FOR[parallel fine propagation]{$\timestepdummy =\pararealit  ,\ldots,\nIntervals-1$}\label{step:beginFineProp}
\STATE\label{step:fineProp1} $\finesolFOM{\timestepdummy +1}{\pararealit+1} =
\fnprop{\tptcrs{\timestepdummy +1}}{\tptcrs{\timestepdummy
}}{\apxsolFOM{\timestepdummy }{\pararealit +1 }}$
\ENDFOR\label{step:endFineProp}
\STATE $\pararealit \leftarrow \pararealit +1$
\ENDWHILE
\reviewerB{
\FOR{$\timestepdummy =\pararealit  ,\ldots,\nIntervals-1$}
\STATE\label{step:fineProp2} $
\apxsolFOM{\timestepdummy +1}{\pararealit  } 
\leftarrow
\finesolFOM{\timestepdummy +1}{\pararealit} $
\ENDFOR}
\end{algorithmic}
\end{algorithm}


\section{Data-driven time parallelism}\label{sec:proposed}

The objective of this work is to devise inputs to Algorithm \ref{alg:parareal}
that leverage the availability of data that inform the \textit{time evolution} of the
state. Our two primary points of focus are (1) to devise an initialization
method that yields an accurate initial guess, and (2) to develop a coarse
propagator that is fast, accurate, and stable. In particular, we aim to improve upon
the performance of existing techniques, which generally employ coarse
propagators and initialization techniques that do not
exploit time-evolution data that may be available. 

Our critical assumption is that we 
have access to \textit{time-evolution bases}
$\timebasisj\in\orthogonalMatrix{\nTimestepsFine}{\dimBasisj}$,
$j\innat{\ndof}$
with $\dimBasisj\leq \nTimestepsFine$ that describe the \textit{time evolution} of the $j$th state
$\stateNPEntry{j}$.
Here, $\orthogonalMatrix{\ndof}{\timestepdummy } \ourReReading{\subset} \R^{\ndof \times \timestepdummy}$ 
denotes the Stiefel manifold,
\ourReReading{i.e., 
the set of all real-valued $\ndof\times \timestepdummy$ matrices with orthonormal columns.}
Subsequent sections will describe how these bases can be
computed in the case of parameterized ODEs (Section
\ref{sec:FOM}) and projection-based reduced-order models (Section
\ref{sec:ROM}); for now, we simply assume that these bases are available 
\ourReReading{and for ease of notation all have identical dimension $\dimBasisj$.}

\subsection{Global forecasting}\label{sec:forecasting}
We begin by summarizing the data-driven forecasting method proposed in
Ref.\
\cite{Carlberg_carlberg2015decreasing}.
Given bases $\timebasisj$, $j\innat{\ndof}$ and a time
instance $\currentTimeInstance\innat{\totalTS}$, the forecasting approach
approximates the time evolution of state variable $\stateNPEntry{\tevix}$ via
gappy POD using the basis $\timebasisj$ and the value of
$\stateNPEntry{\tevix}$ at the most recent $\memory$ time instances. Here,
$\memory\in\natNo$ with $\memory\geq \dimBasisj$ denotes the `memory', which will
be considered a global variable in this manuscript. First,
the method computes \ourReReading{the gappy POD approximation} $\fccoeffArg{\tptfn{\currentTimeInstance}}{\stateNPEntry{\tevix}}$,
defined as
\begin{align}\label{eq:genSolutionGlobalForecast}
	\fccoeff\ourReReading{(\stateContinuousDummyEntryNo;\timeVarArg{\ourReReading{i}
	})=} \underset{\ourReReading{\fcarg \in \range{\timebasisj}}}{\arg\min} \|
\sampleMat{\ourReReading{i} }{\prevtpts} \fcarg -
\sampleMat{\ourReReading{i} }{\prevtpts} \unrollfunc{\stateContinuousDummyEntryNo} \|_2
=\ourReReading{\timebasisj}[\sampleMat{\ourReReading{i} }{\prevtpts}
\timebasisj]^+\sampleMat{\ourReReading{i} }{\prevtpts}
\unrollfunc{\stateContinuousDummyEntryNo}\ourReReading{,
\ourReReading{\quad i \innatZero{\nptfn-\memory}},\,y \in \setOfTimeDepFunctions,}
\end{align}
where the superscript + denotes the Moore--Penrose pseudoinverse,
\ourReReading{$\range{{\bm A}}$ denotes the range of the matrix
${\bm A}$,} and $\fccoeff:\setOfTimeDepFunctions\times
\timeInstanceSet\rightarrow\RR{\ourReReading{\nptfn}}$. Here,
the \textit{sampling matrix} 
\ourReReading{$\sampleMat{i}{\memory} \defeq
\bmat{\unitvec{i+1} \ \cdots \ \unitvec{i+\memory}
}^{T} \in \{0,1\}^{\memory\times\nptfn}$} extracts entries $i+1$ through
$i+\memory$ of a given vector 
\ourReReading{and $\unitvec{i}\in \{0,1\}^{\nptfn}$ denotes the $i$th canonical 
unit vector.} 
Further,
$\unrollfuncNo:\setOfTimeDepFunctions\rightarrow\RR{\nTimestepsFine}$
centers and `unrolls' \ourReReading{a time-dependent variable} according to the time discretization as
\begin{equation}
\unrollfuncNo:{\stateContinuousDummyEntryNo} \mapsto \left[
\stateContinuousDummyEntryNo(\tptfn{1}) -
\stateContinuousDummyEntryNo(\tptfn{0})\
\cdots\
\stateContinuousDummyEntryNo(\tptfn{\nTimestepsFine})- \stateContinuousDummyEntryNo(\tptfn{0})
\right]^T.
\end{equation}

\begin{figure}[htbp] 
\centering 
\subfigure[\reviewerB{Global forecast.
 Here,
 the time evolution basis vectors (i.e., columns of $\timebasisj$) are denoted by thin
colored
lines, the state entry $\stateNPEntry{j}$ is denoted by a thick black line, the sampled
state $\sampleMat{i}{}\unrollfunc{\stateNPEntry{j}} +
\stateNPEntry{j}(\tptfn{0})$ is denoted by +
markers, the forecast
$\forecastFunctionj{\cdot}{\timeVarArg{i}}{\stateNPEntry{j}}$ 
with $i =5$ is plotted as dashed line and 
the forecast
$\forecastFunctionj{\timeVarArg{j}}{\timeVarArg{i}}{\stateNPEntry{j}}$ 
for $j =19$ is denoted by a $\circ$ marker. 
}]{
\includegraphics[width=0.45\textwidth]{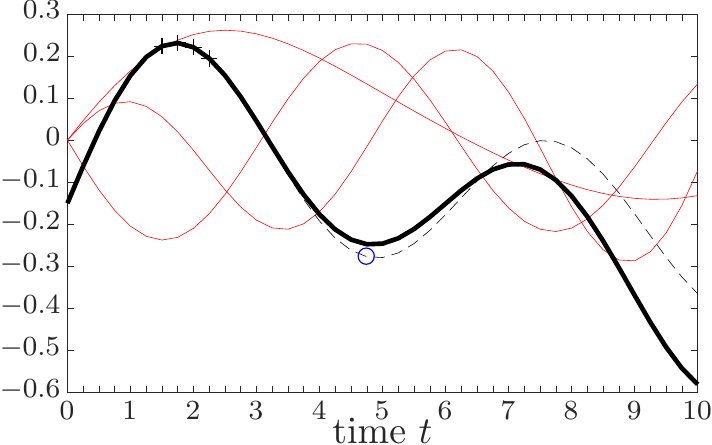} 
\label{fig:globalForecastIllustration} 
}$\,\,\,\,$
\subfigure[\reviewerB{Local forecast with $\nptcrs=3$ coarse time intervals. 
	Here, the time evolution basis vectors (i.e., columns of $\timebasisjm$) are denoted by
 thin
colored
lines, the state entry $\stateNPEntry{j}$ is denoted by a thick black line, the sampled
state 
$\sampleMat{0}{}\unrollmfunc{\stateNPEntry{j}}+ 
\stateNPEntry{j}(\tptcrs{\timestepdummy })$ is denoted by +
markers, and the forecast
$\forecastFunctionjm{\tptcrs{\timestepdummy +1}}{\tptcrs{\timestepdummy }}{\stateNPEntry{j}}$
is
denoted by $\circ$ markers. 
}]{
\includegraphics[width=0.45\textwidth]{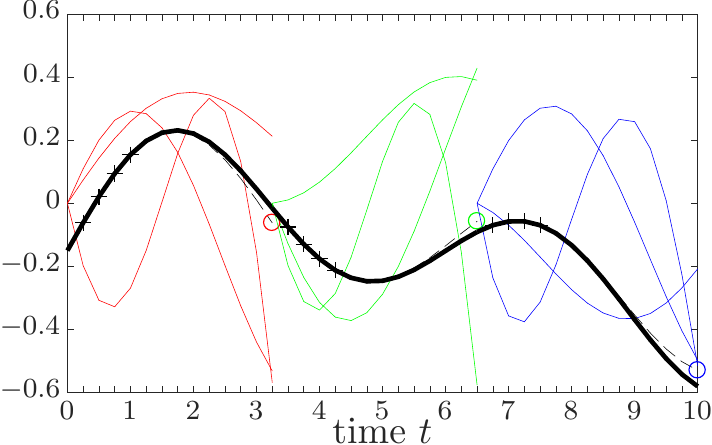} 
\label{fig:localForecastIllustration} 
}
\caption{\reviewerB{Illustration of global and local forecasting with memory $\memory=4$,
time-evolution basis dimension $\dimBasisj=3$, final time $\tintend = 10$,
and time step $\tssfnk=0.25$.}}
\end{figure}

Then, the forecast 
at a given time instance $\timestepdummy $, which aims to approximate the value
${\stateNPEntry{j}}(\tptfn{\timestepdummy })$, is set to $
\forecastFunctionj{\tptfn{\timestepdummy }}{\tptfn{\currentTimeInstance}}{\stateNPEntry{j}}
$, where we have
defined the function that forecasts the time-dependent variable $\stateContinuousDummyEntryNo$
to
time $\timeVarArg{k}$ using its value at times 
$\timeVarArg{\ourReReading{i} +\ourReReading{\ell}}$, $\ourReReading{\ell}\innat{\memory}$ as
 \begin{align} \label{eq:globalForecast}
\forecastFunctionjNo:\forecastFunctionjArgs{\timeVarArg{k}}{\timeVarArg{\ourReReading{i} }}{\stateContinuousDummyEntryNo}
\mapsto\stateContinuousDummyEntryNo(0) +
\unitvec{k}^T\fccoeffArg{\timeVarArg{\ourReReading{i}}}{\stateContinuousDummyEntryNo}\ourReReading{,}
\ourReReading{\quad i \innatZero{\nptfn-\memory},\, k \innatZero{\nptfn},\,y \in \setOfTimeDepFunctions,}
  \end{align} 
	with $\forecastFunctionjNo:\setOfTimeDepFunctions\times
	\timeInstanceSet\times \timeInstanceSet\rightarrow \RR{}$. Figure
	\ref{fig:globalForecastIllustration} illustrates the global-forecasting
	method graphically\reviewerA{, and Algorithm \ref{alg:global_fore} provides an algorithmic
	description of the method such that
$$\forecastFunctionj{\timeVarArg{k}}{\timeVarArg{i
}}{\stateContinuousDummyEntryNo} =
\globalForeAlg(\timebasisj,\stateContinuousDummyEntryNo(0),\{\stateContinuousDummyEntryNo(\timeVarArg{\ell})\}_{\ell=i+1}^{i+\memory},i,k) .$$}
\begin{algorithm}[tb]
	\caption{\reviewerA{\globalForeAlg\ (algorithmic description of the global
	forecast \eqref{eq:globalForecast})}}
\begin{algorithmic}[1]\label{alg:global_fore}
	\reviewerA{
\REQUIRE time-evolution basis $\timebasisj\in\orthogonalMatrix{\nTimestepsFine}{\dimBasisj}$, initial state 
$\stateContinuousDummyEntryNo(0)\in\RR{}$, 
sampled state
$\{\stateContinuousDummyEntryNo(\timeVarArg{\ell})\}_{\ell=i+1}^{i+\memory}\subset\RR{}$,
initial sampling time index $i \in\natZero{\nptfn-\memory}$, 
forecast time index
$k \in\natZero{\nptfn}$
\ENSURE global forecast $\forecastFunctionj{\timeVarArg{k
}}{\timeVarArg{i }}{\stateContinuousDummyEntryNo}\in\RR{}$
\STATE Solve gappy POD linear least-squares problem \eqref{eq:genSolutionGlobalForecast}
for $\fccoeff(\stateContinuousDummyEntryNo;\timeVarArg{i
})\in\RR{\nptfn}$, noting that $\sampleMat{i }{\prevtpts}
	\unrollfunc{\stateContinuousDummyEntryNo}=\left[
\stateContinuousDummyEntryNo(\timeVarArg{i + 1}) -
\stateContinuousDummyEntryNo(0), \cdots, 
\stateContinuousDummyEntryNo(\timeVarArg{i + \memory}) -
\stateContinuousDummyEntryNo(0)
\right]^T$. Note that the gappy POD approximation over all time is $\unrollfunc{\forecastFunction{\cdot}{\timeVarArg{i
 }}{\stateContinuousDummyEntryNo}} = \ones{\nptfn}
 \stateContinuousDummyEntryNo(0) + \fccoeff(\stateContinuousDummyEntryNo;\timeVarArg{i
	})$, where
 $\ones{m}$ denotes an $m$-vector of ones.
\STATE Extract forecast at time instance $\timeVarArg{k }$ via Eq.~\eqref{eq:globalForecast} as
 ${\forecastFunctionj{\timeVarArg{k }}{\timeVarArg{i
 }}{\stateContinuousDummyEntryNo}} = 
 \unitvec{k}^T\unrollfunc{\forecastFunctionj{\cdot}{\timeVarArg{i
 }}{\stateContinuousDummyEntryNo}}
 = 
 \stateContinuousDummyEntryNo(0) + \unitvec{k}^T\fccoeff(\stateContinuousDummyEntryNo;\timeVarArg{i
	}).$ 
	}
\end{algorithmic}
\end{algorithm}

The approach proposed in Ref.\
\cite{Carlberg_carlberg2015decreasing} employed
the forecast 
$
\forecastFunctionj{\tptfn{\timestepdummy }}{\tptfn{\currentTimeInstance}}{\stateNPEntry{j}}
$,
$j\innat{\ndof}$ as an
initial guess for the Newton solver at time $\tptfn{\timestepdummy }$ for $\timestepdummy  >
\currentTimeInstance+\memory$ obtained
after discretizing the ODE associated with a ROM using a linear multistep scheme.\footnote{The
method was also generalized to handle Runge--Kutta schemes and second-order ODEs;
the only difference in these cases is that the forecast is constructed for the
\textit{unknown} variable computed at each time step, which can correspond to
the velocity or acceleration depending on the ODE and time integrator.}
Instead, this work considers employing this forecasting strategy to define
both the \textit{initialization} and \textit{coarse propagator} as inputs to
parareal Algorithm \ref{alg:parareal}. We now propose a \textit{local}
variant of this global forecasting method that operates within a single coarse
time interval.
   
\subsection{Local forecasting}\label{sec:localforecast}
\begin{algorithm}[tb]
\caption{\localalgNo}
\begin{algorithmic}[1]\label{alg:localBasis}
	\REQUIRE time-evolution basis $\ourReReading{\timebasisj}\in\orthogonalMatrix{\nTimestepsFine}{\dimBasisj}$,
coarse-time-interval index $\timestepdummy \in\natZero{\nIntervals-1}$, energy criterion
$\energyCrit\in\left[0,1\right]$
\ENSURE local time-evolution basis 
$\ourReReading{\timebasisjm}\in\orthogonalMatrix{\nptfnk}{\dimBasism}$ 
\STATE $\ourReReading{\timebasisjm}\leftarrow\sampleMatLocal{\timestepdummy }\ourReReading{\timebasisj} $ \COMMENT{Extract
values on $\timestepdummy$th coarse time interval}
\STATE $\ourReReading{(\U,\Sig,\V)} = \texttt{thin\_SVD}(\ourReReading{\timebasisjm})$ \COMMENT{Compute (thin) singular value decomposition}
\STATE\label{step:truncateSVD} $\ourReReading{\timebasisjm}\leftarrow \vectomat{\uvec}{\dimBasism}$, where 
$\dimBasism = \min_{i\in\ourReReading{\energyCritSet}(\energyCrit)} i$, 
$\ourReReading{\energyCritSet{(\energyCrit)}} \defeq \{i\ |\
	\sum_{\ell=1}^i\singularValueArg{\ell}/\sum_{k
	=1}^\dimBasis\singularValueArg{k}\geq
\energyCrit\}$, \\ $\ourReReading{\Sig}=\mathrm{diag}(\singularValueArg{1},\ldots,
\singularValueArg{\dimBasis})$, 
$\ourReReading{\U} = \vectomat{\ourReReading{\uvec}}{\dimBasis}$. \COMMENT{Truncate}
\end{algorithmic}
\end{algorithm}
The proposed \textit{local} forecasting approach relies on 
local time-evolution bases
$\timebasisjm\in\orthogonalMatrix{\nptfnArg{\timestepdummy }}{\dimBasisjm}$,
$j\innat{\ndof}$ that inform the time evolution of
 the $j$th state
$\stateNPEntry{j}$ over time interval $[\tptcrs{\timestepdummy },\tptcrs{\timestepdummy +1}]$. Given a (global) time-evolution basis $\timebasisj$,
these local bases $\timebasisjm$, $\timestepdummy \innatZero{\nCoarseIntervals-1}$ can be computed via Algorithm \ref{alg:localBasis} as
$\timebasisjm = \localalg{\timebasisj}{\timestepdummy }{\energyCrit}$,
where $\energyCrit\in[0,1]$ defines a statistical `energy criterion' and we have defined 
$\sampleMatLocal{\timestepdummy }\defeq
[\unitvec{\crstofn{\timestepdummy }+1}\ \cdots\ \unitvec{\crstofn{\timestepdummy +1}}]^T - 
\ones{\nTimestepsIntervalArg{\timestepdummy }}[\unitvec{\crstofn{\timestepdummy }}]^T
\in\{0,1\}^{\nTimestepsIntervalArg{\timestepdummy }\times
\nTimestepsFine}$
as the matrix that samples entries associated with the $\timestepdummy $th
coarse time interval and subtracts the initial value on that time interval. Here,
$\ones{i}$ denotes an
$i$-vector of ones.
 Note that truncation in Step \ref{step:truncateSVD} of Algorithm \ref{alg:localBasis}
 ensures that the local basis $\timebasisjm$ will have full column rank.
 Using these local time-evolution bases (which have zero values at the
beginning of their respective time intervals), we can define the local
forecast using a similar construction to that of Section
\ref{sec:forecasting}. In particular, the linear least-squares problem for the
locally defined \ourReReading{gappy POD approximation}  becomes
\begin{align}\label{eq:defineGenCoords}
\begin{split}
	\fccoeffm\ourReReading{(\stateContinuousDummyEntryNo;\timeVarArg{i}) =}&
	\underset{\fcarg \in \ourReReading{\range{\timebasisjm}}}{\arg \min} \|
\sampleMat{i-\crstofn{\timestepdummy }}{\prevtpts}  \fcarg -
\sampleMat{i-\crstofn{\timestepdummy }}{\prevtpts} \unrollmfunc{\stateContinuousDummyEntryNo} \|_2
=
\timebasisjm[\sampleMat{i-\crstofn{\timestepdummy }}{\prevtpts}
\timebasisjm]^+\sampleMat{i-\crstofn{\timestepdummy }}{\prevtpts}
\unrollmfunc{\stateContinuousDummyEntryNo}
\end{split}
\end{align}
for $\timeVarArg{i+\ell}\in\timeInstanceSetInterval{\timestepdummy }\ourReReading{,}\ \forall
\ell\innat{\memory}$,
with $\fccoeffm:\setOfTimeDepFunctionsArg{\timestepdummy }\times
\timeInstanceSetInterval{\timestepdummy }\rightarrow\RR{\nfinepercoarse}$, where 
$\setOfTimeDepFunctionsArg{\timestepdummy }$, $\timestepdummy \innatZero{\nIntervals-1}$ denotes the set of functions from
$\timeInstanceSetInterval{\timestepdummy }$
to $\RR{}$.
Here,
the function
$\unrollmfuncNo:\setOfTimeDepFunctionsArg{\timestepdummy }\rightarrow\RR{\nTimestepsIntervalArg{\timestepdummy }}$ locally centers and unrolls 
\ourReReading{a time-dependent variable}
over the
$\timestepdummy $th time interval as
\begin{equation}\label{eq:localUnroll}
\unrollmfuncNo:{\stateContinuousDummyEntryNo} \mapsto
\left[
\stateContinuousDummyEntryNo(\tptfn{\crstofn{\timestepdummy }+1})- \stateContinuousDummyEntryNo(\Tm)\
\cdots\
\stateContinuousDummyEntryNo(\Tmp)- \stateContinuousDummyEntryNo(\Tm)
\right]^T.
\end{equation}
Note that if $\stateContinuousDummyEntryNo\in\setOfTimeDepFunctions$, then
$\unrollmfunc{\stateContinuousDummyEntryNo} =
\sampleMatLocal{\timestepdummy }\unrollfunc{\stateContinuousDummyEntryNo}$.
 Then, the function
 $\forecastFunctionjmNo$
 that forecasts a local time-dependent variable to time $\timeVarArg{\pararealit }$
 using the value of the variable at times $\timeVarArg{i+p}$,
 $p\innat{\memory}$ can be defined algebraically as
 \begin{align} \label{eq:localForecast}
 \begin{split}
\forecastFunctionjmNo&:\forecastFunctionjmArgs{\timeVarArg{\pararealit }}{\timeVarArg{i}}{\stateContinuousDummyEntryNo}
\mapsto\stateContinuousDummyEntryNo(\Tm) +
\unitvec{\pararealit -\crstofn{\timestepdummy
}}^T\fccoeffmArg{\timeVarArg{i}}{\stateContinuousDummyEntryNo},\quad\timeVarArg{\pararealit
},\timeVarArg{i+\ell}\in\timeInstanceSetInterval{\timestepdummy }\ourReReading{,}\
\forall \ell\innat{\memory}
 \end{split}
  \end{align} 
	with 
	$\forecastFunctionjmNo:\setOfTimeDepFunctionsArg{\timestepdummy }\times\timeInstanceSetInterval{\timestepdummy }\times\timeInstanceSetInterval{\timestepdummy }\rightarrow\RR{}$.
Figure
	\ref{fig:localForecastIllustration} illustrates the local-forecasting
	method graphically\reviewerA{, and Algorithm \ref{alg:local_fore} provides an algorithmic
	description of the method such that
$$\forecastFunctionjm{\timeVarArg{k}}{\timeVarArg{i
}}{\stateContinuousDummyEntryNo} =
\localForeAlg(\timestepdummy,\timebasisjm,\stateContinuousDummyEntryNo(\Tm),\{\stateContinuousDummyEntryNo(\timeVarArg{\ell})\}_{\ell=i+1}^{i+\memory},i,k).$$}
\begin{algorithm}[tb]
	\caption{\reviewerA{\localForeAlg\ (algorithmic description of the local
	forecast \eqref{eq:localForecast})}}
\begin{algorithmic}[1]\label{alg:local_fore}
	\reviewerA{
\REQUIRE time interval index $\timestepdummy\in\natZero{\nptcrs-1}$, time-evolution basis
$\timebasisjm\in\orthogonalMatrix{\nptfnArg{\timestepdummy }}{\dimBasisjm}$,
state at beginning of time interval
$\stateContinuousDummyEntryNo(\Tm)\in\RR{}$, 
sampled state
$\{\stateContinuousDummyEntryNo(\timeVarArg{\ell})\}_{\ell=i+1}^{i+\memory}\subset\RR{}$,
initial sampling time index $i
\in\{\nfinepercoarse\timestepdummy,\ldots,\nfinepercoarse(\timestepdummy+1)-\memory\}$, 
forecast time index
$k \in\{\nfinepercoarse\timestepdummy,\ldots,\nfinepercoarse(\timestepdummy+1)\}$
\ENSURE local forecast 
$\forecastFunctionjm{\timeVarArg{k }}{\timeVarArg{i}}{\stateContinuousDummyEntryNo}\in\RR{}$
\STATE Solve gappy POD linear least-squares problem \eqref{eq:defineGenCoords}
for $\fccoeffm(\stateContinuousDummyEntryNo;\timeVarArg{i})$, noting that 
$
\sampleMat{i-\crstofn{\timestepdummy }}{\prevtpts}
\unrollmfunc{\stateContinuousDummyEntryNo} = 
\left[
\stateContinuousDummyEntryNo(\timeVarArg{i + 1}) -
\stateContinuousDummyEntryNo(\Tm), \cdots,  
\stateContinuousDummyEntryNo(\timeVarArg{i + \memory}) -
\stateContinuousDummyEntryNo(\Tm)
\right]^T
$. Note that the gappy POD approximation over $\timestepdummy$th time interval
 is $\unrollmfunc{\forecastFunctionjm{\cdot}{\timeVarArg{i}}{\stateContinuousDummyEntryNo}}
 = \ones{\nfinepercoarse}
 \stateContinuousDummyEntryNo(\Tm) + \fccoeffm(\stateContinuousDummyEntryNo;\timeVarArg{i
	})$.
	\STATE Extract forecast at time instance $\timeVarArg{k }$ via Eq.~\eqref{eq:localForecast} as
 ${\forecastFunctionjm{\timeVarArg{k }}{\timeVarArg{i
 }}{\stateContinuousDummyEntryNo}} = 
 \unitvec{k -\crstofn{\timestepdummy
}}^T\unrollmfunc{\forecastFunctionjm{\cdot}{\timeVarArg{i}}{\stateContinuousDummyEntryNo}}
 = 
 \stateContinuousDummyEntryNo(\Tm) + \unitvec{k -\crstofn{\timestepdummy
}}^T\fccoeffm(\stateContinuousDummyEntryNo;\timeVarArg{i
	}).$
	}
\end{algorithmic}
\end{algorithm}

\subsection{Coarse propagator: local forecast}\label{sec:coarseprop}

	 We aim to employ the local forecasting approach to construct a
	data-driven coarse propagator to be used in the parareal Algorithm \ref{alg:parareal}. In particular, we propose 
       to construct a propagator
	that maps the state evaluated at the \textit{first $\memory$ \ourReReading{fine time instances}}
	of a given coarse time interval to an approximation of the state at the
	\textit{final time} of the coarse time
	interval. However, inspired by the multigrid interpretation of parareal, we
	acknowledge that the role of the coarse propagator is to reduce
	\textit{large-wavelength errors}; thus we allow the technique to apply this
	propagation only to a restriction of the
	state.\footnote{Numerical experiments highlight the
	importance of this (see Figure \ref{fig:nY_Test_plot}).} That is, we set the coarse propagation of the $j$th
	element of a \textit{restricted} time-dependent vector to be the mapping
\begin{align}\label{eq:coarsePropDefLocalPropFirst}
\begin{split}
\stateContinuousVar\mapsto\forecastFunctionjm{\Tmp}{\Tm}{\restrictjArg{\stateContinuousVar}}
:(\setOfTimeDepFunctionsArg{\timestepdummy })^\ndof\rightarrow\RR{},
\end{split}
\end{align}
where \reviewerB{$\restrictT\defeq\left[\restrictTEntry{1}\ \cdots\
\restrictTEntry{\nrestrict}\right]\in\RR{\ndof\times\nrestrict}$} with
$\nrestrict\ourReReading{\innat{\ndof}}$ denotes a \reviewerB{(linear)} restriction
operator with associated prolongation operator 
\reviewerB{$\prolongate\in\RR{\ndof\times\nrestrict}$.}
	Note that the time-evolution bases should therefore be constructed to capture the
	time-evolution of the \textit{restricted} time-dependent variable. Possible
	choices for the restriction operator include projection onto
	large-wavelength Fourier modes or onto a set of high-energy POD modes; 
	the latter choice is natural for reduced-order models and is explored in the
	numerical experiments.

Introducing a function that maps a vector at the beginning of a coarse time interval
to a function over the fine time discretization within that interval, i.e.,
$
\fineFillInLocalOne{\timestepdummy }:\stateVar\mapsto\fnprop{\cdot}{\tptcrs{\timestepdummy }}{\stateVar}$
with
$\fineFillInLocalOne{\timestepdummy }:\RR{\ndof}\rightarrow(\setOfTimeDepFunctionsArg{\timestepdummy })^\ndof$,
we define coarse propagation of the $j$th
	element of the restricted state on coarse time interval $\timestepdummy $ to be 
\begin{align}\label{eq:coarsePropDefLocalProp}
\begin{split}
\crspropforesymjm&:\stateVar\mapsto
\forecastFunctionjm{\Tmp}{\Tm}{\restrictjArg{{\fineFillInLocal{\timestepdummy }{\stateVar}}}}
,\quad
\timestepdummy \innatZero{\nptcrs-1}
\end{split}
\end{align}
with $\crspropforesymjm:\RR{\ndof}\rightarrow\RR{}$\ourReReading{, which can be expressed algebraically as
\begin{align}\label{eq:coarsePropDefLocalPropAlg}
\begin{split}
\crspropforesymjm&:\stateVar\mapsto
\restrictjArg{\stateVar} + \sum_{i=1}^\memory\forecastScalar{i}{j}
(
\restrictjArg{\fnprop{\tptfn{\crstofn{\timestepdummy }+i}}{\Tm}{\stateVar}}-\restrictjArg{\stateVar}
)
,\quad
\timestepdummy \innatZero{\nptcrs-1}
\end{split}
\end{align}
 with
$\forecastScalar{i}{j}\defeq\unitvec{\nTimestepsIntervalArg{\timestepdummy }}^T\timebasisjm\left[\sampleMat{
0}{\memory}\timebasisjm\right]^+\unitvec{i}\in\RR{}$.}
	We then propose employing a coarse propagator
	$\crspropsym\leftarrow\crspropforesymAll$ with 
	\begin{align}\label{eq:dataDrivenCoarse1}
	\begin{split}
	\crspropforesymAll&:(\stateVar;\tptcrs{\timestepdummy },\tptcrs{\timestepdummy  +
	1})\mapsto\prolongateArg{\left[
	\crspropforesym{1}{\timestepdummy }(\stateVar)\
	\cdots\
	\crspropforesym{\nrestrict}{\timestepdummy }(\stateVar)
	\right]^T}
	\end{split}
	\end{align}
 with $\crspropforesymAll:\RR{\ndof}\times \timeInstanceSetCoarse\times
	\timeInstanceSetCoarse\rightarrow\RR{\ndof}$,
which can be expressed
algebraically as
	\begin{equation}\label{eq:dataDrivenCoarse2}
	\crspropforesymAll:(\stateVar;\tptcrs{\tixcrs},\tptcrs{\tixcrs +
	1})\mapsto
	\prolongateArg{\restrictArg{\stateVar}} + \prolongateArg{\sum_{i=1}^\memory\linearOpForecast{i}\left[
\restrictArg{\fnprop{\tptfn{\crstofn{\timestepdummy }+i}}{\Tm}{\stateVar}}-\restrictArg{\stateVar}
	\right]}.
	\end{equation}
Here, we have defined
$\linearOpForecast{i}\defeq\diag(\forecastScalar{i}{1},\ldots,\forecastScalar{i}{\nrestrict})\in\RR{\nrestrict\times\nrestrict}$.

\subsection{Initialization: local and global forecasts}\label{sec:initialSeed}
 Initialization in Step \ref{step:initialize} of Algorithm
\ref{alg:parareal} is typically executed by sequentially applying the coarse
propagator, i.e., 
\begin{equation}\label{eqn:standard_initialization}
\apxsolFOM{\timestepdummy +1}{0} = \crsprop{\tptcrs{\timestepdummy +1}}{\tptcrs{\timestepdummy }}{\apxsolFOM{\timestepdummy }{0}},\quad
\timestepdummy \innatZeroseq{\nTimeIntervals-1}. 
\end{equation}
 This approach could be applied with the proposed local-forecasting coarse
propagator 
$\crspropsym\leftarrow\crspropforesymAll$. However, we can also consider an alternative
initialization that is both computationally less expensive and more stable (as
will be further discussed in Remark \ref{remark:initialSeedStable}). 

In particular, we 
consider performing initialization by forecasting the state from
\textit{the first $\memory$ time steps} of the first time interval to
\textit{all coarse time instances} using the
global time-evolution bases $\timebasisj$. That is, we can perform
initialization via \textit{global forecasting} as
\begin{equation} \label{eq:globalForecastInitialize}
\apxsolFOM{\timestepdummy +1}{0} = \crspropforeGlobal{\tptcrs{\timestepdummy +1}}{\apxsolFOM{0}{0}},\quad
\timestepdummy \innatZeroseq{\nTimeIntervals-1},
\end{equation} 
where we have defined
\begin{align} 
\begin{split}
\crspropforesymAllGlobal
&:(\stateVar;\tptcrs{\tixcrs})\mapsto\prolongateArg{[\forecastFunctionFour{\tptcrs{\timestepdummy }}{\tptcrs{0}}{\restrictEntryArg{1}{\fineFillInOne{\stateVar}}}{1}\
\cdots\
\forecastFunctionFour{\tptcrs{\timestepdummy }}{\tptcrs{0}}{\restrictEntryArg{\nrestrict}{\fineFillInOne{\stateVar}}}{\nrestrict}]^T}
\end{split}
\end{align} 
\ourReReading{with} $\crspropforesymAllGlobal:\RR{\ndof}\times 
	\timeInstanceSetCoarse\rightarrow\RR{\ndof}$\ourReReading{, and 
	$
\fineFillInOne{\stateVar}:\stateVar\mapsto\fnprop{\cdot}{\tptfn{0}}{\stateVar}$
with
$\fineFillInOne{\stateVar}:\RR{\ndof}\rightarrow(\setOfTimeDepFunctions)^\ndof$}, which can be expressed algebraically as
 \begin{equation}
\crspropforesymAllGlobal
:(\stateVar;\tptcrs{\tixcrs})\mapsto
\prolongateArg{\restrictArg{\stateVar}} + \prolongateArg{\sum_{i=1}^\memory\linearOpForecastGlobal{i}\left[
\restrictArg{\fnprop{\tptfn{i}}{\tptcrs{\timestepdummy }}{\stateVar}}-\restrictArg{\stateVar}
	\right]}.
\label{eq:globalForecastInitializeAlgebraic}
  \end{equation} 
Here, we have defined
$\linearOpForecastGlobal{i}\defeq\diag(\forecastScalarGlobal{i}{1},\ldots,\forecastScalarGlobal{i}{\nrestrict})\in\RR{\nrestrict\times\nrestrict}$
and
$\forecastScalarGlobal{i}{j}\defeq\unitvec{\crstofn{\timestepdummy
}}^T\timebasisj\left[\sampleMat{
0}{\memory}\timebasisj\right]^+\unitvec{i}\in\RR{}$.

%
%
\section{Analysis}\label{sec:analysis}

We now analyze the proposed data-driven time-parallel methodology to derive
insight into the coarse-propagator error (Section \ref{sec:error}), the
method's theoretical speedup (Section \ref{sec:speedup}), the
method's stability (Section \ref{sec:stability})
\ourReReading{ and convergence aspects (Section \ref{sec:convergence})}.  All norms in this section
refer to the Euclidean norm unless otherwise specified. Appendix
\ref{sec:proofs} contains all proofs.
\subsection{Coarse-propagator error analysis}\label{sec:error}
We first analyze the error of the coarse propagator with respect to the fine
propagator.

\subsubsection{General case}\label{sec:errorGeneral}
We 
introduce the following assumptions:
\begin{Assumption}[series=assumption]
\item\label{ass:prolongOrtho} The restriction and prolongation operators have
counterparts
\reviewerB{$\restrictOrthT\in\RR{\ndof\times(\ndof-\nrestrict)}$}
and 
\reviewerB{$\prolongateOrth\in\RR{\ndof\times(\ndof-\nrestrict)}$}, respectively, 
that satisfy $\stateVar =
\prolongateArg{\restrictArg{\stateVar}} +
\prolongateOrthArg{\restrictOrthArg{\stateVar}}$, $\forall \stateVar\in\RR{\ndof}$.
\item \label{ass:lipschitzProlong} The prolongation 
operators are bounded
\ourReReading{by constants $\lipschitzProlong,\lipschitzProlongOrth \in \R$ }, i.e.,
$\|\prolongateArg{\stateVar}\|\leq\lipschitzProlong\|\stateVar\|,\ \forall\stateVar\in\RR{\nrestrict}$
and
$\|\prolongateOrthArg{\stateVar}\|\leq\lipschitzProlongOrth\|\stateVar\|,\ \forall\stateVar\in\RR{\ndof-\nrestrict}$.
\end{Assumption}

\begin{theorem}\label{lem:coarseError}
If Assumptions \ref{ass:prolongOrtho} and \ref{ass:lipschitzProlong} hold, then
\begin{align} \label{eq:fineCoarsePropError}
\begin{split}
\|\fnprop{\Tmp}{\Tm}{\stateVar} -
\crspropfore{\Tmp}{\Tm}{\stateVar}
\|
\leq
&\lipschitzProlongOrth\|\restrictOrthArg{\fnprop{\Tmp}{\Tm}{\stateVar}}\| \\
& 
+\lipschitzProlong\sum_{j=1}^{\nrestrict}
\normedQuantityGappyLarger{j}{\timestepdummy }\Big\|
(\identityArg{\nptfnArg{\timestepdummy }} - \timebasisjm[\timebasisjm]^T)\left[
\begin{array}{c}
\restrictjArg{\fnprop{\tptfn{\crstofn{\timestepdummy }+1}}{\Tm}{\stateVar}}
-\restrictjArg{\stateVar}\\
\vdots\\
\restrictjArg{\fnprop{\Tmp}{\Tm}{\stateVar}}
-\restrictjArg{\stateVar}
\end{array}
\right]
\Big\|.
\end{split}
\end{align} 
where $\normedQuantityGappyLarger{j}{\timestepdummy }\defeq 1/\minSingular{\sampleMat{
0}{\memory}\timebasisjm}\geq 1$.\footnote{Note that
$\normedQuantityGappyLarger{j}{\timestepdummy }\geq 1$ because appending a row to a matrix
cannot decrease its minimum singular value, and $\minSingular{\timebasisjm}=1$
because $\timebasisjm\in\orthogonalMatrix{\nptfnArg{\timestepdummy }}{\dimBasisjm}$.}
\end{theorem}
\begin{remark}[Interpolation v.~oversampling]\label{rem:interpOversampling}
As the memory $\memory$ increases, the stability constants
$\normedQuantityGappyLarger{j}{\timestepdummy }$ in inequality
\eqref{eq:fineCoarsePropError} monotonically decrease. This occurs
because increasing the memory has the effect of appending a row to the matrix
$\sampleMat{ 0}{\memory}\timebasisjm$, which cannot
decrease its minimum singular value. This highlights the stabilizing effect of
employing a least-squares approach (i.e., gappy POD) as opposed to an
interpolation approach (i.e., EIM/DEIM) in the forecast: oversampling can
reduce a bound for the error between the fine and coarse propagators.
\end{remark}

\begin{remark}[Restriction tradeoff]\label{rem:restrictionTradeoff}
Increasing the dimension of the restriction operator (i.e., the number of
variables included in the forecast $\nrestrict$)
\reviewerB{decreases} the first term in bound
\eqref{eq:fineCoarsePropError}. However, doing so also increases
the second term, as the number of terms in the summation increases. This
latter effect is exacerbated when the time evolution of higher-index
solution components (i.e., $\restrictjArg{\stateVar}$ for large $j$) is not
well captured by the associated time-evolution bases
(i.e., $\timebasisjm$ for large $j$); this can occur, for example, if higher-index solution
components associate with high-frequency solution modes, as is the case when
the restriction operator associates with a projection onto a low-frequency
Fourier or POD (see Section
\ref{sec:POD}) basis.
These two effects comprise the tradeoff that should be considered when
selecting the dimension of the restriction operator $\nrestrict$ in practice.
\end{remark}

\subsubsection{Ideal case}\label{sec:idealError} 
We now show that the coarse propagator is exact (i.e., incurs no error with
respect to the fine propagator) under the following `ideal conditions':
\begin{Assumption}[resume=assumption]
\item\label{ass:subspace} The time evolution of the restricted state is an
element of the subspace spanned by the time-evolution basis (i.e.,
$\unrollfunc{\restrictjArg{\stateNP}}\in\range{\timebasisj}$,
$j\innat{\nrestrict}$).
\item \label{ass:notruncate} The local bases are constructed with no truncation (i.e.,
$\energyCrit = 1.0$ in Algorithm \ref{alg:localBasis}).
\item\label{ass:isomorphic} The original and restricted state spaces are isomorphic
(i.e., $\prolongate\restrict= \restrict\prolongate=
\identityArg{\ndof}$ with $\nrestrict = \ndof$).
\end{Assumption}
\begin{lemma}[Local-subspace condition]\label{lem:localSubspace}
If Assumptions \ref{ass:subspace} and \ref{ass:notruncate} hold, then 
$\unrollmfunc{\restrictjArg{\stateNP}}\in\range{\timebasisjm},\;
j\innat{\ndof},\ \timestepdummy \innatZero{\nTimeIntervals-1}.$
\end{lemma}

\begin{theorem}[\ourReReading{Exact coarse propagator}]\label{thm:exactCoarse}
Under Assumptions~\ref{ass:subspace}, \ref{ass:notruncate}, and \ref{ass:isomorphic}, the coarse
propagator is exact when applied to the state, i.e.,
 $
\crspropfore{\Tmp}{\Tm}{\stateNP(\Tm)} = \fnprop{\Tmp}{\Tm}{\stateNP(\Tm)} =\stateNP(\Tmp)
,\ \timestepdummy \innatZero{\nptcrs-1}.
  $
\end{theorem}

\subsection{Speedup
analysis}\label{sec:speedup} This section analyzes the theoretical speedup
of the method under various conditions. Section \ref{sec:speedupGeneral}
provides the theoretical speedup of the methodology achieved for a given
number of parareal iterations when both the local forecast (Theorem
\ref{thm:localForecastGen}) and global forecast (Theorem
\ref{thm:globalForecastGen}) are employed for initialization. Section
\ref{sec:idealSpeedups} derives theoretical speedups for the method under
`ideal conditions' for both the local-forecast (Theorem
\ref{thm:idealSpeedup}) and global-forecast (Theorem
\ref{thm:idealSpeedupGlobal}) initializations. \ourReReading{Appendix
\ref{sec:idealSpeedupsNewtonSolver} shows that}
 the proposed method can produce super-ideal
theoretical speedups when the forecast is also employed for providing initial
guesses for the Newton solver in the case of implicit fine propagators and
nonlinear dynamical systems.

Each theoretical result employs a subset of the following assumptions:
\begin{Assumption}[resume=assumption]
\item\label{ass:localinitialize} Initialization in Step \ref{step:initialize} of Algorithm
\ref{alg:parareal} is computed via local forecasting (i.e.,
Eq.~\eqref{eqn:standard_initialization} with
$\crspropsym\leftarrow\crspropforesymAll$).
\item\label{ass:globalinitialize} Initialization in Step \ref{step:initialize} of Algorithm
\ref{alg:parareal} is computed via global forecasting (i.e., Eq.~\eqref{eq:globalForecastInitialize}).
\item\label{ass:integrateCostDominant} The wall time incurred by
computing time advancement with the fine propagator
$\fnprop{\tptfn{m+1}}{\tptfn{m}}{\stateNP(\tptfn{m})}$ dominates all other
costs and parallel overhead.
\end{Assumption}
Further, all speedup results assume that the number of processors is equal to the
number of coarse time intervals $\nptcrs$.
\subsubsection{General case}\label{sec:speedupGeneral}
\begin{theorem}[\textit{Speedup}: local-forecast
initialization]\label{thm:localForecastGen}
If Assumptions \ref{ass:localinitialize} and \ref{ass:integrateCostDominant}
hold, then the proposed method (which employs the local forecast for
initialization and coarse propagation) \ourReReading{upon convergence in
$\pararealItConverge$ parareal iterations} realizes a speedup of
\begin{equation}\label{eq:speedupLocal}
\speedupLocal{\pararealItConverge} \defeq 
\frac{\nptfn}{\left[\memory(\nptcrs-\frac{\pararealItConverge}{2})+
	\nptfnArg{\timestepdummy }\right](\pararealItConverge+1)
-\memory(\pararealItConverge+1)}.
\end{equation}
\end{theorem}

\begin{theorem}[\textit{Speedup}: global-forecast
initialization]\label{thm:globalForecastGen}

If Assumptions \ref{ass:globalinitialize} and \ref{ass:integrateCostDominant}
hold, then the proposed method (which employs the global forecast for
initialization and the local forecast for coarse propagation) 
\ourReReading{upon convergence in
$\pararealItConverge$ parareal iterations}
realizes a speedup of
\begin{equation}\label{eq:speedupGlobal}
\speedupGlobal{\pararealItConverge} \defeq 
\frac{\nptfn}{
\memory(1 + \indicator{\pararealItConverge> 0})	+
\memory\pararealItConverge(\nptcrs-\frac{1}{2}(1+\pararealItConverge))
+(\pararealItConverge+1)\nptfnArg{\timestepdummy
}-\memory\pararealItConverge}.
\end{equation}
Here, the indicator function is defined as
$\indicator{A}= 1$ if $A$ is true, while $\indicator{A}= 0$ otherwise.
\end{theorem}

\begin{remark}[Memory tradeoff: iteration count and
speedup]\label{rem:memoryTradeoff}
Eqs.~\eqref{eq:speedupLocal} and \eqref{eq:speedupGlobal} demonstrate that
increasing the memory $\memory$ can reduce the speedup of the
methodology, assuming the number of iterations $\pararealItConverge$ needed
for convergence is constant. However, as discussed in Remark
	\ref{rem:interpOversampling}, increasing the memory also leads to a
	non-increasing bound for the error between coarse and fine propagators,
	which can (in practice) promote convergence, thereby reducing the number of
	iterations $\pararealItConverge$. These two effects constitute the tradeoff
	that should be considered when selecting the memory $\memory$ in practice.
\end{remark}

\begin{remark}[Reuse of sampled state]
We note that the $\memory$ applications of the fine propagator employed by the
local-forecast coarse propagator to sample the restricted state can be reused
during the subsequent fine propagation; this leads to speedup improvements
as manifested in terms
$-\memory(\pararealItConverge+1)$ and $-\memory\pararealItConverge$ in the
denominators of 
Eqs.~\eqref{eq:speedupLocal} and \eqref{eq:speedupGlobal}, respectively. This
is also an important aspect of the practical implementation of the
local-forecast coarse propagator.
\end{remark}

\subsubsection{Ideal case}\label{sec:idealSpeedups} We now derive theoretical speedups for the
method under `ideal conditions', i.e., when Assumption \ref{ass:subspace}
holds.
\begin{theorem}[\textit{Ideal-conditions speedup}: local-forecast initialization]\label{thm:idealSpeedup}
If Assumptions 
\ref{ass:subspace}, \ref{ass:notruncate}, \ref{ass:isomorphic}, and \ref{ass:localinitialize} hold,
then the proposed method converges after initialization
(i.e., $\pararealItConverge=0$ in Algorithm
\ref{alg:parareal}). Further, if Assumption \ref{ass:integrateCostDominant} holds, then the method realizes a speedup of
$
\speedupLocal{0} = \nptfn/\left((\nptcrs-1)\memory +
\nptfnArg{\timestepdummy }\right).
$
\end{theorem}

\begin{figure}[htbp] 
\centering 
\subfigure[Local-forecast initialization]{
\includegraphics[width=0.45\textwidth]{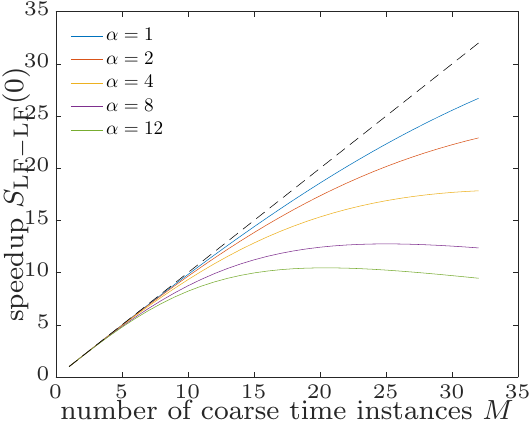} 
\label{fig:theoretical_speedup_local}
}
\subfigure[Global-forecast initialization]{
\includegraphics[width=0.45\textwidth]{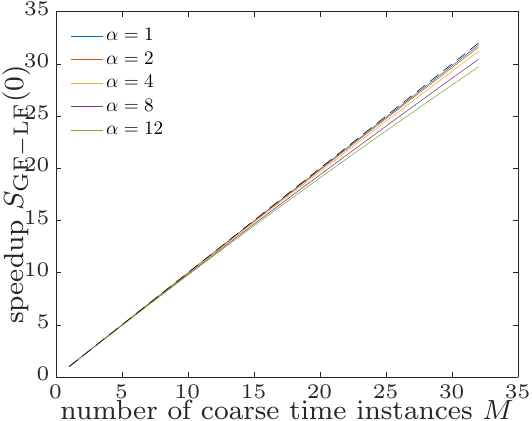} 
\label{fig:theoretical_speedup_global}
}
\caption{Ideal-conditions speedup. Plot
corresponds to $\nptfn=5000$ fine time instances, setting the number of processors
equal to the number of coarse time instances $\nptcrs$.} 
\label{fig:theoretical_speedup}
\end{figure}

Figure \ref{fig:theoretical_speedup_local} provides a visualization of this
theoretical speedup for specific values of method parameters.  First, note
that the `serial bottleneck' of time evolution is apparent from this result:
the speedup degrades as the number of coarse time instances $\nptcrs$
increases. This is due to the requirement of computing $\memory$ fine
propagations in serial across coarse time intervals for this initialization
method. Second, note that the memory $\memory$ has an appreciable effect on
the speedup; keeping this value as low as possible without compromising
convergence is thus desirable.

\begin{theorem}[\textit{Ideal-conditions speedup}: global-forecast
initialization]\label{thm:idealSpeedupGlobal}
If Assumptions \ref{ass:subspace}, \ref{ass:isomorphic}, and \ref{ass:globalinitialize}
 hold,
%
then the proposed method converges after parareal initialization (i.e., 
$\pararealItConverge=0$ in Algorithm \ref{alg:parareal}). Further, if
Assumption \ref{ass:integrateCostDominant} holds, then the method realizes a
speedup of 
\begin{equation}\label{eq:idealSpeedupGlobal}
 \speedupGlobal{0} = \frac{\nptfn}{\memory + \nptfnArg{\timestepdummy }}.
 \end{equation}
\end{theorem}

Figure \ref{fig:theoretical_speedup_global} visualizes this
theoretical speedup in the case of global-forecast initialization. As compared
with local-forecast initialization, note that the theoretical speedup
realizable by the global forecast is much closer to ideal. Further,
it is more stable as discussed in Remark \ref{remark:initialSeedStable}.

\reviewerB{To illustrate the full potential of the proposed approach,
	Appendix \ref{sec:idealSpeedupsNewtonSolver} demonstrates that super-ideal speedups
	can be realized when the proposed methodology is combined with the method
	presented in Ref.\
\cite{Carlberg_carlberg2015decreasing} for defining Newton-solver
	initial guesses via forecasting.}

\subsection{Stability analysis}\label{sec:stability}
We begin by providing a general proof for stability of the parareal
recurrence; we then derive specific quantities needed to demonstrate stability
when the proposed forecast is employed as a coarse propagator.
These results employ a subset of the following
assumptions:
\begin{Assumption}[resume=assumption]
\item\label{ass:stable}The fine propagator is stable\footnote{Note that this
	assumption implies that $\stateVar = 0$ is an equilibrium point. The
following analysis also holds when $\|\stateVar\|$ is replaced by $\|\stateVar
- \stateVarEquilibrium\|$ with $\stateVarEquilibrium$ an equilibrium point. The same applies to the bound in Eq.~\eqref{eq:coarseBoundGen}.}, i.e., 
$
\|\fnprop{\timeVarArg{j}}{\timeVarArg{i}}{\stateVar}\|\leq (1 +
\stabilityConstant(\timeVarArg{j}-\timeVarArg{i}))\|\stateVar\|,\
\forall\stateVar\in\RR{\ndof}.
$
\item \label{ass:lipschitzRestrict} The restriction
operators and prolongation operator counterparts are bounded 
\ourReReading{by constants $\lipschitzOrth$, $\lipschitzRestrictEntry{j} \in \R$}, i.e.,
$\|\prolongateOrthArg{\restrictOrthArg{\stateVar}}\|\leq\lipschitzOrth\|\stateVar\|,\
\forall\stateVar\in\RR{\ndof}$ \ourReReading{and}
$\|\restrictjArg{\stateVar}\|\leq\lipschitzRestrictEntry{j}\|\stateVar\|,\ \forall\stateVar\in\RR{\ndof}$.
\end{Assumption}

The following lemma follows some elements of the stability analysis performed
in Ref.\ 
\cite{chen2014use}.
\begin{lemma}[General parareal stability]\label{lem:pararealStability}
If constants $\constCoarse$ and $\stabilityCoarse$ exist such that the coarse
propagator can be bounded as
\begin{equation}\label{eq:coarseBoundGen}
\|\crsprop{\Tmp}{\Tm}{\stateVar}\|\leq\constCoarse(1 +
\stabilityCoarse\reviewerB{\tsscrs})\|\stateVar\|
\end{equation} 
and constants $\constFineCoarse$ and $\stabilityFineCoarse$ exist such that
the difference between the coarse and fine propagators can be bounded as
\begin{equation}\label{eq:fineCoarseBoundGen}
\|\fnprop{\Tmp}{\Tm}{\stateVar}-
\crsprop{\Tmp}{\Tm}{\stateVar}\|\leq\constFineCoarse(1 +
\stabilityFineCoarse\reviewerB{\tsscrs})\|\stateVar\|,
\end{equation} 
then the parareal recurrence \eqref{eq:parareal} is stable, as it satisfies
\begin{align}\label{eq:pararealStableOne}\quad\quad
\|\apxsolFOM{\tixcrs}{\itvar}\|\leq&
	(\constCoarse)^\timestepdummy \exp(\stabilityCoarse  \timestepdummy \tsscrs)\|\initstateNP\| \!+\! 
	\sum_{j=1}^\timestepdummy (\constCoarse)^{\timestepdummy -j}\constFineCoarse\exp(((\timestepdummy -j)\stabilityCoarse \!+\!
\stabilityFineCoarse)\tsscrs)\|\apxsolFOM{j-1}{\itvar-1}\|,\
\timestepdummy \!\innat{\nptcrs}, \pararealit \!\innat{\timestepdummy }\\
\label{eq:pararealStableTwo}
\|\apxsolFOM{\tixcrs}{\itvar}\|\leq&
\sum_{j=0}^{\tixcrs} { \tixcrs\choose
j}(\constCoarse)^j(\constFineCoarse)^{\tixcrs-j}\exp\left(\left(j\stabilityCoarse +
(\timestepdummy -j)\stabilityFineCoarse\right)\tsscrs \right)\|\initstateNP\|,\quad
\timestepdummy \innat{\nptcrs},\ \pararealit  = \timestepdummy .
\end{align}
\end{lemma}

We now derive the quantities $\constCoarse$, $\stabilityCoarse$,
$\constFineCoarse$, and $\stabilityFineCoarse$ from Lemma
\ref{lem:pararealStability} that are specific to the proposed coarse
propagator $\crspropforesymAll$.

\ourReReading{
\begin{lemma}[Stability of proposed coarse propagator]\label{lem:coarseStable}
Under Assumptions \ref{ass:lipschitzProlong}, \ref{ass:stable}, and \ref{ass:lipschitzRestrict} 
\begin{align}\label{eq:coarseboundTwo}
\|\crspropforesymAll(\stateVar;\tptcrs{\tixcrs},\tptcrs{\tixcrs +
	1})\|\leq\sqrt{\nrestrict}
\lipschitzProlong\max_{j\innat{\nrestrict}}\lipschitzRestrictEntry{j}\constCoarseLFnj\left(1 +
\stabilityCoarseLFnj\stabilityConstant\tsscrs
 \right)\|\stateVar\|,
\end{align}
where we have defined
$\constCoarseLFnj
\defeq\left(\normedQuantityGappySmall{j}{\timestepdummy }
+
\normedQuantityGappy{j}{\timestepdummy }\sqrt{\memory}\right)\geq 1
$,
$\stabilityCoarseLFnj\defeq\frac{(\memory/\nfinepercoarse)\normedQuantityGappy{j}{\timestepdummy }\sqrt{\memory}}{\normedQuantityGappySmall{j}{\timestepdummy }
+
\normedQuantityGappy{j}{\timestepdummy }\sqrt{\memory}}\leq 1
$,
$\normedQuantityGappySmall{j}{\timestepdummy }\defeq|1-\sum_{i=1}^\memory\forecastScalar{i}{j}|$ and
 $ \normedQuantityGappy{j}{\timestepdummy } \defeq
\sqrt{\sum_{i=1}^\memory(\forecastScalar{i}{j})^2}
$.
Hence, \ourReReading{we} have proven stability in 
\ourReReading{the} sense of Eq.~\eqref{eq:coarseBoundGen}.
\end{lemma}
}

\begin{figure}[htbp] 
\centering 
\subfigure[scaling of terms in stability constant]{
\includegraphics[width=0.3\textwidth]{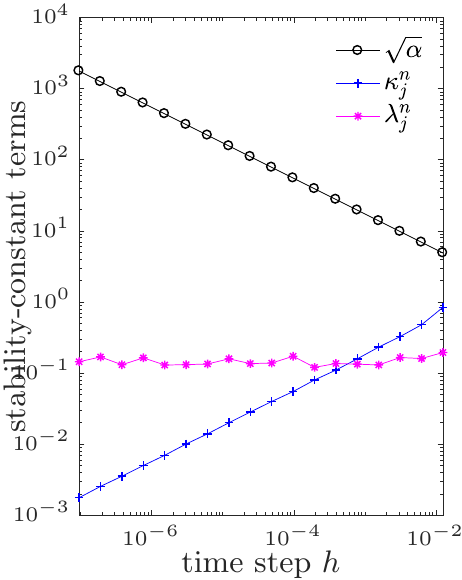} 
\label{fig:terms_converge}
}
\subfigure[scaling of the stability constants]{
\includegraphics[width=0.29\textwidth]{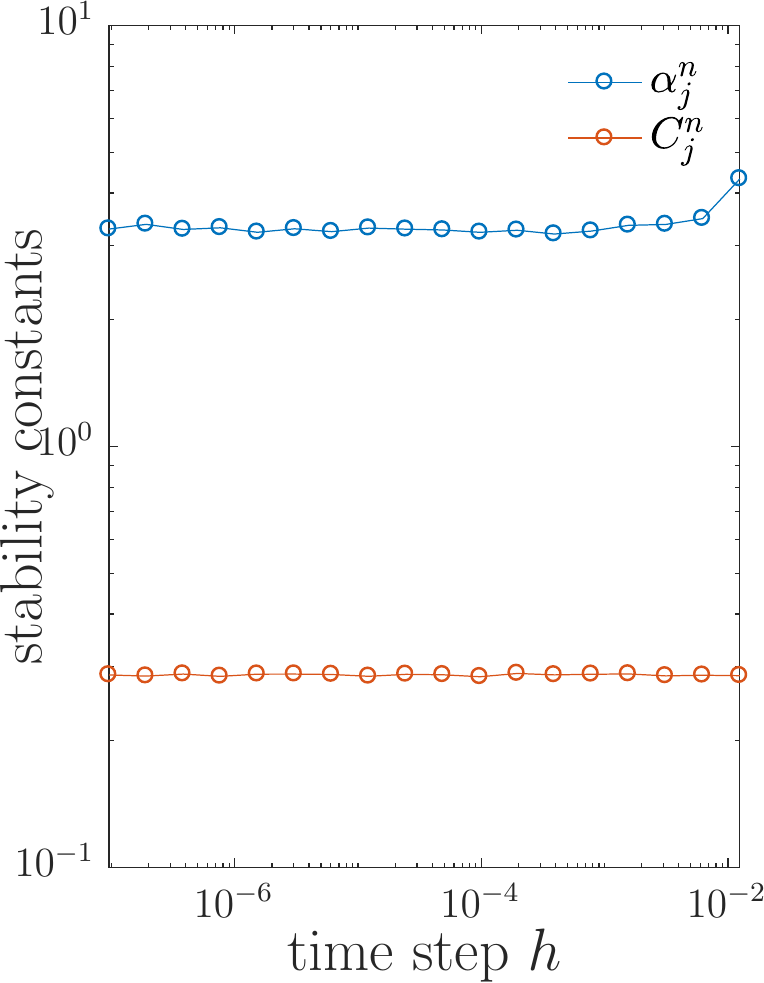} 
\label{fig:stab_constant_converge}
}
\caption{Scaling of the stability constants in bound \eqref{eq:coarseboundTwo}.
Data correspond to
a coarse time interval $\reviewerB{\tsscrs}=1$,
a \ourReReading{time-sampling fraction of $\memory/\nfinepercoarse = 0.3$};
reported values correspond to averages taken over 50 random orthogonal
matrices $\timebasisjm$.
Note
that the stability \ourReReading{constants approach constant values} as the fine time step
\ourReReading{decreases}, which \ourReReading{implies that} bound \eqref{eq:coarseboundTwo} \ourReReading{is} independent of the time discretization for
sufficiently small time steps.} 
\label{fig:stab_constant_converge_overall}
\end{figure}

\begin{remark}[Bound dependence on
discretization]\label{eq:discretizationDependence}
For a fixed 
\ourReReading{coarse time step $\tsscrs$ and time-sampling fraction $\memory/\nfinepercoarse$},
the only \ourReReading{quantities} in bound \eqref{eq:coarseboundTwo} that \ourReReading{depend} on the underlying
(fine) time discretization (i.e., $\tssfn$,
$\nTimestepsIntervalArg{\timestepdummy }$) \ourReReading{are} the stability constants \ourReReading{$\constCoarseLFnj$ and $\stabilityCoarseLFnj$}.
We now assess the dependence of \ourReReading{these stability constants} on the time
discretization.
For a fixed sampling time interval,
the stability constants \ourReReading{$\constCoarseLFnj$ and $\stabilityCoarseLFnj$} approach constant values as the
time step approaches zero. This can be seen from the scaling of the terms
that compose the constant as the (fine) time step $\tssfn$
decreases: $\sqrt{\memory}$ increases with an exponential power
of 1/2, $\normedQuantityGappy{j}{\timestepdummy }$ decreases with an exponential power of 1/2,
and $\normedQuantityGappySmall{j}{\timestepdummy }$ is constant. The
\ourReReading{second of these trends arises} from the fact that \ourReReading{the columns of the matrix}
$\timebasisjm\in\orthogonalMatrix{\nptfnArg{\timestepdummy }}{\dimBasisjm}$ remain
orthogonal when the fine time step $\tssfn$ changes.
Figure \ref{fig:stab_constant_converge_overall} reports a numerical investigation of these trends.
\end{remark}
 

\begin{remark}[Superior stability of global-forecasting initialization to
local-forecasting initialization]\label{remark:initialSeedStable}
We now consider the implications of Lemma \ref{lem:coarseStable} in terms of
the two initialization methods proposed in Section \ref{sec:initialSeed}. 
The first proposal involved applying the local forecast for initialization,
i.e., computing the initial values $\apxsolFOM{\timestepdummy }{0}$, $\timestepdummy \in\nat{\nptcrs}$ via Eq.~\eqref{eqn:standard_initialization} with 
$\crspropsym\leftarrow\crspropforesymAll$.
Applying inequality \eqref{eq:coarseboundTwo} to Eq.~\eqref{eqn:standard_initialization} with 
$\crspropsym\leftarrow\crspropforesymAll$ leads to the following stability
result for the computed initial values:
$
\|\apxsolFOM{\timestepdummy }{0}\|\leq(\stabilityCoarseProp)^\timestepdummy \|\initstateNP\|
,\
\timestepdummy \innat{\nTimeIntervals}.
$
Here, \ourReReading{
$\stabilityCoarseProp \defeq
\sqrt{\nrestrict}
\lipschitzProlong\max_{j\innat{\nrestrict},n\innatZero{\nptcrs-1}}\lipschitzRestrictEntry{j}\constCoarseLFnj\left(1 +
\stabilityCoarseLFnj\stabilityConstant\tsscrs
 \right)$};
that is, the stability factor associated with local-forecast initialization
grows exponentially in the number of coarse time instances $\timestepdummy $. This
phenomenon can be
interpreted as follows: small errors in a local forecast can be amplified by
subsequent local forecasts, as these are performed sequentially.

On the
other hand, by comparing 
Eqs.~\eqref{eq:globalForecastInitialize} and
\eqref{eq:dataDrivenCoarse1}, one can note that
global-forecast initialization \eqref{eq:globalForecastInitialize} is equivalent
to applying the local forecast with global time-evolution bases $\timebasisj$
over a time interval $\Tmp-\tptcrs{0}$. Thus,
the stability of the global-forecast initialization can be derived directly
from inequality \eqref{eq:coarseboundTwo} applied with these modifications as
\begin{equation}\label{eq:stabilityGlobal}
\|\apxsolFOM{\timestepdummy +1}{0}\|\leq\stabilityCoarsePropGlobal\|\initstateNP\|
,\quad
\timestepdummy \innatZero{\nTimeIntervals-1}.
	\end{equation}
Here, we have defined
\ourReReading{
$$ \stabilityCoarsePropGlobal \defeq
\sqrt{\nrestrict}
\lipschitzProlong\max_{j\innat{\nrestrict},n\innatZero{\nptcrs-1}}\lipschitzRestrictEntry{j}\constCoarseGFnj\left(1 +
\stabilityCoarseGFnj\stabilityConstant\tsscrs
 \right),
$$
$\constCoarseGFnj
\defeq\left(\normedQuantityGappySmallGlobal{j}{\timestepdummy }
+
\normedQuantityGappyGlobal{j}{\timestepdummy }\sqrt{\memory}\right)\geq 1
$,
$\stabilityCoarseGFnj\defeq\frac{(\memory/\nfinepercoarse)\normedQuantityGappyGlobal{j}{\timestepdummy }\sqrt{\memory}}{\normedQuantityGappySmallGlobal{j}{\timestepdummy }
+
\normedQuantityGappyGlobal{j}{\timestepdummy }\sqrt{\memory}}\leq 1
$,
$\normedQuantityGappySmallGlobal{j}{\timestepdummy }\defeq|1-\sum_{i=1}^\memory\forecastScalarGlobal{i}{j}|$, and
 $ \normedQuantityGappyGlobal{j}{\timestepdummy } \defeq
\sqrt{\sum_{i=1}^\memory(\forecastScalarGlobal{i}{j})^2}
$.
} Inequality \eqref{eq:stabilityGlobal} shows
that the stability factor associated with global-forecast initialization does not grow with the number of
coarse time instances; it depends on the coarse time instance only through the
\ourReReading{quantities $\normedQuantityGappySmallGlobal{j}{\timestepdummy }$ and 
$ \normedQuantityGappyGlobal{j}{\timestepdummy }$}, which should not grow with
$\timestepdummy $. This
phenomenon can be
interpreted as follows: small forecasting errors cannot be amplified,
as a single forecast is employed for the entire time interval.
\end{remark}

\begin{lemma}[Stability of difference between fine and proposed coarse propagators]\label{lem:fineCoarseStable}
If 
Assumptions \ref{ass:prolongOrtho}, \ref{ass:lipschitzProlong}, 
\ref{ass:stable},
and
\ref{ass:lipschitzRestrict}
hold, then\ourReReading{
\begin{align} 
\begin{split}
\|\fnprop{\Tmp}{\Tm}{\stateVar} -
\crspropfore{\Tmp}{\Tm}{\stateVar}
\|\leq 
(\lipschitzOrth+
\sqrt{\nrestrict}\lipschitzProlong
\max_{j\innat{\nrestrict}}\lipschitzRestrictEntry{j}(\constCoarseLFnj+1))\left(1 +
\stabilityConstant\tsscrs
 \right)\|\stateVar\|.
\end{split}
\end{align} 
}
\end{lemma}

\begin{theorem}[\ourReReading{Parareal stability with proposed coarse
	propagator}]\label{thm:finalStability}
Under Assumptions \ref{ass:prolongOrtho}, \ref{ass:lipschitzProlong},
\ref{ass:stable}, and
\ref{ass:lipschitzRestrict},
employing the proposed coarse propagator in the parareal recurrence (i.e.,
$\crspropsym\leftarrow\crspropforesymAll$ in Eq.~\eqref{eq:parareal}) yields
a stable recurrence, as the iterates satisfy
Eqs.~\eqref{eq:pararealStableOne} and \eqref{eq:pararealStableTwo} with
\ourReReading{
\begin{gather}
\constCoarse = \constCoarseComplexNoPar,\quad \stabilityCoarse =
\stabilityCoarseComplexNoPar,\\ \constFineCoarse =
\constFineCoarseComplexNoPar,\quad  \stabilityFineCoarse =
\stabilityFineCoarseComplexNoPar.
\end{gather}
}
\end{theorem}

Remark \ref{eq:discretizationDependence} shows that, for a fixed
\ourReReading{coarse time step $\tsscrs$ and time-sampling fraction $\memory/\nfinepercoarse$}, the stability constants \ourReReading{$\constCoarseLFnj$ and $\stabilityCoarseLFnj$ approach
constant values} as the time step $\tssfn$ approaches zero. As all other quantities
in the stability bounds
\eqref{eq:pararealStableOne}--\eqref{eq:pararealStableTwo} (with coefficients
specified in Theorem \ref{thm:finalStability}) are
independent of the underlying fine time discretization \ourReReading{(i.e.,
$\tssfn$ and
$\nTimestepsIntervalArg{\timestepdummy }$)}, we know the stability result is not
sensitive to the selected fine time step if it is taken to be sufficiently
small.

\reviewerA{
\subsection{Convergence analysis}\label{sec:convergence}
Recall that the proposed approach merely defines alternative techniques for
initialization and coarse propagation for the parareal method. Thus, one might
expect that existing convergence results derived for the parareal methods
still hold in the present context. However, this is not always the case, as
many existing results assume that the coarse propagator corresponds to a time
integrator with a known order of accuracy
\cite{lions2001parareal,bal2005convergence,bal2002parareal,gander2008nonlinear};
the local-forecast coarse propagator cannot be straightforwardly assigned such
an order of accuracy, as its error is bounded by an expression that does not
explicitly depend on 
the coarse time step $\tsscrs$ (see Theorem
\ref{lem:coarseError}).

Instead, we can make use of existing convergence results that require only
\textit{a general
definition} of the coarse propagator. One such example is the convergence
analysis of 
Ref.~\cite{gander2007analysis}, which assumes a
fixed coarse time step $\tsscrs$ and assesses convergence as the number of
parareal iterations increases. We proceed by describing how the proposed
initialization and coarse propagator affect these convergence results.

Following Ref.\
\cite{gander2007analysis}, we now consider a simplified problem
setting that relies on the following assumptions:
\begin{Assumption}[resume=assumption]
\item\label{ass:linearscalar} Problem~\eqref{eq:ODE} is scalar and linear
i.e.,
$\ndof = 1$ and $\fscalar:
(\stateVarScalar;\tvar)\mapsto\linearFscalar\stateVarScalar$ with
$\linearFscalar\in\RR{}$.
\item \label{ass:coarsePropLinear} The coarse propagator is a linear operator, i.e.,
$\crspropsym:(\stateVarScalar;\tptcrs{i},\tptcrs{j})\mapsto(\linearCrspropscalar)^{j-i}\stateVarScalar$
with $\linearCrspropscalar\in\RR{}$.
\item \label{ass:coarsePropStable} The coarse propagator is in its region of absolute stability
such that $|\linearCrspropscalar|<1$.
\end{Assumption}
We note that under Assumption \ref{ass:linearscalar}, problem
\eqref{eq:ODE} simplifies to
\begin{align}
\label{eq:ODEscalarlinear} 
\frac{d}{dt}
{\stateNPExactscalar} (\tvar) &= \linearFscalar \stateNPExactscalar,\quad
\stateNPExactscalar (0) = \initstateNPscalar \ ,
\end{align}
where $\stateNPExactscalar:\RRplus \rightarrow \RR{}$ and the fine propagator
becomes a linear
operator satisfying
$\fnpropsym:(\stateVarScalar;\timeVarArg{i},\timeVarArg{j})\mapsto(\linearFnpropscalar)^{j-i}\stateVarScalar$
with $\linearFnpropscalar\in\RR{}$.  For example, the backward-Euler fine
propagator becomes
$\fnpropsymBE:(\stateVarScalar;\timeVarArg{i},\timeVarArg{j})\mapsto(\linearFnpropBEscalar)^{j-i}\stateVarScalar$
with $\linearFnpropBEscalar\defeq(1-\tssfn\linearFscalar)^{-1}$. 
Theorem 4.5 (with Corollary 4.6) of that reference is repeated below
(including modifications discussed in Section 4.5 of Ref.\
\cite{gander2007analysis}) in
the current notation.
\begin{theorem}[Parareal convergence (Theorem 4.5 and Corollary 4.6 of
	\cite{gander2007analysis})]\label{thm:pararealConvergence}
	If $\finalT < \infty$ and Assumptions \ref{ass:linearscalar} and
	\ref{ass:coarsePropLinear}	
	hold, then 
	\begin{equation}\label{eq:convergeGeneral}
		\max_{n\innat{\nptcrs}}|\stateNPscalar(\tptcrs{n}) - 
		\apxsolFOMscalar{n}{k}
		|
		\leq
		|(\linearFnpropscalar)^\nfinepercoarse-\linearCrspropscalar|^k\|(\Mmat(\linearCrspropscalar))^k\|_\infty
		\max_{n\innat{\nptcrs}}|\stateNPscalar(\tptcrs{n}) - 
		\apxsolFOMscalar{n}{0}|,
	\end{equation}
	where $\Mmat\in\RR{\nptcrs\times\nptcrs}$ is a Toeplitz matrix whose
	elements are defined by the value of the elements in the first column, which
	are
\begin{align} 
	\MmatEntry{i}{1}(\beta) = \begin{cases} 0,\quad &i=1\\
	\beta^{i-2},\quad &i\in\{2,\ldots,\nptcrs\}
\end{cases}
\end{align} 
If additionally Assumption \ref{ass:coarsePropStable} holds, the recurrence converges superlinearly
as
\begin{equation}\label{eq:convergeSuperlinear}
		\max_{n\innat{\nptcrs}}|\stateNPscalar(\tptcrs{n}) - 
		\apxsolFOMscalar{n}{k}
		|
		\leq
		\frac{|(\linearFnpropscalar)^\nfinepercoarse-\linearCrspropscalar|^k\prod_{j=1}^k(\nptcrs-j)}{k!}		\max_{n\innat{\nptcrs}}|\stateNPscalar(\tptcrs{n}) - 
		\apxsolFOMscalar{n}{0}|.
	\end{equation}
\end{theorem}

We now describe how our prescribed coarse propagator can be integrated in this
convergence result.
First, we note that under Assumption \ref{ass:linearscalar}, the proposed coarse
propagator is characterized by 
$\restrictT=\prolongate= 1$ and 
$\nrestrict = \ndof = 1$.
We now collect assumptions related to the proposed coarse propagator:
\begin{Assumption}[resume=assumption]
\item\label{ass:samebasis} The same local basis is employed for every coarse
	time interval, i.e., $\timebasismArg{i} = \timebasismArg{j} =
	\timebasis\in
	\orthogonalMatrix{\nptfnArg{\timestepdummy }}{\dimBasisjm}
	$,
	$i,j\innatZero{\nptcrs-1}$.
\item \label{ass:coarseContract} The forecast satisfies the inequality 
	$|1+\sum_{i=1}^\memory\forecastScalarScalar{i}
 [(\linearFnpropscalar)^i-1]| < 1$ with 
$\forecastScalarScalar{i}\defeq\unitvec{\nTimestepsIntervalArg{\timestepdummy }}^T\timebasis\left[\sampleMat{
	0}{\memory}\timebasis\right]^+\unitvec{i}\in\RR{}$, $i\innat{\memory}$.
\end{Assumption}
We now show that the parareal recurrence executed with the proposed coarse
propagator converges superlinearly under the stated conditions.
\begin{corollary}[Superlinear parareal convergence using the proposed coarse
	propagator]\label{cor:convergenceProposed}
Under Assumptions \ref{ass:linearscalar} and \ref{ass:samebasis}, 
the proposed coarse propagator is linear and satisfies $
\crspropforesymAll:(\stateVarScalar;\tptcrs{\timestepdummy },\tptcrs{\timestepdummy  +
	1})\mapsto
	\linearCrspropscalarLF\stateVarScalar
	$ with $\linearCrspropscalarLF \defeq 1+\sum_{i=1}^\memory\forecastScalarScalar{i}
 [(\linearFnpropscalar)^i-1]$.
Further, the error in
the parareal recurrence executed with the proposed coarse propagator
$\crspropsym\leftarrow\crspropforesymAll$
satisfies
\begin{align}\label{eq:corollaryOursConvergenceFirst}
		\max_{n\innat{\nptcrs}}|\stateNPscalar(\tptcrs{n}) - 
		\apxsolFOMscalar{n}{k}
		|
		&\leq
		(\linearCrspropscalarLFError)^k\|(\Mmat(\linearCrspropscalarLF))^k\|_\infty
		\max_{n\innat{\nptcrs}}|\stateNPscalar(\tptcrs{n}) - 
		\apxsolFOMscalar{n}{0}|.
	\end{align}
	where
	$\linearCrspropscalarLFError\defeq|(\linearFnpropscalar)^\nfinepercoarse-
	\linearCrspropscalarLF|=|(\linearFnpropscalar)^\nfinepercoarse-1-\sum_{i=1}^\memory\forecastScalarScalar{i}
	[(\linearFnpropscalar)^i-1]|
	$.
	If Assumption \ref{ass:coarseContract} additionally holds, 
	then the recurrence converges superlinearly as
	\begin{align}\label{eq:corollaryOursConvergenceSecond}
		\max_{n\innat{\nptcrs}}|\stateNPscalar(\tptcrs{n}) - 
		\apxsolFOMscalar{n}{k}
		|
		&\leq
		\frac{(\linearCrspropscalarLFError)^k\left(\prod_{j=1}^k(\nptcrs-j)\right)}{k!}
		\max_{n\innat{\nptcrs}}|\stateNPscalar(\tptcrs{n}) - 
		\apxsolFOMscalar{n}{0}|.
	\end{align}
\end{corollary}
\begin{remark}[Role of accuracy in convergence]
	Inequalities
	\eqref{eq:corollaryOursConvergenceFirst}--\eqref{eq:corollaryOursConvergenceSecond}
	demonstrate the effect of coarse-propagation and initial-seed accuracy on
	convergence. In particular, the term $\linearCrspropscalarLFError$ represents the
	error the coarse propagator incurs with respect to the fine propagator; this
	is precisely the quantity we aim to minimize with the proposed coarse
	propagator. In fact, Theorem \ref{lem:coarseError} bounds this error, and 
	Theorem \ref{thm:exactCoarse} demonstrates that this error is zero under
	`ideal conditions'. Further, the error incurred by the
	initial seed appears as $\max_{n\innat{\nptcrs}}|\stateNPscalar(\tptcrs{n}) - 
		\apxsolFOMscalar{n}{0}|$ in these results. This is the term we aim to
		minimize by applying the proposed local and global forecasting methods for
		initialization; this quantity is also zero under the `ideal conditions'
		stated in Theorem \ref{thm:exactCoarse}.
\end{remark}
}

\section{Computing forecasting ingredients via SVD/POD}
\label{sec:POD} 
We now describe how the three ingredients that define the proposed
methodology---the time-evolution bases $\timebasisj$, $j\innat{\ndof}$, the
restriction operator $\restrictT$, and the prolongation operator
$\prolongate$---can be computed using the POD method.  Section \ref{sec:FOM} describes this for the case of
parameterized ODEs, while Section \ref{sec:ROM} specializes this for the case
of POD-based ROMs.

\subsection{Parameterized ODEs}\label{sec:FOM}

We first introduce a parameterized variant of the governing initial-value
ordinary-differential-equation (ODE) problem \eqref{eq:ODE}.

\begin{align}
\label{eq:ODEParam}
\frac{d}{dt}{\stateExact} (\tvar, \parv) =
\fArg{\stateExact}{\tvar}{\parv},\quad
\stateExact(0, \parv) = \initstate (\parv) \ ,
\end{align}
where
$\parv \in \parsp \subset \R^{\pardim}$ denotes the
parameters,
$\stateExact:\RRplus \times \parsp\rightarrow \R^{\ndof}$ denotes the
(parameterized) state implicitly
defined as the exact solution to problem~\eqref{eq:ODEParam},
$\ivpfunc: \R^{\ndof} \times \RRplus \times \parsp\rightarrow
\R^{\ndof}$ with $(\stateVar;\timeVar,\paramVar)\mapsto
\ivpfunc(\stateVar;\timeVar,\paramVar)$ denotes the
velocity, and
$\initstate:\parsp\rightarrow\RR{\ndof}$
denotes the initial state. Analogously to Eq.~\eqref{eq:propagatorExact}, we
define 
$\state(\cdot,\parv):\timeVar\mapsto\fnprop{\timeVar}{0}{\initstate(\parv)}$ for 
\ourReReading{$\timeVar \in \timeInstanceSet$}
as the associated numerical solution with 
$\state(\cdot,\parv)\in(\setOfTimeDepFunctions)^\ndof$.

\begin{algorithm}[tb]
\caption{\podalgNo}
\begin{algorithmic}[1]\label{alg:podTemporal}
\REQUIRE 
training parameter instances $\paramTrain\subset \paramDomain$,
$\energyCrit\in\left[0,1\right]$
\ENSURE POD state basis $\podStateBasis\in\orthogonalMatrix{\ndof}{\nrestrictPOD}$,
POD time-evolution bases
$\podTemporalBasisj\in\orthogonalMatrix{\nTimestepsFine}{\ntrain}$,
$j=1,\ldots,\nrestrictPOD$
\FOR[collect snapshots]{$i=1,\ldots,\ntrain$}\label{step:snapshotsFirst}
\STATE Numerically solve Eq.~\eqref{eq:ODEParam} with
$\parv\leftarrow\paramTraini{i}$ to obtain snapshots \\
$\parsnpsht{i}\defeq
[\state(\tptfn{1}, \paramTraini{i})- \initstate(\paramTraini{i})\ \cdots\ 
\state(\tptfn{\nptfn}, \paramTraini{i})- \initstate(\paramTraini{i})
] \in \RR{\ndof \times \nptfn}
$
\ENDFOR\label{step:snapshotsLast}
\STATE\label{step:svdPOD} $(\U,\Sig,\V) = \texttt{thin\_SVD}(\left[\parsnpsht{1}\ \cdots\
\parsnpsht{\ntrain}\right])$ \COMMENT{Compute (thin) singular value decomposition}
\STATE\label{step:truncateSVDPOD} $\podStateBasis\leftarrow
\vectomat{\uvec}{\nrestrictPOD}$, where 
$\nrestrictPOD = \min_{i\in\energyCritSet(\energyCrit)} i$, 
$\energyCritSet(\energyCrit) \defeq \{i\ |\
\sum_{j=1}^i\singularValue_j/\sum_{k=1}^{\nptfn\ntrain}\singularValue_k\geq
\energyCrit\}$, \\ $\Sig=\mathrm{diag}(\singularValue_1,\ldots,
\singularValue_{\nptfn\ntrain})$. \COMMENT{Compute truncated state basis}
\FOR[Extract temporal bases from right singular
vectors]{$j=1,\ldots,\nrestrictPOD$}
\STATE \ourReReading{$(\podTemporalBasisQj,\podTemporalBasisRj) = 
\texttt{thin\_QR}(
[\podTemporalBasisVec{j}{1}\ \cdots\
\podTemporalBasisVec{j}{\ntrain}]
)
$ with
$\podTemporalBasisVec{j}{i}\defeq[v_{j,\nptfn(i-1)+1}\ \cdots\ v_{j,\nptfn i
}]^T$}
\COMMENT{Compute (thin) QR factorization}
\STATE \ourReReading{$\podTemporalBasisj\leftarrow \podTemporalBasisQj$}
\ENDFOR
\end{algorithmic}
\end{algorithm}
The ingredients required for the proposed methodology can be computed in a
data-driven manner via the POD method by executing the following steps:
\begin{enumerate} 
\item\label{step:odeDataTrain} Given training parameter instances
$\paramTrain\subset\paramDomain$ and energy criterion
$\energyCrit\in[0,1]$, execute Algorithm \ref{alg:podTemporal} to obtain 
POD state basis $\podStateBasis\in\orthogonalMatrix{\ndof}{\nrestrictPOD}$ and
POD time-evolution bases
$\podTemporalBasisj\in\orthogonalMatrix{\nTimestepsFine}{\ntrain}$,
$j\innat{\nrestrictPOD}$.
\item Set the forecasting time-evolution bases equal to the POD time-evolution
bases $\timebasisj\leftarrow\podTemporalBasisj$, $j\innat{\nrestrictPOD}$.
Note that $\nrestrict = \nrestrictPOD$ and 
$\dimBasisj = \ntrain$.
\item\label{step:restrictprolongate} Define the restriction and prolongation operators as
$\restrictT\leftarrow\podStateBasis$ and $\prolongate\leftarrow\podStateBasis$,
respectively.
\end{enumerate}
This approach is sensible, as numerous studies have shown that POD tends to
truncate solution modes associated with high-frequency temporal behavior
\cite{carlbergGalDiscOpt}. Thus, the resulting restriction operator
will ensure that forecasting is applied only to the long-temporal-wavelength
solution components. \ourReReading{We note that this approach is equivalent to
	computing `tailored' temporal subspaces \cite{choiCarlberg} via the sequentially truncated
high-order SVD \cite{vannieuwenhoven2012new}.}

\reviewerA{
\begin{remark}[Ideal predictive case for parameterized linear ODEs]
For illustration, consider a variant of the initial-value ODE problem
\eqref{eq:ODEParam} wherein the velocity is linear in the state but
independent of time and the parameters,
i.e., 
$\f: (\stateVar;\timeVar,\paramVar)\mapsto\linearF\stateVar$, the initial
condition exhibits separable parameter dependence, i.e., 
$
\initstate (\parv) =
\sum_{i=1}^{\nfunctions}\paramFunctioni(\parv)\initialStateContribi
$ with $\paramFunctioni:\parsp\rightarrow\RR{}$ and 
$\initialStateContribi\in\RR{\ndof}$,
$i\innat{\nfunctions}$ linearly independent, and the parameter set is
unbounded, i.e., $\parsp = \R^{\pardim}$.
Then, problem \eqref{eq:ODEParam} becomes
\begin{align}
\label{eq:ODEParamLinear}
\frac{d}{dt}{\stateExact} (\tvar, \parv) =
\linearF\stateExact(\tvar, \parv),\quad
\stateExact (0, \parv) =
\sum_{i=1}^{\nfunctions}\paramFunctioni(\parv)\initialStateContribi.
\end{align}
In this case, the fine propagator is also linear and can be written as
$\fnpropsym:(\stateVar;\timeVarArg{i},\timeVarArg{j})\mapsto(\linearFnprop)^{j-i}\stateVar$.
For example, the backward-Euler fine propagator becomes
$\fnpropsymBE:(\stateVar;\timeVarArg{i},\timeVarArg{j})\mapsto(\linearFnpropBE)^{j-i}\stateVar$
with $\linearFnpropBE\defeq(\identityArg{\ndof}-\tssfn\linearF)^{-1}$.
Therefore, the discrete
solution is simply
\begin{align}\label{eq:discreteFormLinear}
	\stateNP(\tptfn{j},\parv)&= \fnprop{\timeVarArg{j}}{0}{\state (0, \parv)} =
(\linearFnprop)^{j}\state (0, \parv) 
= 
(\linearFnprop)^{j}\sum_{i=1}^{\nfunctions}\paramFunctioni(\parv)\initialStateContribi,\quad
j\innatZero{\nptfn}
\end{align}
Now, assume that  $\ntrain=\nfunctions$ training parameter instances
$\paramTrainNfunctions$
are employed such that the matrix
$\mappingTrainingToBasis\in\RR{\nfunctions\times\nfunctions}$ with elements
$\mappingTrainingToBasisEntry{i}{j}\defeq\paramFunctioni(\paramTraini{j})$ is
invertible. Then, we have 
$
\left[
	\initialStateContribArg{1}\ \cdots\ \initialStateContribArg{\nfunctions}
\right] = 
\left[
	\state(0,\paramTraini{1})\ \cdots\ \state(0,\paramTraini{\nfunctions})
\right]\mappingTrainingToBasis^{-1}
$ and Eq.~\eqref{eq:discreteFormLinear} becomes
\begin{align}\label{eq:discreteFormLinearSub}
	\stateNP(\tptfn{j},\parv)&=
	(\linearFnprop)^{j}\sum_{i=1}^{\nfunctions}\paramFunctionModi(\parv)\state(0,\paramTraini{i})
= 
\sum_{i=1}^{\nfunctions}\paramFunctionModi(\parv)(\linearFnprop)^{j}\state(0,\paramTraini{i})= 
\sum_{i=1}^{\nfunctions}\paramFunctionModi(\parv)
\state(\timeVarArg{j},\paramTraini{i}),\quad
j\innatZero{\nptfn}
\end{align}
where
$\paramFunctionModi(\parv)\defeq\sum_{j=1}^\nfunctions[\mappingTrainingToBasis^{-1}]_{ij}\paramFunctionArg{j}(\parv)$.
Therefore, we have 
$
\stateNP(\tptfn{j},\parv) - \stateNP(0,\parv) = 
\sum_{i=1}^{\nfunctions}\paramFunctionModi(\parv)
(\state(\timeVarArg{j},\paramTraini{i}) - \state(0,\paramTraini{i})
),
$
or equivalently
\begin{equation}\label{eq:conditionLinear}
	\unrollfunc{\stateNPEntry{j}(\cdot,\parv)}\in\range{\left[
		\unrollfunc{
			\stateNPEntry{j}(\cdot,\paramTraini{1})
		}\ \cdots\ 
		\unrollfunc{\stateNPEntry{j}(\cdot,\paramTraini{\nfunctions})}
\right]},\quad j\innat{\ndof},\ \forall \parv\in\parsp.
\end{equation}
Thus, employing $\restrictT=\prolongate=\identityArg{\ndof}$ (such that $\nrestrict = \ndof$) with 
$\timebasisj = \left[
		\unrollfunc{
			\stateNPEntry{j}(\cdot,\paramTraini{1})
		}\ \cdots\ 
		\unrollfunc{\stateNPEntry{j}(\cdot,\paramTraini{\nfunctions})}
	\right]$, $j\innat{\nrestrict}$ in this case ensures that Assumptions \ref{ass:subspace} and
	\ref{ass:isomorphic} hold. Then, if 
	the local bases are constructed with no truncation (i.e., Assumption
	\ref{ass:notruncate} holds), the coarse propagator is exact (see Theorem
	\ref{thm:exactCoarse}). Further, 
	the proposed method converges
	in $\pararealItConverge=0$ iterations if initialization is computed either
	via local forecasting (i.e., Assumption \ref{ass:localinitialize}
	holds; see Theorem
	\ref{thm:idealSpeedup}) or 
	via global forecasting (i.e.,
	Assumption \ref{ass:globalinitialize} holds; see Theorem
	\ref{thm:idealSpeedupGlobal}).\footnote{\reviewerA{We note that rather than employing $\restrictT=\prolongate=\identityArg{\ndof}$ and
$\timebasisj = \left[ \unrollfunc{ \stateNPEntry{j}(\cdot,\paramTraini{1}) }\
\cdots\ \unrollfunc{\stateNPEntry{j}(\cdot,\paramTraini{\nfunctions})}
\right]$, it can also be shown that executing the Steps
\ref{step:odeDataTrain}--\ref{step:restrictprolongate} above with $\energyCrit =
1.0$ to compute
the forecasting bases $\timebasisj$, the restriction operator $\restrictT$,
and the prolongation operator $\prolongate$ also leads to ideal convergence
(i.e., convergence in $\pararealItConverge = 0$ iterations) for both
local-forecast and global-forecast initialization.}}
So this is an example where the forecast of the proposed method is equivalent
to the fine propagator for all parameters $\parv \in \R^{\pardim}$; hence, it
is an ideal predictive coarse propagator.
\end{remark}
}

\subsection{POD-based reduced-order model}\label{sec:ROM}

Projection-based model reduction aims to reduce the cost of numerically
solving Eq.~\eqref{eq:ODEParam} by reducing the dimensionality of the governing
equations. To achieve this, these techniques employ a `trial basis'
$\tribasmtx \in \RRstar{\ndof \times \reddim}$ with reduced state dimension
$\reddim\leq \ndof$, and subsequently approximate the state as $\apxstateExact:
(\tvar, \paramVar) \mapsto \initstate(\parv)+\tribasmtx \redstateExact (\tvar,
\paramVar)$.
Here, 
$\RRstar{m\times n}$ denotes the 
set of full-column-rank $m\times n$ real-valued matrices (i.e., the noncompact
Stiefel manifold), and
the reduced state
$\redstateExact:\RRplus\times\RR{\pardim}\rightarrow\RR{\reddim}$ satisfies
\begin{align}\label{eq:ODEROMParam}
\frac{d}{dt}\redstateExact (\tvar, \parv) = \redfArg{\redstateExact}{\tvar}{\parv},\quad
\redstateExact (0, \parv) = \zero,
\end{align}
where $\redf:(\redstateVar;\timeVar,\paramVar)\mapsto
(\testbasis(\redstateVar;\timeVar,\paramVar)^T
\trialbasis)^{-1} \projmtx(\redstateVar;\timeVar,\paramVar)^T
\fArg{\initstate(\paramVar) + \tribasmtx \redstateVar}{
\timeVar}{ \paramVar}$ denotes the reduced velocity and 
$\testbasmtx:\RR{\reddim}\times\RRplus{}\times\RR{\nparam}\rightarrow
\RRstar{\ndof\times \reddim}$ denotes the `test basis'. Note that
Eq.~\eqref{eq:ODEROMParam} enforces the ODE residual 
$\frac{d}{dt} \apxstateExact(t,\parv) - \fArg{\apxstateExact(t,\parv)}{
t}{ \parv}$
to be orthogonal to $\range{\testbasisTypical}$.
The test basis can be set equal to the trial basis (i.e.,
$\testbasisTypical = \trialbasis$)---which is referred to as Galerkin
projection---or can be chosen to
minimize the discrete residual arising after time discretization (e.g.,
$\testbasisTypical = [\alpha_0\identity -\timestep\beta_0\partial \f/\partial
\stateVar(\initstate+\trialbasis\redstate;\tvar,\param)]\trialbasis$ for linear multistep
schemes, where  $\alpha_0$ and $\beta_0$
are coefficients for a given scheme), which is referred to as least-squares Petrov--Galerkin
projection \cite{CarlbergGappy,carlberg2013gnat,carlbergGalDiscOpt}, for example.
Again, we
define 
$
\redstate(\cdot,\parv):\timeVar\mapsto\fnprop{\timeVar}{0}{\initredstate(\parv)},
$
as the associated numerical solution with 
$\redstate(\cdot,\parv)\in(\setOfTimeDepFunctions)^\reddim$.

When the trial basis $\tribasmtx$ is computed via POD, 
both the trial basis and the
proposed method's ingredients can be computed by
executing the following steps:
\begin{enumerate} 
\item\label{step:romDataTrain} Given training parameter instances
$\paramTrain\subset\paramDomain$ and energy criterion
$\energyCrit\in[0,1]$, execute Algorithm \ref{alg:podTemporal} to obtain 
POD state basis $\podStateBasis\in\orthogonalMatrix{\ndof}{\nrestrictPOD}$ and
POD time-evolution bases
$\podTemporalBasisj\in\orthogonalMatrix{\nTimestepsFine}{\ntrain}$,
$j\innat{\nrestrictPOD}$.
\item\label{step:romDataSetTrial} Set the trial basis equal to the POD state basis
$\tribasmtx\leftarrow\podStateBasis$; note that $\reddim = \nrestrictPOD$.
\item Set the forecasting time-evolution bases equal to the truncated POD time
evolution bases such that only the first
$\nrestrict$ (with $\nrestrict\leq \nrestrictPOD=\reddim$) POD
modes are employed for forecasting: 
$\timebasisj\leftarrow\podTemporalBasisj$,
$j\innat{\nrestrict}$. Note that $\dimBasisj = \ntrain$.
\item\label{step:romDataFinal} Define the restriction and prolongation operators as
$\restrictT\leftarrow[\unitvec{1}\ \cdots\
\unitvec{\nrestrict}]\in\{0,1\}^{\reddim\times\nrestrict}$ and
$\prolongate\leftarrow[\unitvec{1}\ \cdots\ \unitvec{\nrestrict}]\in\{0,1\}^{
	\reddim\times \nrestrict}$,
respectively.
\end{enumerate}
\begin{remark}[Negligible additional
cost and effective use of right singular vectors]
Steps \ref{step:romDataTrain}--\ref{step:romDataSetTrial} above are already
required when the trial basis is computed via POD. Thus, in this case, the
ingredients required for the proposed method can be obtained with negligible
additional computational cost, as the dominant costs in Steps
\ref{step:romDataTrain}--\ref{step:romDataFinal} above are incurred in Step 1.
In particular, these dominant costs comprise (1) collecting snapshots (Steps
\ref{step:snapshotsFirst}--\ref{step:snapshotsLast} in Algorithm
\ref{alg:podTemporal}) and (2) computing the singular value decomposition
(Step \ref{step:svdPOD} in Algorithm \ref{alg:podTemporal}). Thus, one can
interpret the proposed methodology as providing a technique to
\textit{effectively use the right singular vectors}, which are already
available for POD-based reduced-order models after computing the SVD in Step \ref{step:svdPOD} of Algorithm
\ref{alg:podTemporal}.
\end{remark}
\begin{remark}[General reduced-order models]
When the trial basis is not computed via POD, the approach described in
Section \ref{sec:FOM} can be employed, as the reduced-order-model ODE
\eqref{eq:ODEROMParam} has the same structure as the parameterized ODE
\eqref{eq:ODEParam}. In this case, the snapshot collection required in Step
\ref{step:odeDataTrain} incurs a small computational cost, as Steps
\ref{step:snapshotsFirst}--\ref{step:snapshotsLast} of Algorithm
\ref{alg:podTemporal} entails numerically solving only the reduced-order-model
ODE \eqref{eq:ODEROMParam} at parameter instances $\paramTrain$.
\end{remark}

\section{Numerical experiments}\label{sec:experiments}

This section compares the performance of several choices for parareal initialization
and coarse propagation 
in the context of model reduction applied to a parameterized
Burgers' equation. Here, the backward-Euler scheme is employed as the time
integrator that defines the fine propagator; that is, we employ
$\fnpropsym\leftarrow\fnpropsymBE$.
In particular, we consider:
 \begin{itemize} 
	 \item \reviewerB{Four} methods for performing initialization in Step
 \ref{step:initialize} of Algorithm \ref{alg:parareal}: 
\begin{enumerate} 
 \item[(\BELabel)] the backward-Euler
 scheme (Eq.~\eqref{eqn:standard_initialization}
 with $\crspropsym\leftarrow\crspropBE$),
where the coarse
propagator $\crspropBE$ \reviewerB{is first-order accurate and} implicitly satisfies 
 $
\crspropBE(\stateVar;\tptcrs{\tixcrs},\tptcrs{\tixcrs + 1}) - \stateVar -
\tsscrsk\fArgNoParam{\crspropBE(\stateVar;\tptcrs{\tixcrs},\tptcrs{\tixcrs +
1})}{\ourReReading{\tptcrs{\tixcrs + 1}}}=0
	$,
 \item[\reviewerB{(\CNLabel)}]\reviewerB{the Crank--Nicolson
 scheme (Eq.~\eqref{eqn:standard_initialization}
 with $\crspropsym\leftarrow\crspropCN$), 
where the coarse
propagator $\crspropCN$ is second-order accurate and implicitly satisfies 
 $$
\crspropCN(\stateVar;\tptcrs{\tixcrs},\tptcrs{\tixcrs + 1}) - \stateVar -
\frac{1}{2}(\tsscrsk\fArgNoParam{\crspropCN(\stateVar;\tptcrs{\tixcrs},\tptcrs{\tixcrs +
1})}{\tptcrs{\tixcrs + 1}}
+
\tsscrsk\fArgNoParam{\stateVar}{\tptcrs{\tixcrs }}
)=0, $$ }
 \item[(\localForeLabel)] local forecasting (Eq.~\eqref{eqn:standard_initialization}
 with $\crspropsym\leftarrow\crspropforesymAll$), and
 \item[(\globalForeLabel)] global forecasting
 (Eq.~\eqref{eq:globalForecastInitialize}).
\end{enumerate} 
\item \reviewerB{Three} coarse propagators:
\begin{enumerate} 
\item[(\BELabel)] the backward-Euler scheme
	($\crspropsym\leftarrow\crspropBE$),
\item[(\CNLabel)] \reviewerB{the Crank--Nicolson scheme
	($\crspropsym\leftarrow\crspropCN$), and}
\item[(\localForeLabel)]\label{localForeProp} local forecasting
($\crspropsym\leftarrow\crspropforesymAll$).
\end{enumerate}
 \end{itemize}
 We refer to method $i$-$j$ as the method where initialization is carried out with method $i$
 and coarse propagation with method $j$; for example, method
 \globalForeLabel-\localForeLabel\ performs initialization using the
 global forecast and employs the local-forecasting coarse propagator.

\subsection{Parameterized Burgers' Equation and model reduction}

We now describe the parameterized Burgers' equation as described in
Ref.\
\cite{tpwl}, which corresponds to the following parameterized initial
boundary value problem for $(x,\tau)\in [0,100]\times[0,25]$:
\begin{align}\label{eq:burgers}
\begin{split}
\frac{\partial u(x,\tau)}{\partial \tau} + \frac{1}{2}\frac{\partial
\left(u^2\left(x,\tau\right)\right)}{\partial x} &= 0.02e^{\parvEntry{2}x}
\end{split}
\end{align}
\ourReReading{with} $u(0,\tau) = \parvEntry{1}, \ \forall \tau\in\left[0,25\right]$, 
$u(x,0) = 1, \ \forall x\in\left[0,100\right]$,
where the parameter domain corresponds to $\parv = (\parvEntry{1},\parvEntry{2})\in\paramDomain = \left[1.5,2.0\right]\times\left[0.02,0.025\right]$.

After applying Godunov's scheme for spatial discretization with 500 control
volumes, 
\eqref{eq:burgers} 
\ourReReading{and the boundary and initial conditions} lead to a
parameterized initial-value ODE problem consistent with
Eq.~\eqref{eq:ODEParam} with $\ndof=500$ degrees of freedom. As described
earlier, we employ the backward-Euler scheme for time discretization
\ourReReading{using} uniform fine time steps $\tssfnk=0.1$, 
which leads to
$\nptfn=250$ fine time instances. Unless otherwise stated, we 
\ourReReading{set} a
parareal termination tolerance $\jumpTolerance=\ourReReading{5\times 10^{-3}}$ in Algorithm
\ref{alg:parareal}.

We compare the time-parallel methods in the POD-based reduced-order-modeling
context as discussed in Section \ref{sec:ROM}. Here, we employ $\ntrain=4$
randomly-selected training points 
$\paramTraini{1}=(1.5331, 0.0249)$,
$\paramTraini{2}=(1.6880, 0.0223)$,
$\paramTraini{3}=(1.9656, 0.0209)$, and
$\paramTraini{4}=(1.8000, 0.0232)$.
We 
\ourReReading{choose} a
reduced-state dimension\footnote{Note that instead of specifying the energy criterion
$\energyCrit$ as suggested in Section \ref{sec:ROM}, we directly specify the
reduced-state dimension  $\reddim$.} of $\reddim=\nrestrictPOD=100$ and use
the least-squares Petrov--Galerkin (LSPG) ROM \cite{CarlbergGappy},
which---for the backward-Euler case---corresponds 
to a test basis of
$\testbasisTypical = [\identity -\timestep\partial \f/\partial
\stateVar(\initstate+\trialbasis\redstate;\tvar,\param)]\trialbasis$.
During the experiments, we will assess the performance of the ROMs and
time-parallel methods at a set of $\nonline=2$ randomly-selected online parameter instances
$\paramOnlinei{1} = (1.6603, 0.0229)$ and $\paramOnlinei{2} = (1.5025,
0.0201)$; that is, we numerically solve the reduced
initial-value ODE problem \eqref{eq:ODEROMParam} for $\parv\in\paramOnline$.


During the experiments, we vary the number of restricted states
$\nrestrict$ and the forecast memory $\memory$.

\subsection{Comparison of initialization and coarse-propagation methods}

We first compare the performance of \reviewerB{multiple combinations} of
initialization methods and coarse
propagators. To achieve this, we set the number of coarse time  instances to
$\nptcrs = 10$ and employ a parareal tolerance of $\jumpTolerance=0$; this
ensures that the parareal method will execute (the maximum value of)
$\pararealItConverge=\nptcrs-1=9$
parareal iterations in Algorithm \ref{alg:parareal}, thereby allowing us to
analyze the complete convergence behavior of all \reviewerB{methods}. For the
forecasting methods, we employ a memory of $\memory = 8$ and restricted-state
dimension $\nrestrict = 8$ (i.e., we forecast only the first 8 POD modes).

\begin{figure}[htb]
\centering 
\subfigure[$\paramOnlinei{1} = (1.6603, 0.0229)$]{
\includegraphics[width=0.48\textwidth]{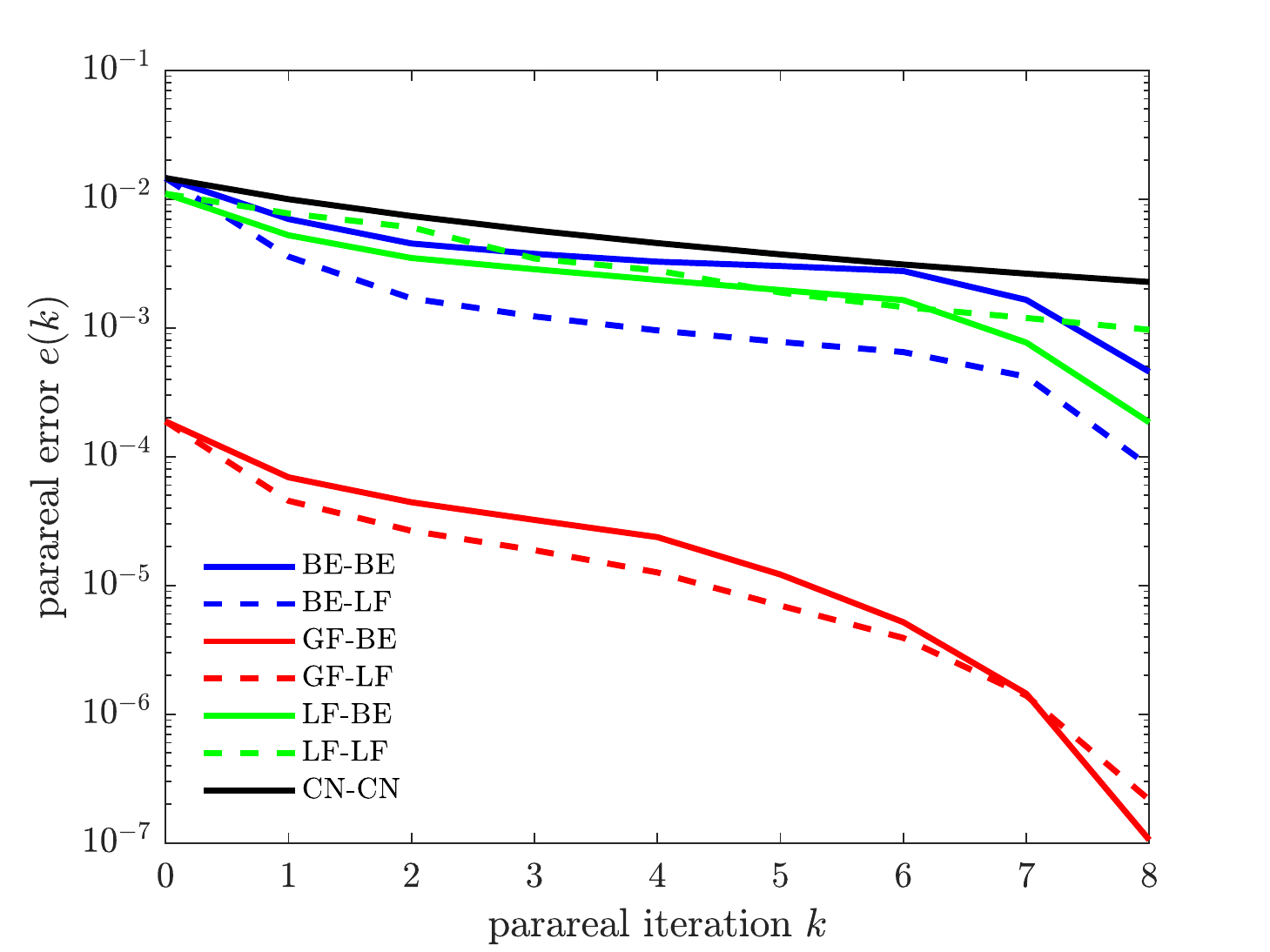} 
\label{fig:methodComp_plot_1}
}
\subfigure[$\paramOnlinei{2} = (1.5025, 0.0201)$]{
\includegraphics[width=0.48\textwidth]{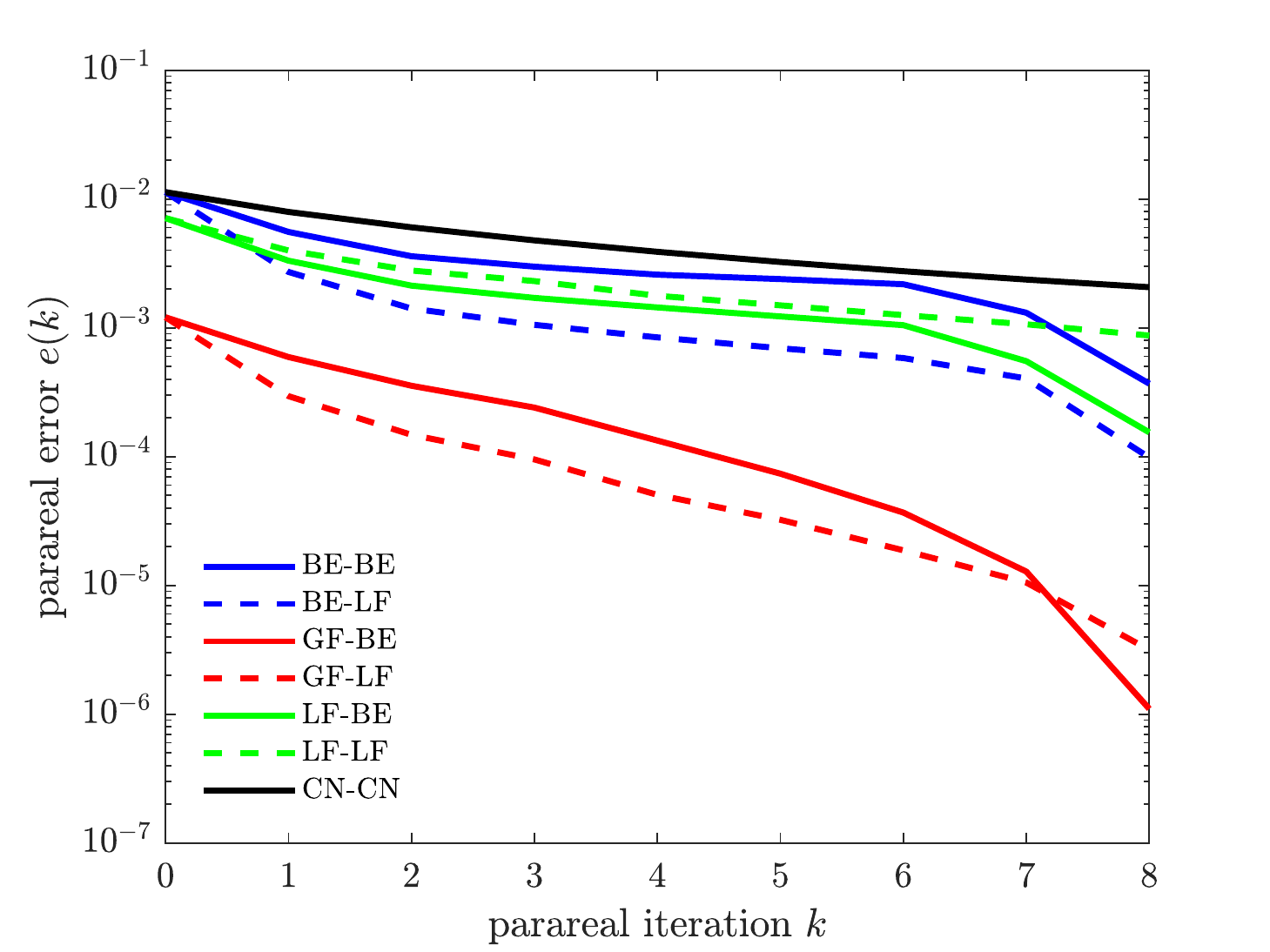} 
\label{fig:methodComp_plot_2}
}
\caption{
	\textit{Comparison of initialization and coarse-propagation methods}. Convergence of
	six methods. 
}
\label{fig:methodComp_plot}
\end{figure}

Figure \ref{fig:methodComp_plot}
reports these results, where the time-parallel error at parareal iteration
$\pararealit$ is defined as 
$
\timeParallelError(\pararealit)\defeq
\reviewerB{\max_{\timestepdummy\in\{\pararealit+1,\ldots,\nIntervals-1\}}
	\|
\fnprop{\tptcrs{\timestepdummy}}{\tptcrs{\timestepdummy
-1}}{\apxsolFOM{\timestepdummy -1}{\pararealit  }}
	-\apxsolFOM{\timestepdummy}{\pararealit
}\|/\|\fnprop{\tptcrs{\timestepdummy }}{\tptcrs{\timestepdummy
-1}}{\apxsolFOM{\timestepdummy-1 }{\pararealit  }}\|
,}
$
\reviewerB{which is a measure of the normalized residual that the parareal
	method is aiming to set to zero
	\cite{gander2007analysis}.\footnote{\reviewerB{In
	Algorithm \ref{alg:parareal}, this corresponds to the tolerance that appears in Step
\ref{step:termination}.}}}
\ourReReading{The figure} 
highlights two important trends. First,
the results empirically support the theoretical result discussed in Remark
\ref{remark:initialSeedStable}: namely, the global-forecast initialization
exhibits superior stability properties to the local-forecast initialization.
In both online parameter instances, global-forecast initialization produces a
very small initial error, while local-forecast initialization
produces a larger initial error despite its use of the same time-evolution
data; backward-Euler \reviewerB{and Crank--Nicolson}
initialization produces a slightly smaller initial error than the local
forecast.
Second, note that the local-forecast propagator 
outperforms the
backward-Euler coarse propagator when either the backward-Euler or
global-forecasting initializations are employed.

\begin{figure}[htbp!] 
%
\centering 
\subfigure[\BELabel-\BELabel]{
\includegraphics[width=0.48\textwidth]{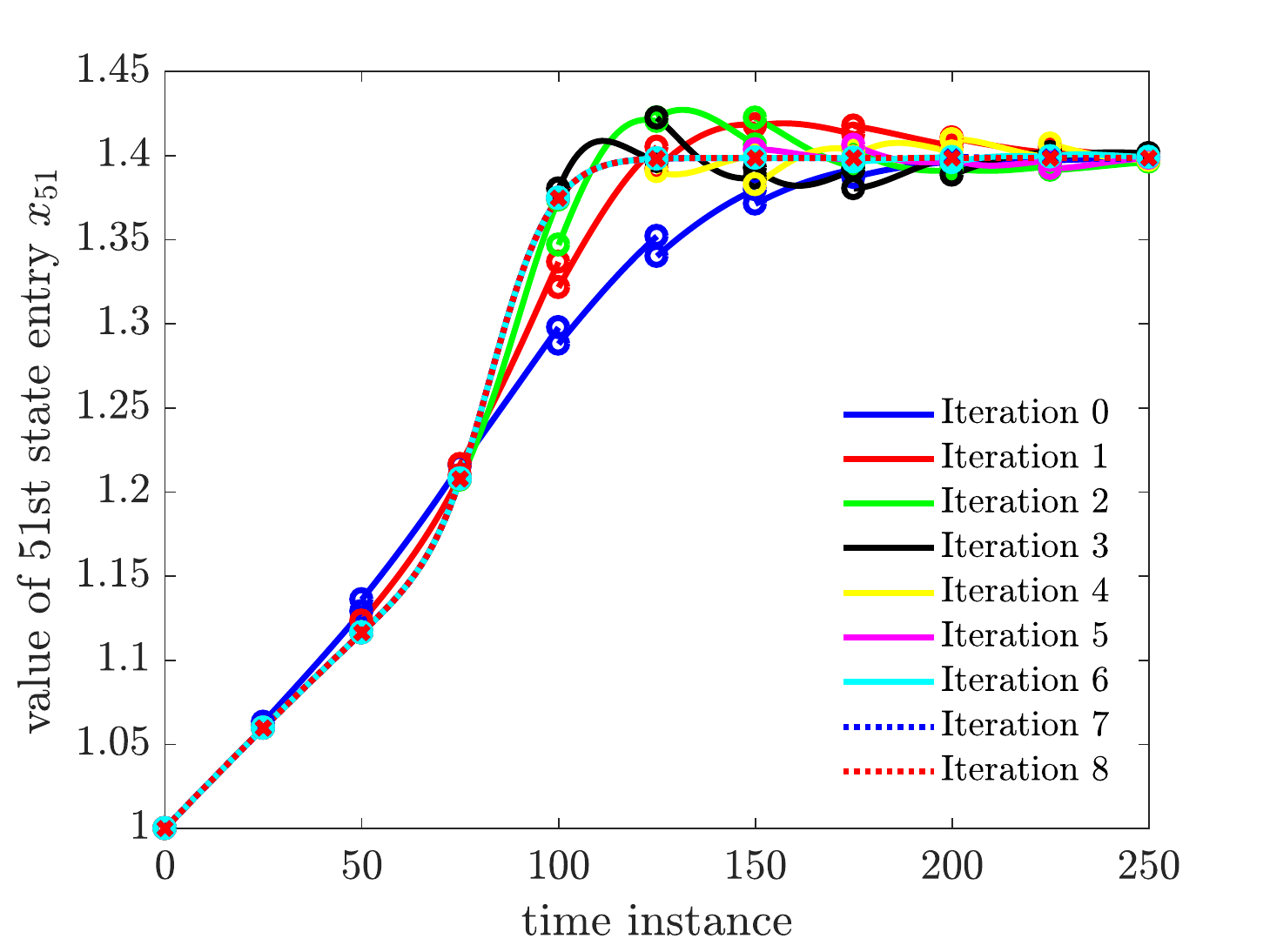} 
\label{fig:jumpsParareal_plot_1}
}
\subfigure[\BELabel-\localForeLabel]{
\includegraphics[width=0.48\textwidth]{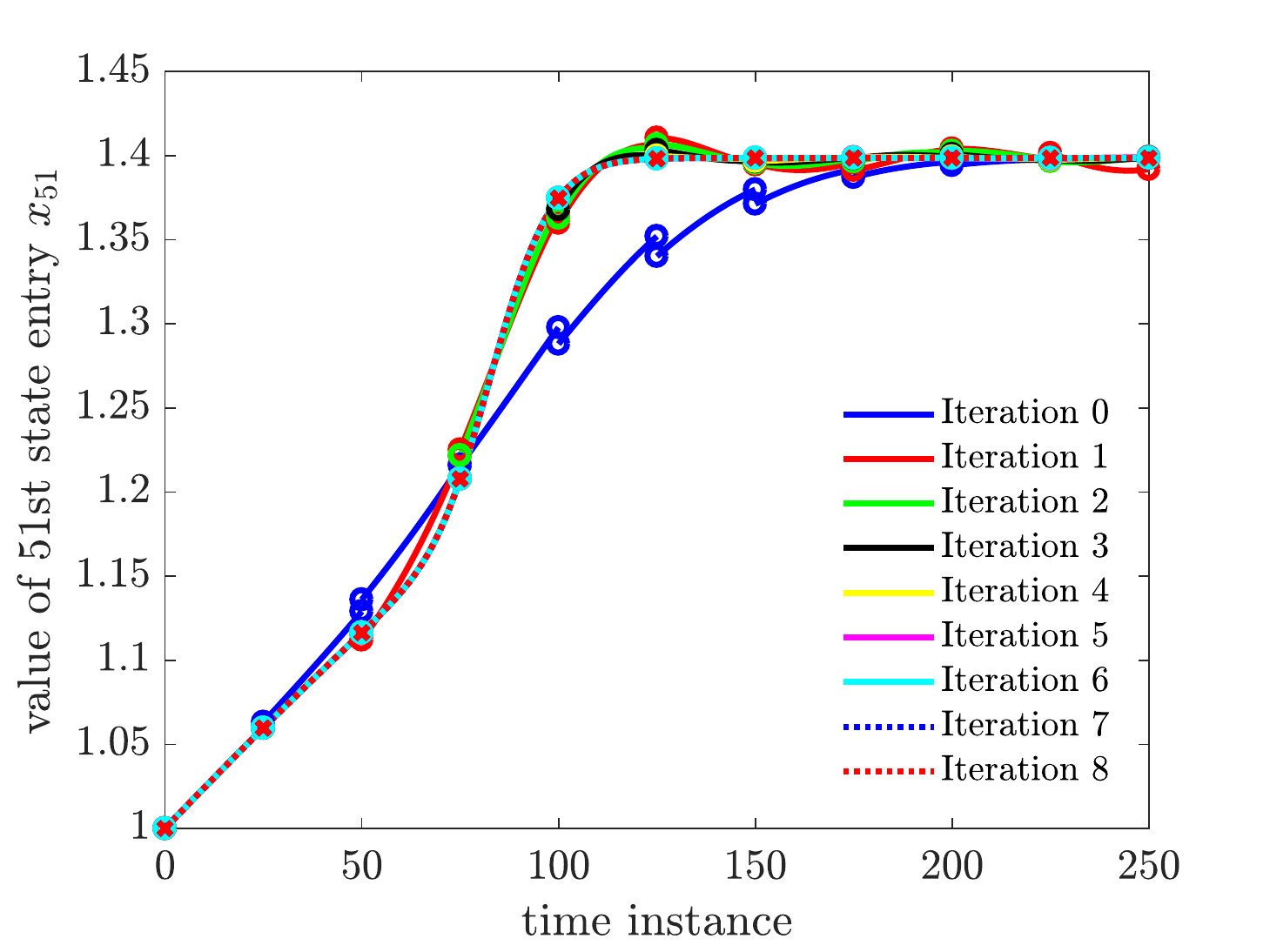} 
\label{fig:jumpsParareal_plot_2}
}
\subfigure[\localForeLabel-\BELabel]{
\includegraphics[width=0.48\textwidth]{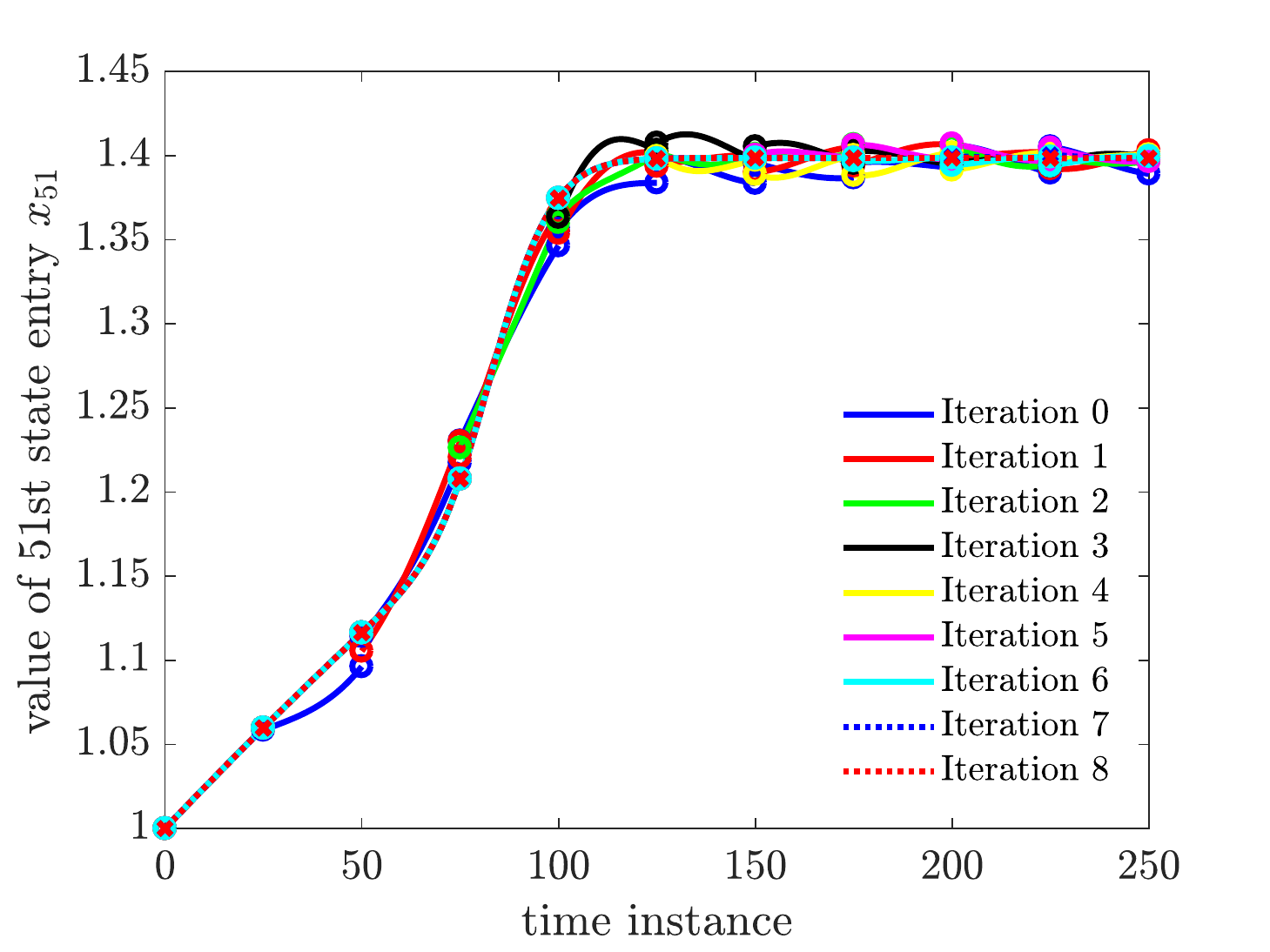} 
\label{fig:jumpsParareal_plot_3}
}
\subfigure[\localForeLabel-\localForeLabel]{
	\includegraphics[width=0.48\textwidth]{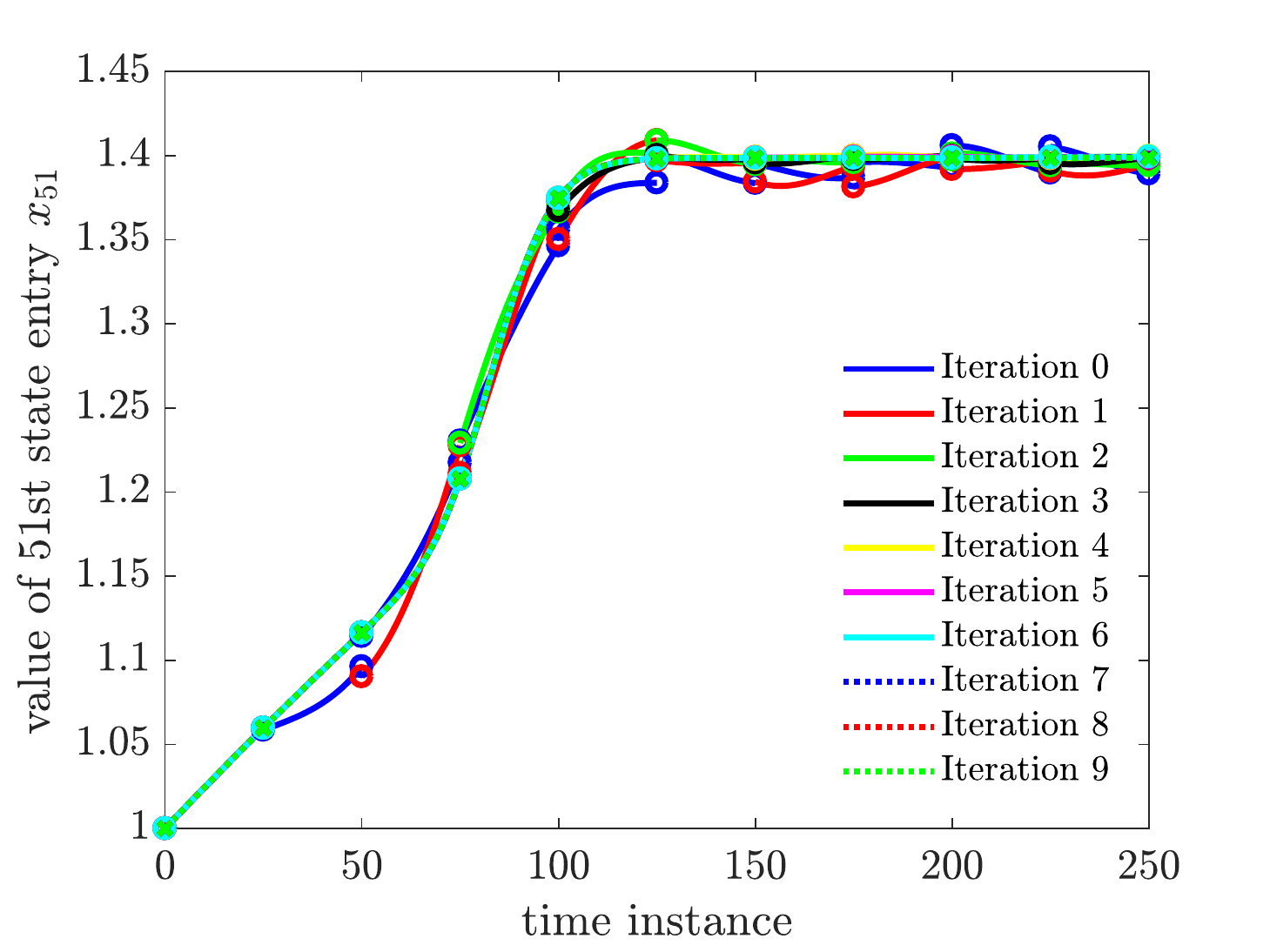} 
\label{fig:jumpsParareal_plot_4}
}
\subfigure[\globalForeLabel-\BELabel]{
\includegraphics[width=0.48\textwidth]{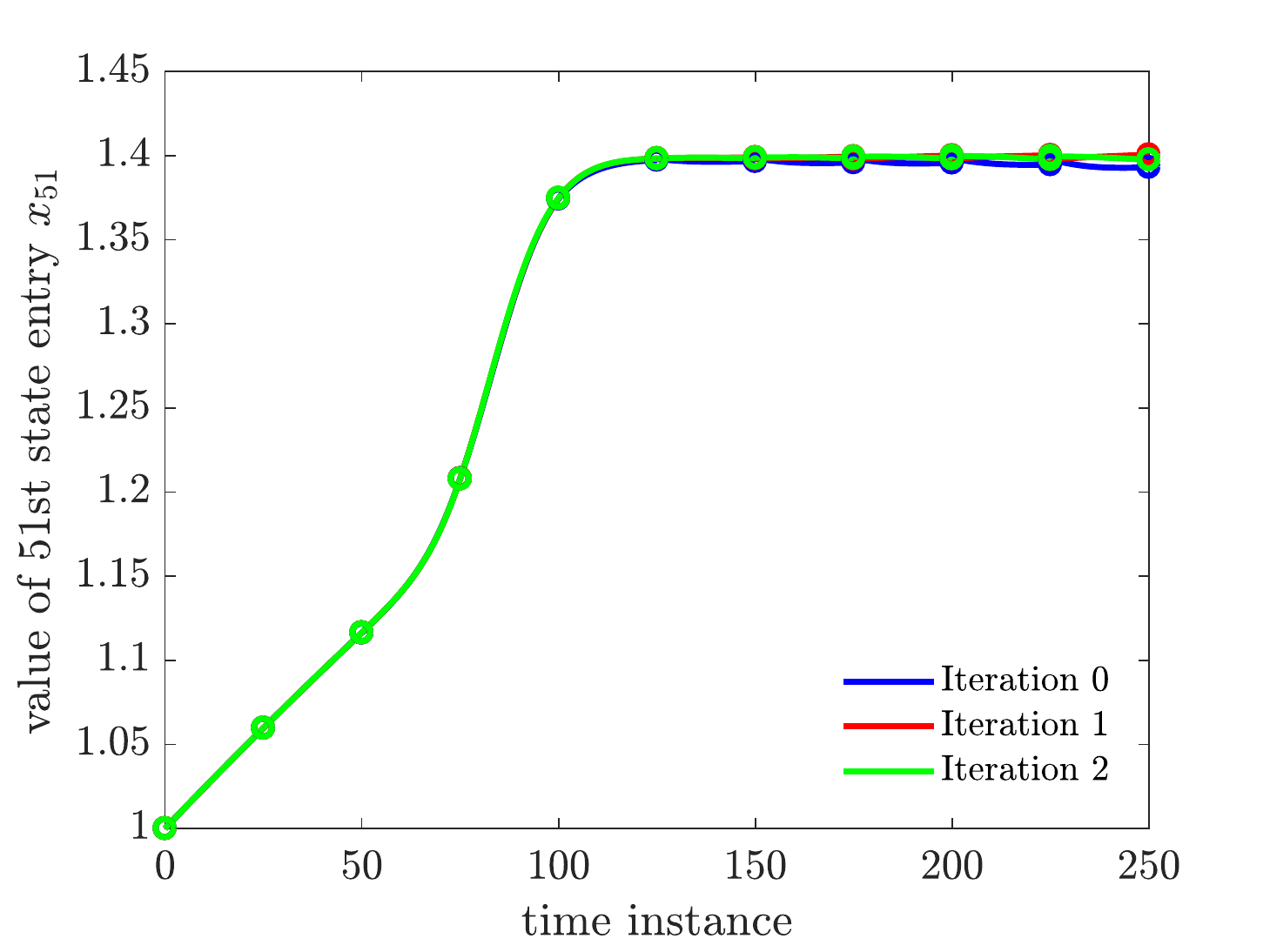} 
\label{fig:jumpsParareal_plot_5}
}
\subfigure[\globalForeLabel-\localForeLabel]{
\includegraphics[width=0.48\textwidth]{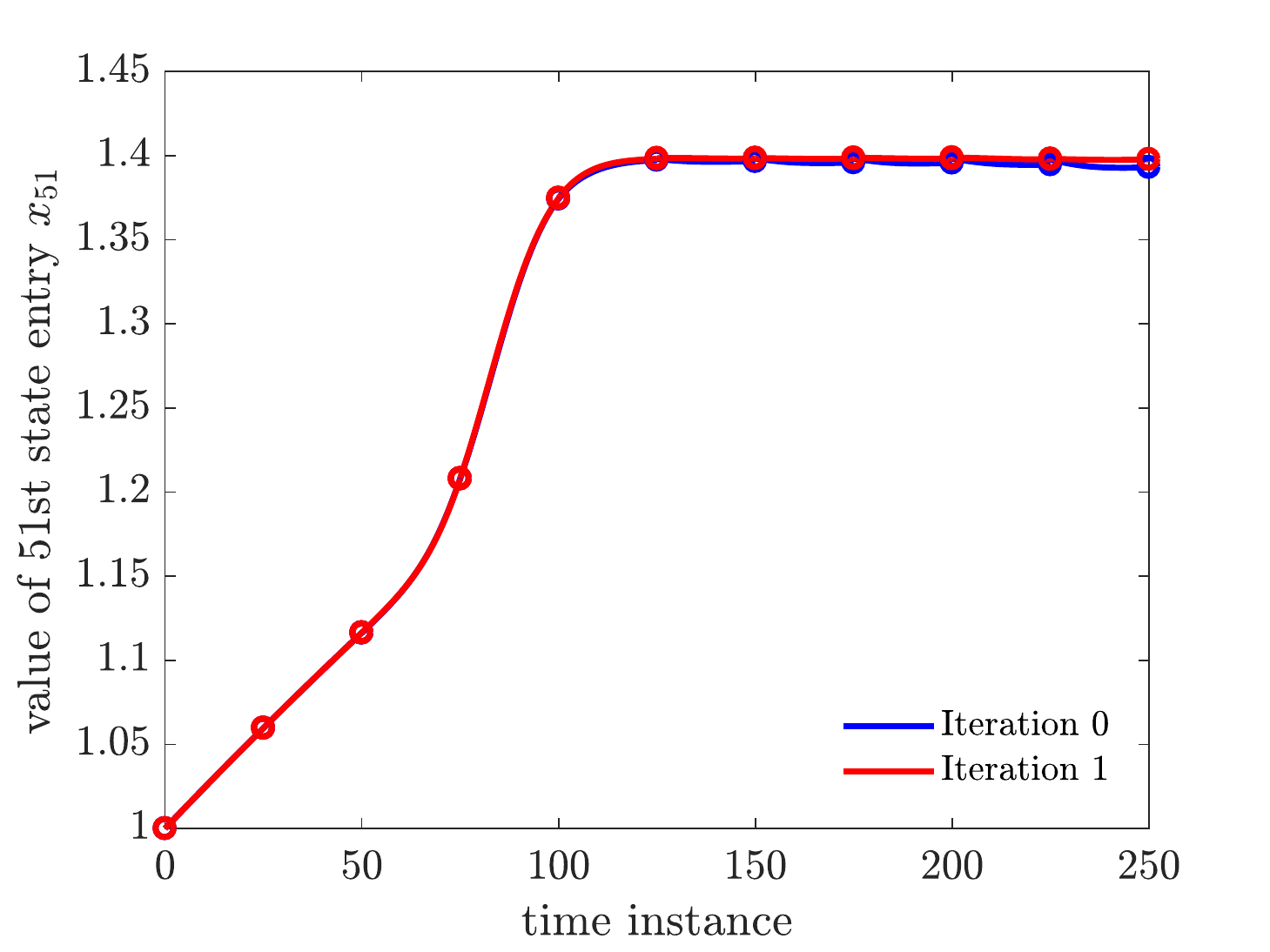} 
\label{fig:jumpsParareal_plot_6}
}
\caption{
	\textit{Comparison of initialization and coarse-propagation methods}.
	Convergence of the 51st entry of the state vector for $\paramOnlinei{2} =
	(1.5025, 0.0201)$ for six methods.\\
}
\label{fig:jumpsParareal_plot}
\end{figure}

\begin{figure}[htbp!] 
\centering 
\subfigure[\BELabel-\BELabel]{
\includegraphics[width=0.48\textwidth]{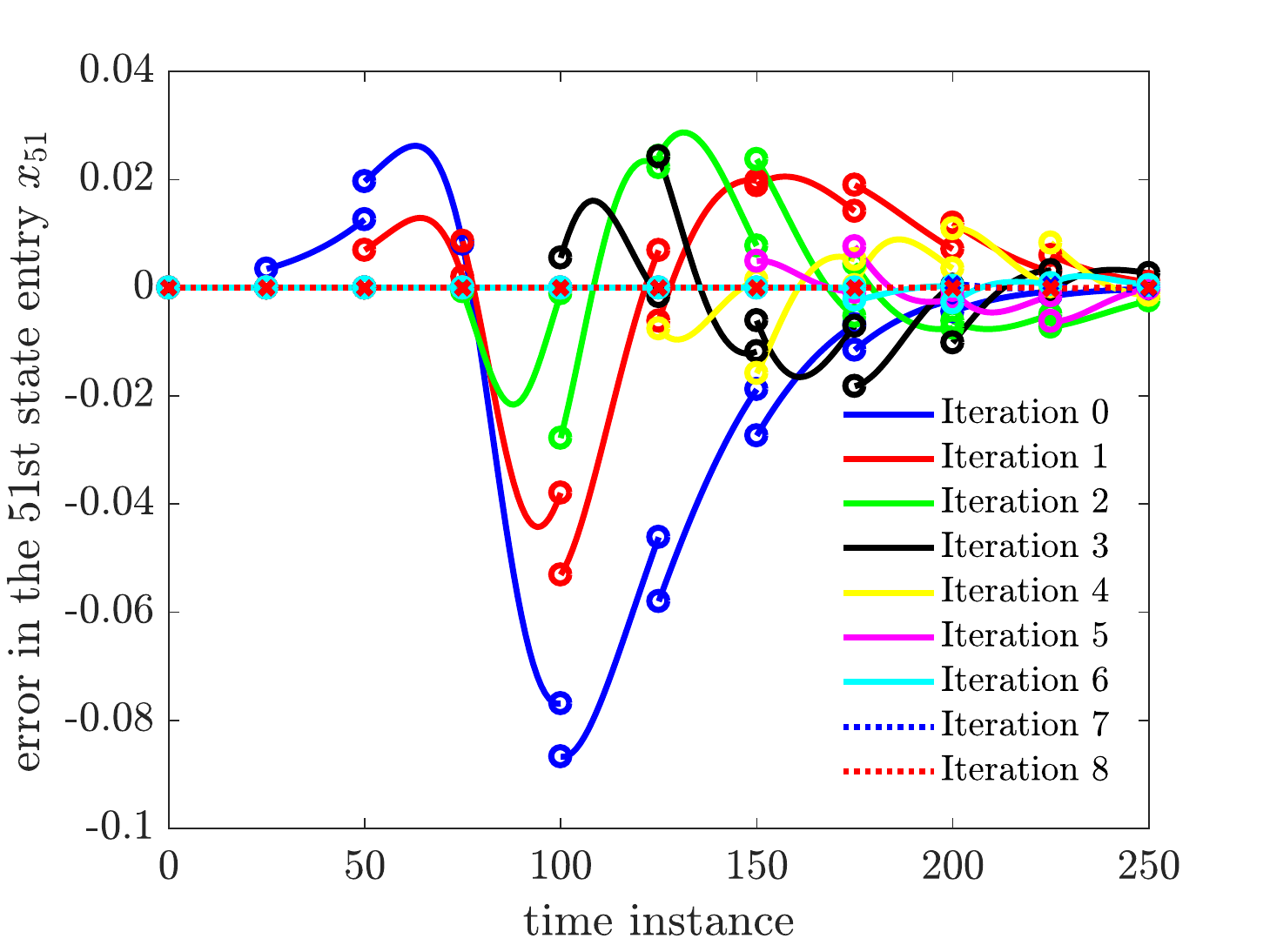} 
\label{fig:ErrorsParareal_plot_1}
}
\subfigure[\BELabel-\localForeLabel]{
\includegraphics[width=0.48\textwidth]{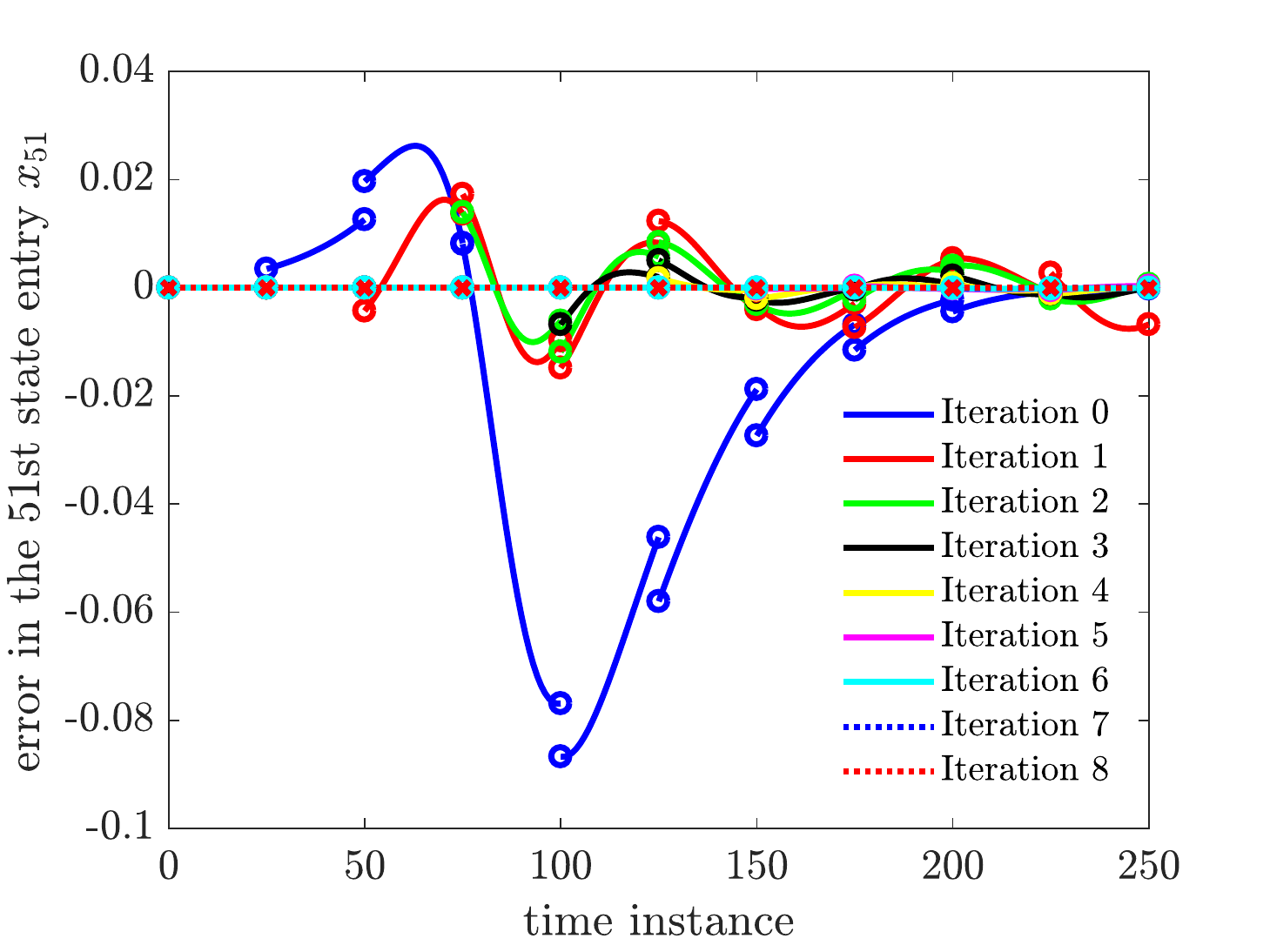} 
\label{fig:ErrorsParareal_plot_2}
}
\subfigure[\localForeLabel-\BELabel]{
\includegraphics[width=0.48\textwidth]{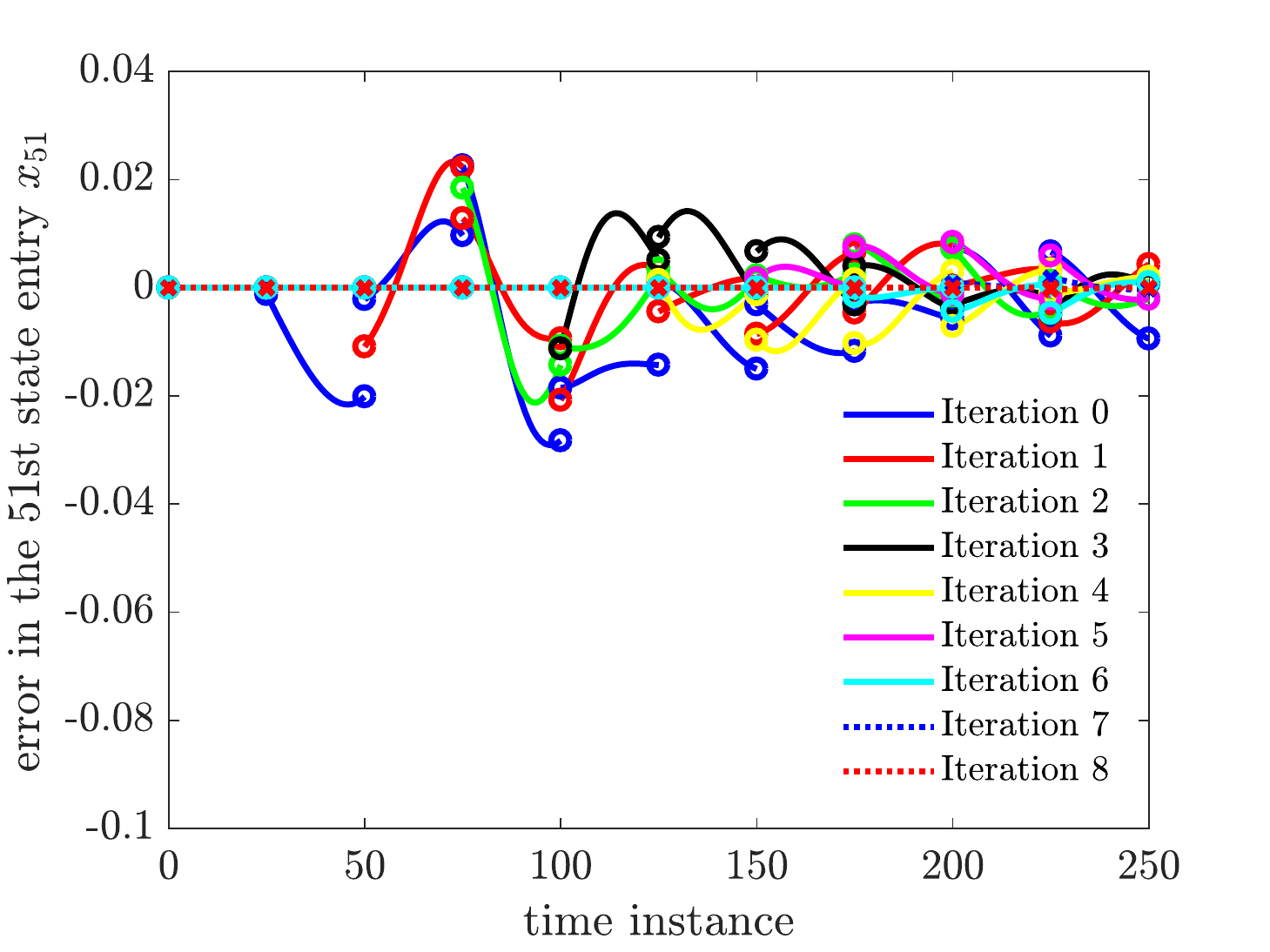} 
\label{fig:ErrorsParareal_plot_3}
}
\subfigure[\localForeLabel-\localForeLabel]{
\includegraphics[width=0.48\textwidth]{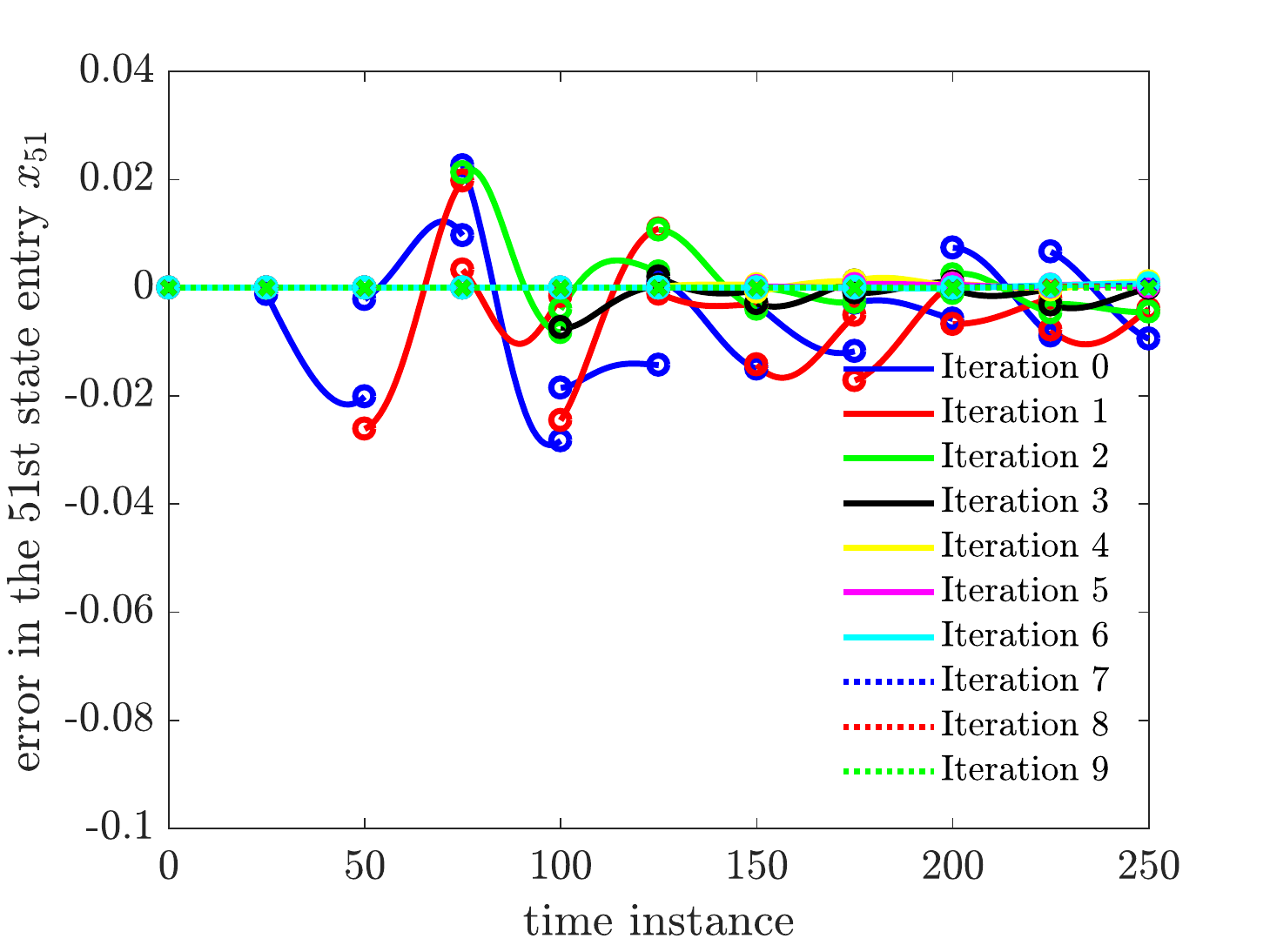} 
\label{fig:ErrorsParareal_plot_4}
}
\subfigure[\globalForeLabel-\BELabel]{
\includegraphics[width=0.48\textwidth]{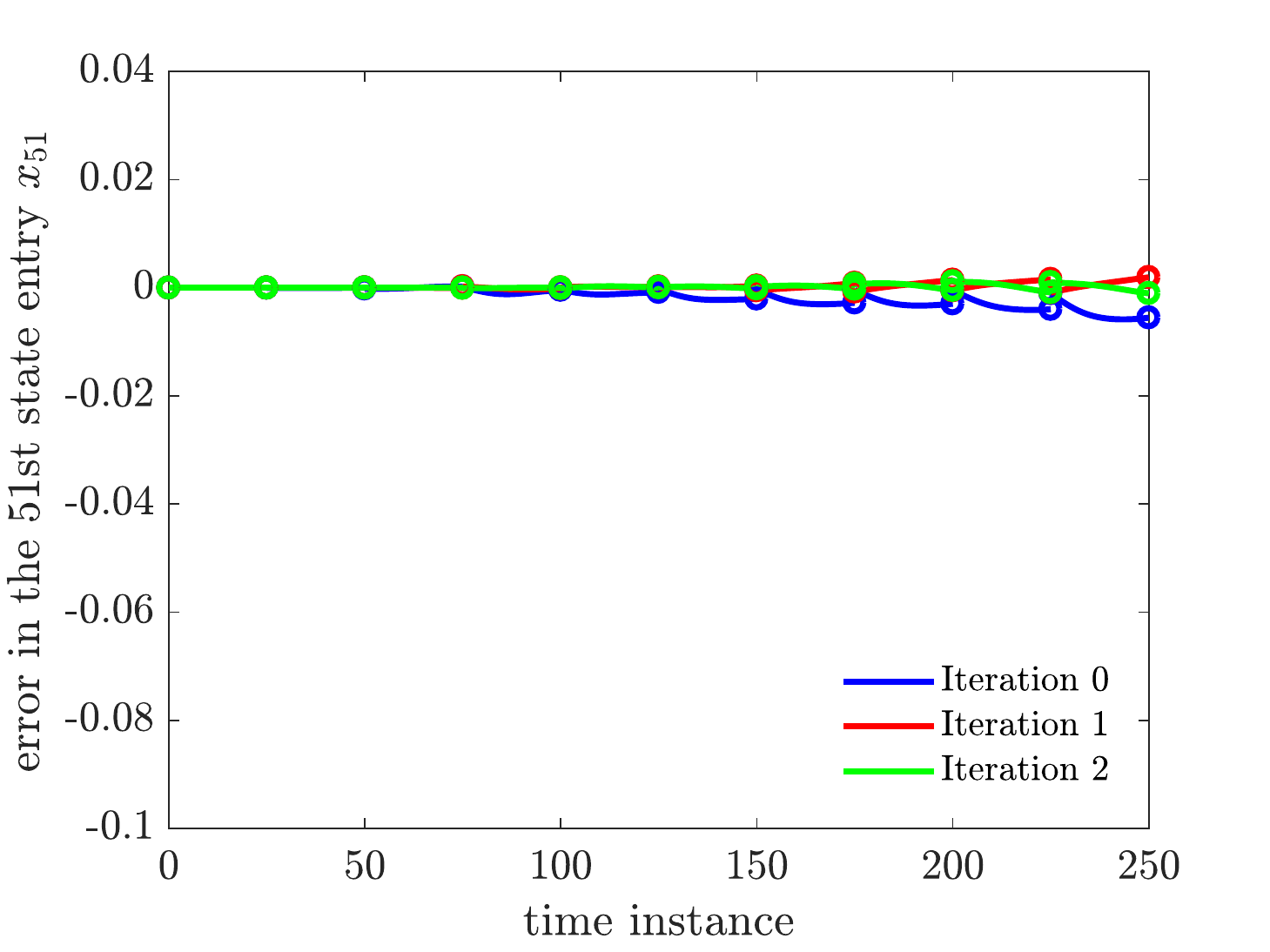} 
\label{fig:ErrorsParareal_plot_5}
}
\subfigure[\globalForeLabel-\localForeLabel]{
\includegraphics[width=0.48\textwidth]{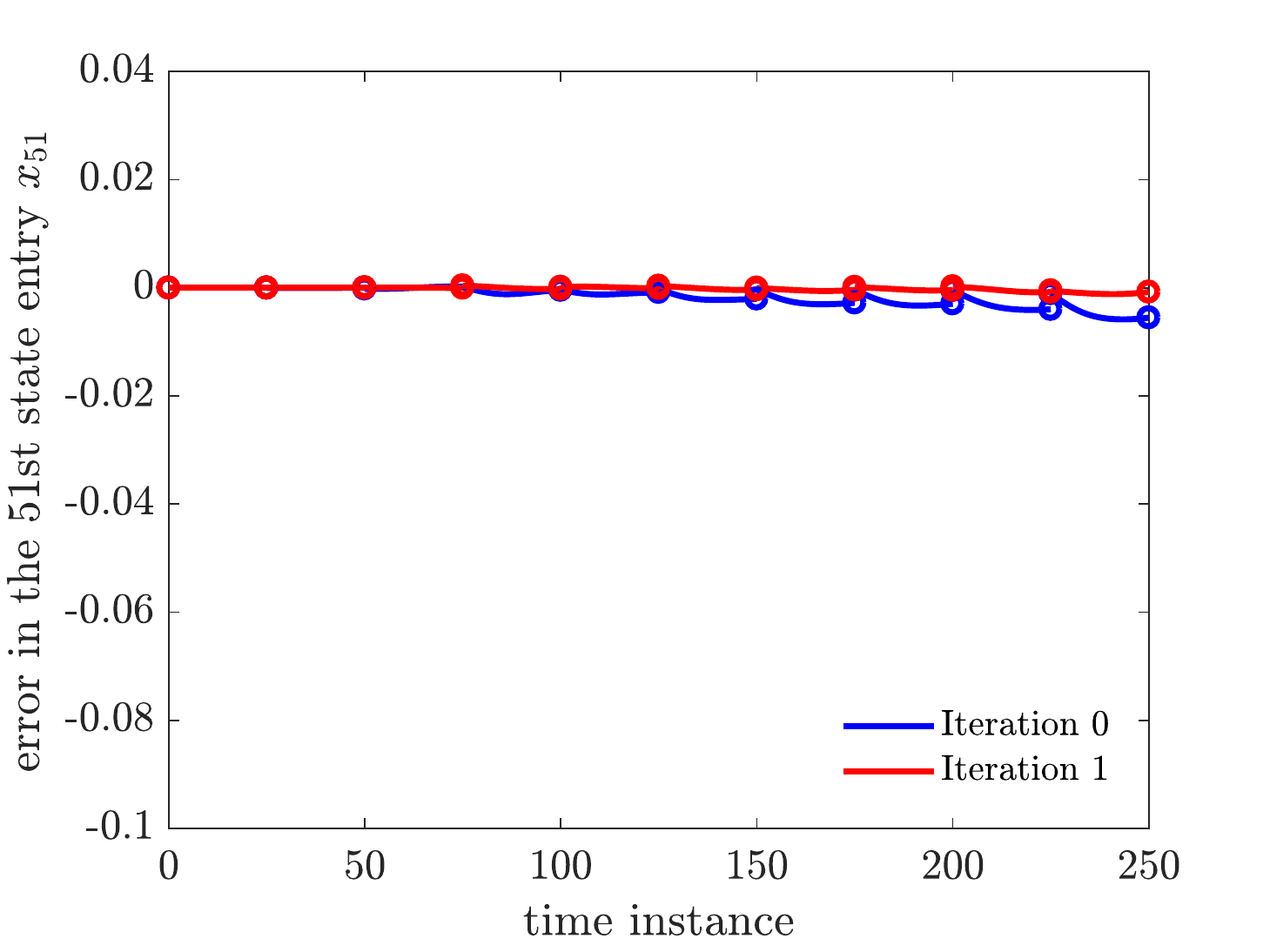} 
\label{fig:ErrorsParareal_plot_6}
}
\caption{
	\textit{Comparison of initialization and coarse-propagation methods}.
	Converence of the error in the 51st entry of the state vector (with respect
	to the serial solution) for $\paramOnlinei{2} =
	(1.5025, 0.0201)$ for six methods. 
}
\label{fig:ErrorsParareal_plot}
\end{figure}

To gain additional insight into the convergence properties of the methods,
Figure \ref{fig:jumpsParareal_plot} reports the convergence of the 51st entry
of the state vector over parareal iterations for online parameter
instance $\paramOnlinei{2}$\reviewerB{, and Figure
\ref{fig:ErrorsParareal_plot} reports convergence of the error in this
quantity}. These results highlight the two trends mentioned
above; specifically, global-forecast initialization leads to a nearly exact
initial solution, local-forecast initialization leads to a very poor initial
solution, and local-forecast coarse propagation reduces errors more quickly
than backward-Euler coarse propagation, even when backward-Euler
initialization is employed\reviewerB{; these plots do not include the CN-CN
results, as they are very similar to the BE-BE results}.

In the remainder of the numerical experiments, we limit our focus to the
typical parareal \reviewerB{methods} \BELabel-\BELabel\ and
\reviewerB{\CNLabel-\CNLabel, as well as} the most promising proposed
data-driven method \globalForeLabel-\localForeLabel. 

\subsection{Ideal case}
\label{ssec:reproduction_case}
This section assesses performance of the method under the `ideal case', i.e.,
when Assumptions \ref{ass:subspace}--\ref{ass:isomorphic} are satisfied as
discussed in Sections \ref{sec:idealError} and \ref{sec:idealSpeedups}.
Here,
we ensure these conditions are met by repeating the training for each online
point (i.e., $\ntrain = 1$ with both the training point set equal to the
online point) and employing $\nrestrict=\reddim$.  Recall that under these conditions, the coarse propagator is
exact (Theorem \ref{thm:exactCoarse}), and the
\globalForeLabel-\localForeLabel\ method should converge after parareal
initialization (hence require \ourReReading{$K=0$} parareal iterations) 
and produce speedups given by Eq.~\eqref{eq:idealSpeedupGlobal}
(Theorem \eqref{thm:idealSpeedupGlobal}).  \reviewerB{Note that these
conditions are `ideal' for the proposed methodology, but not for typical
parareal methods \BELabel-\BELabel\ or \CNLabel-\CNLabel.} We assess memories of $\memory = 1,2,4,6$ and employ a termination
tolerance of $\jumpTolerance = \ourReReading{5\times 10^{-4}}$ in this section only.

 In the remaining experiments, we report the theoretical
 speedups derived in Section \ref{sec:speedup} due to the lack of
 reliability in timings obtained with our Matlab
 implementation.\footnote{Future work will entail implementation in a
 `production' computational-mechanics code
 and assessment of the method in a parallel computing environment.} Here, the speedup for method
\globalForeLabel-\localForeLabel\ 
is provided by
Eq.~\eqref{eq:speedupGlobal}, and the speedup for \reviewerB{methods} \BELabel-\BELabel\
\reviewerB{and \CNLabel-\CNLabel}\footnote{\reviewerB{Because each method is
		characterized by only one implicit stage, we assume that the cost
of Crank--Nicolson is the same as that of backward Euler; the additional
explicit stage for Crank--Nicolson introduces negligible additional cost.}}\
 \reviewerB{are} provided by the following theorem, whose proof can be found in Appendix~\ref{sec:proofs}.
\begin{theorem}[\textit{Speedup}: fine propagator as coarse
propagator]\label{thm:finePropGen}
If the same time integrator is used for both the coarse and fine
propagator and Assumption \ref{ass:integrateCostDominant}
holds, then the parareal method realizes a speedup of
\begin{equation}\label{eq:finePararealSpeedup}
\speedupFine{\pararealItConverge} \defeq 
\frac{\nptfn}{(\nptcrs+ \nptfnArg{\pararealit}-\frac{1}{2}\pararealItConverge)(\pararealItConverge+1) 
}
.
\end{equation}
\end{theorem}

\begin{figure}[htbp!] 
\centering 
\subfigure[$\paramOnlinei{1} = (1.6603, 0.0229)$]{
\includegraphics[width=0.48\textwidth]{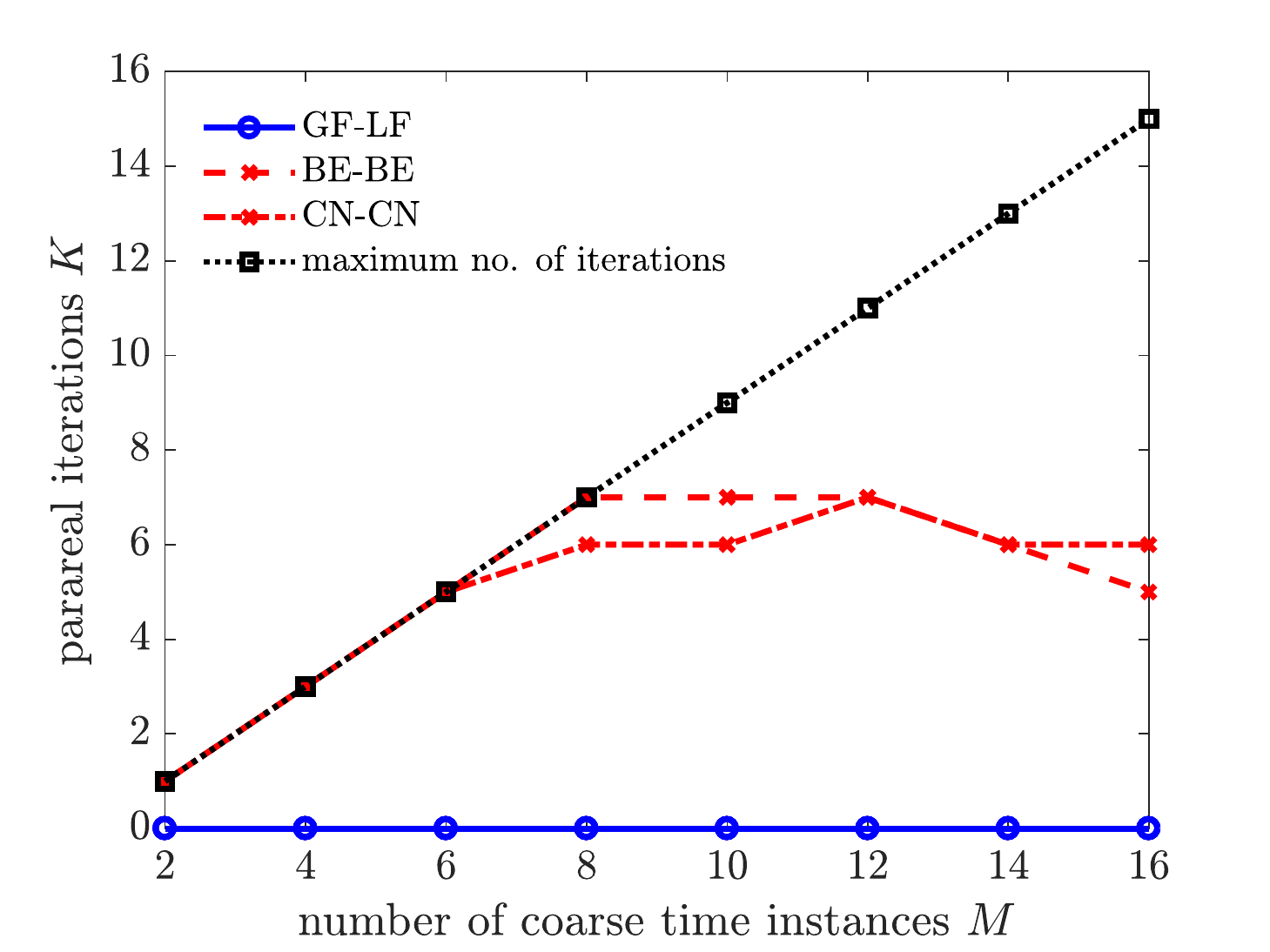} 
\label{fig:repro_iterationVsCPU_plot_1}
}
\subfigure[$\paramOnlinei{2} = (1.5025, 0.0201)$]{
\includegraphics[width=0.48\textwidth]{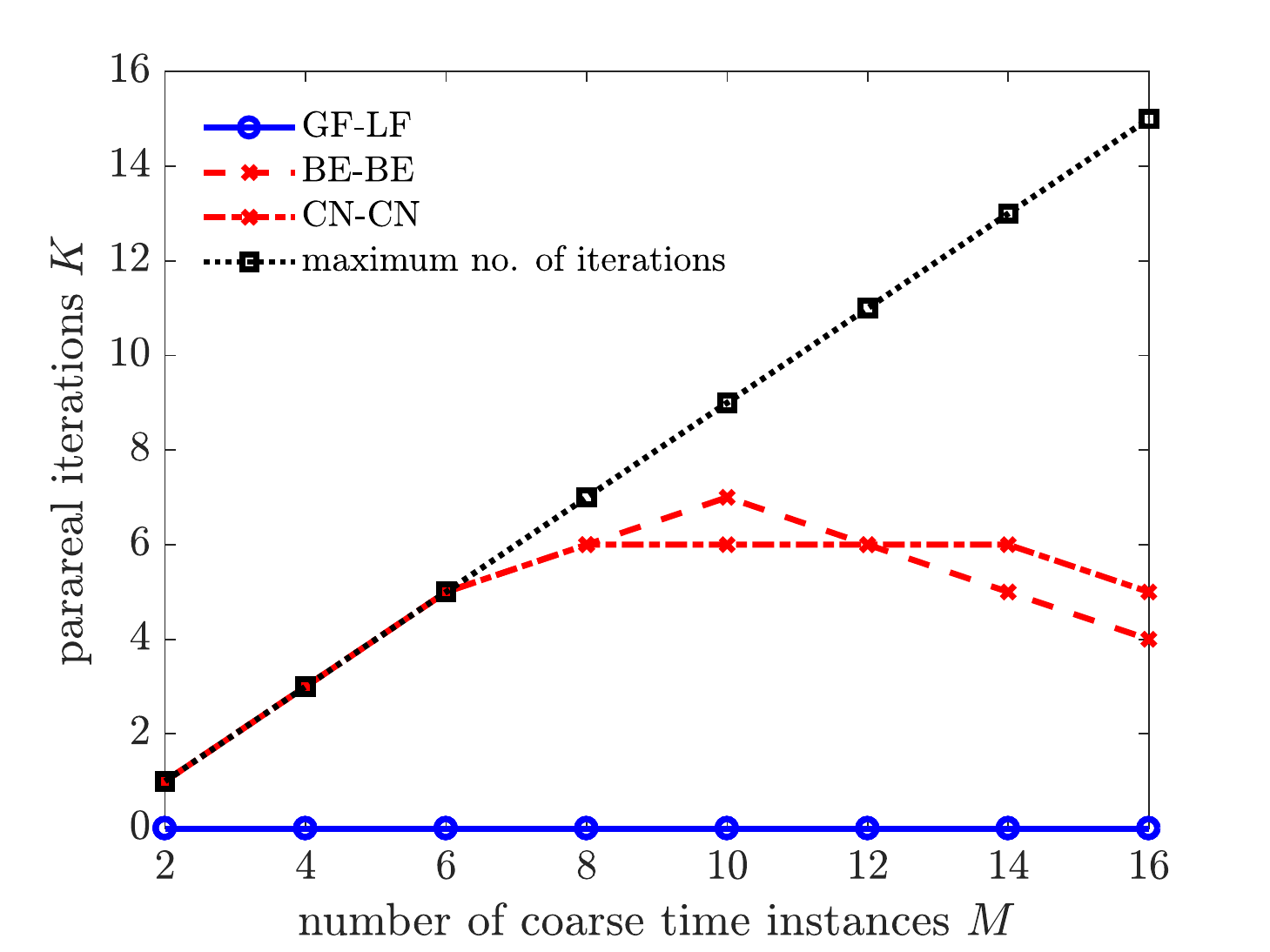} 
\label{fig:repro_iterationVsCPU_plot_2}
}
\subfigure[$\paramOnlinei{1} = (1.6603, 0.0229)$]{
\includegraphics[width=0.48\textwidth]{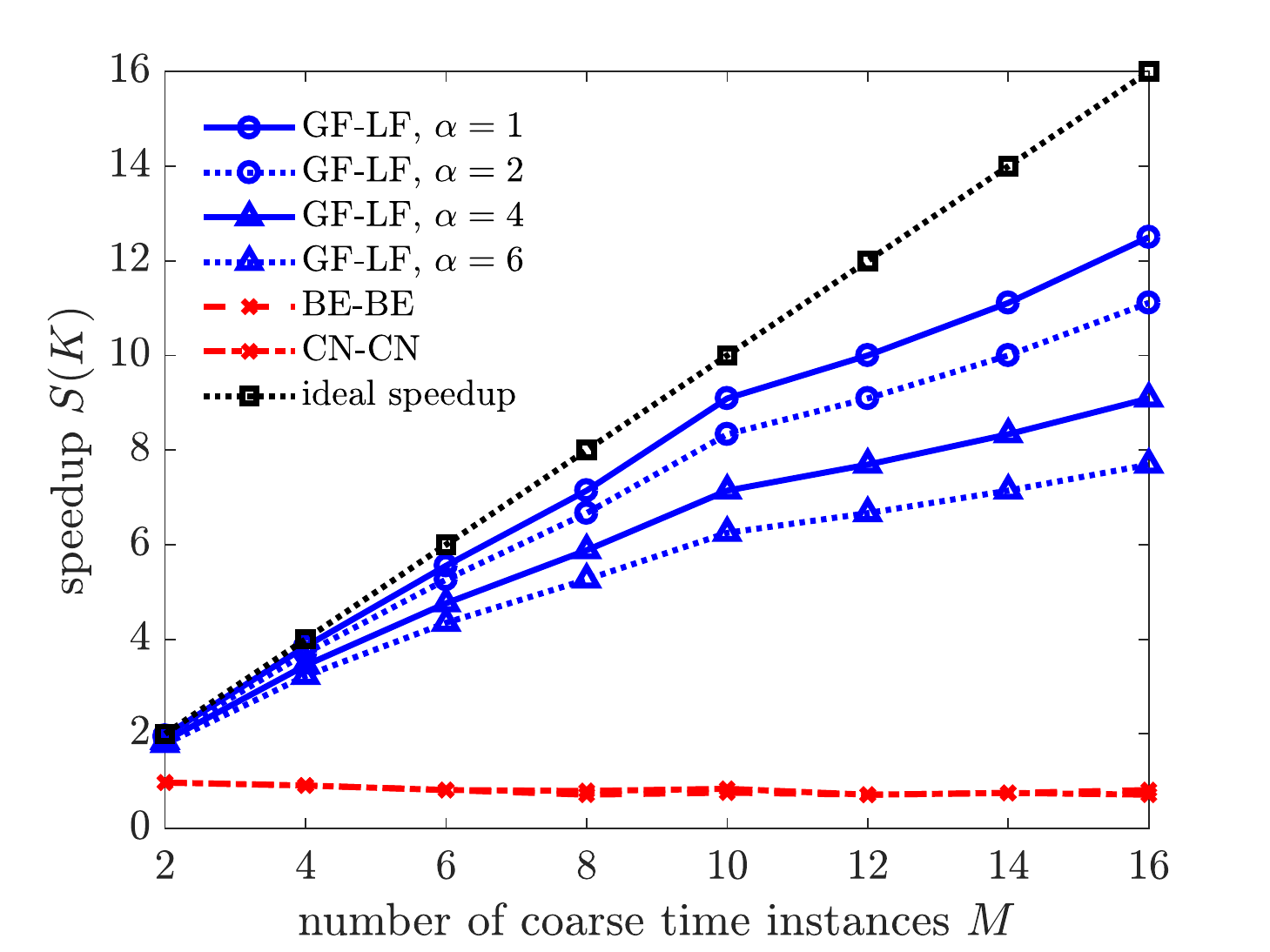} 
\label{fig:repro_speedupVsCPU_ideal_plot_1}
}
\subfigure[$\paramOnlinei{2} = (1.5025, 0.0201)$]{
\includegraphics[width=0.48\textwidth]{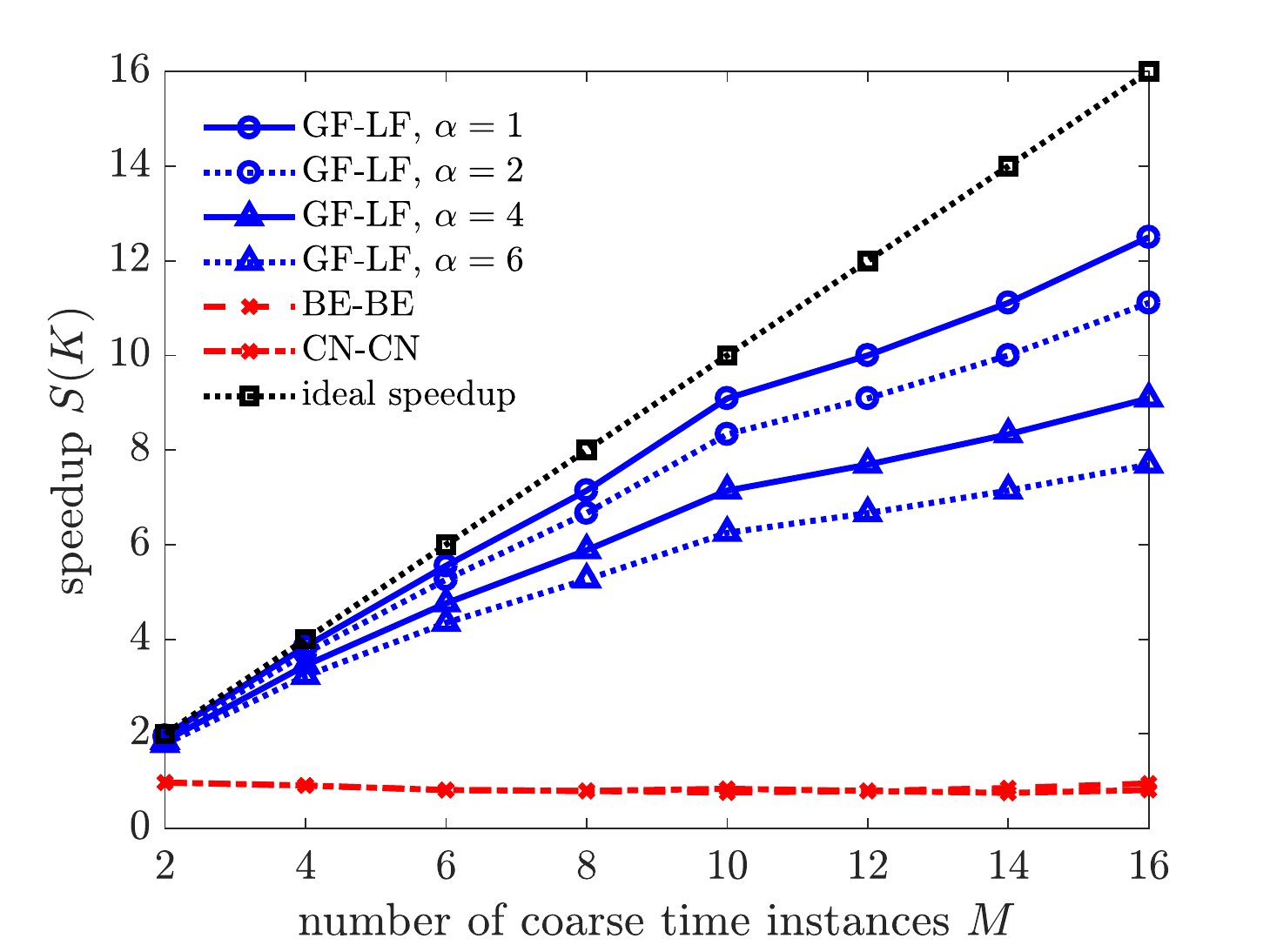} 
\label{fig:repro_speedupVsCPU_ideal_plot_2}
}
\subfigure[$\paramOnlinei{1} = (1.6603, 0.0229)$]{
\includegraphics[width=0.48\textwidth]{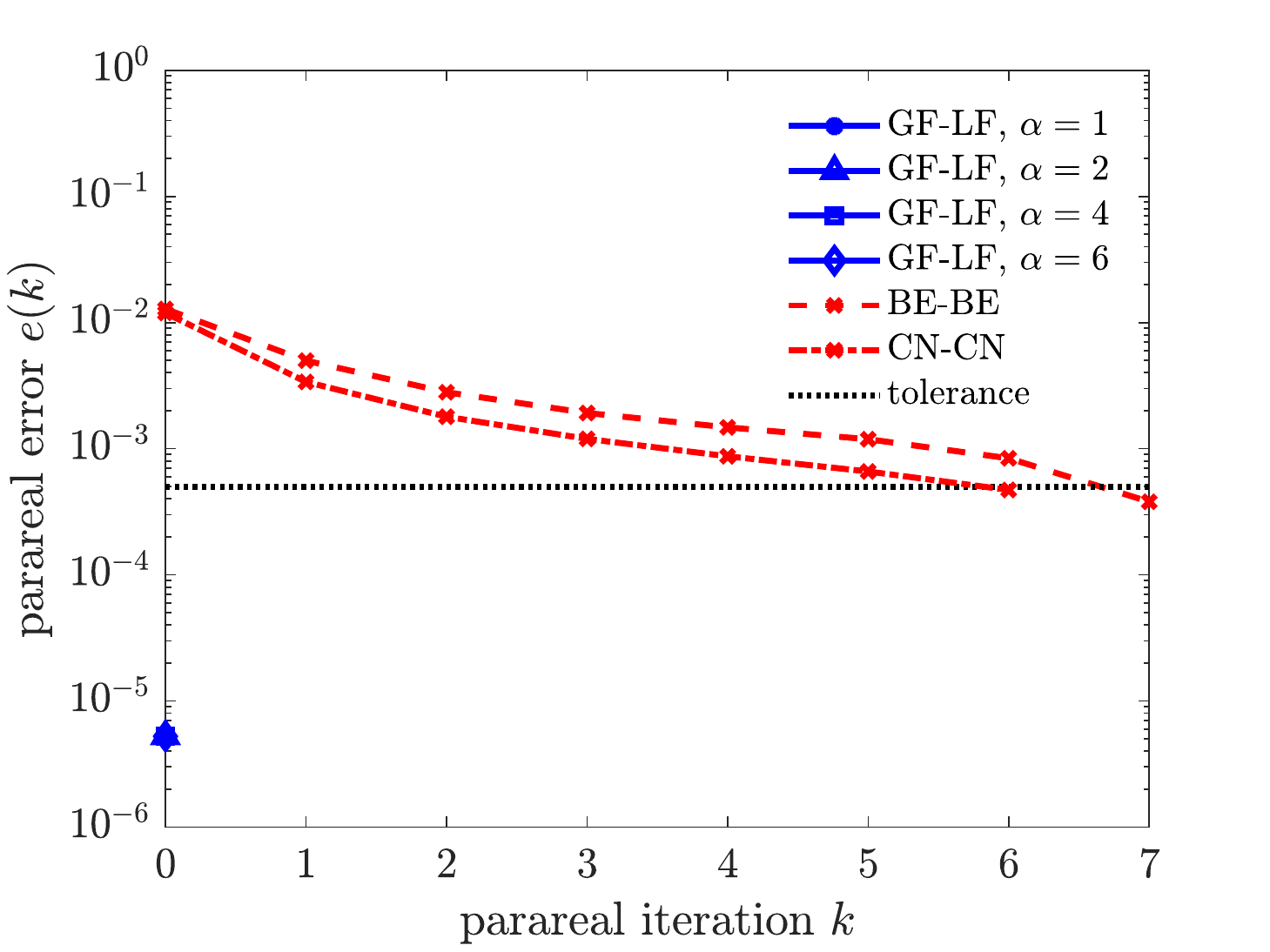} 
\label{fig:repro_timeParallelError_plot_1}
}
\subfigure[$\paramOnlinei{2} = (1.5025, 0.0201)$]{
\includegraphics[width=0.48\textwidth]{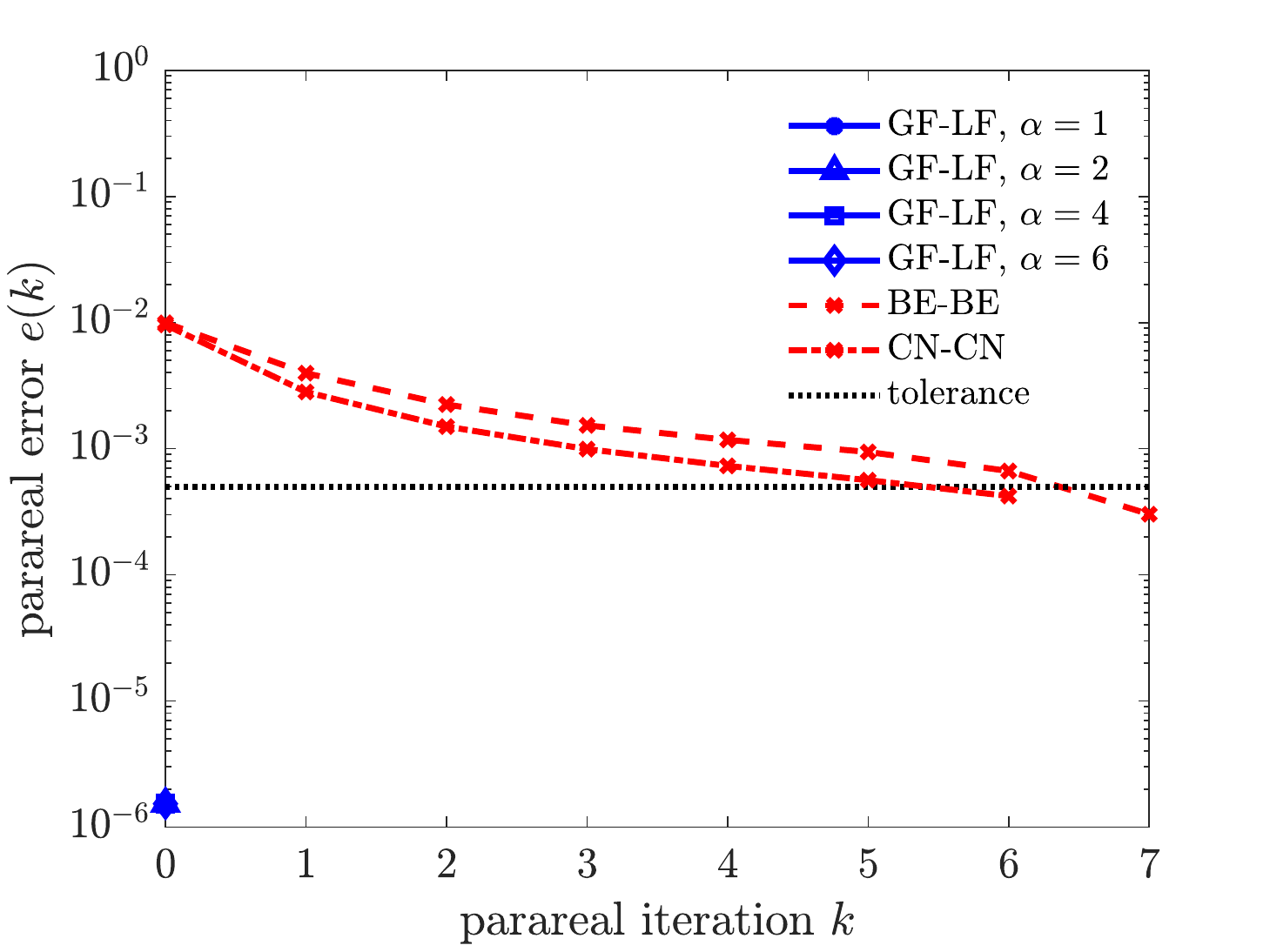} 
\label{fig:repro_timeParallelError_plot_2}
}
\caption{\textit{Ideal case}. Number of parareal iterations $\pararealItConverge$ required for
convergence, theoretical speedups computed via
Eq.~\eqref{eq:finePararealSpeedup} for \BELabel-\BELabel\ and
Eq.~\eqref{eq:speedupGlobal} for \globalForeLabel-\localForeLabel, and
convergence plots. \ourReReading{Figures (e) and (f) associate with $\nptcrs =
10$}. 
}
\label{fig:ideal_case}
\end{figure}

Figures \ref{fig:repro_iterationVsCPU_plot_1}--\ref{fig:repro_iterationVsCPU_plot_2} report the number of parareal
iterations required for convergence when the number of coarse time instances
$\nptcrs$ increases \reviewerB{(and the coarse time step $\tsscrs =
\finalT/\nptcrs$ undergoes an attendant decrease)}. As expected, in all cases, the proposed
\globalForeLabel-\localForeLabel\ method converges in the minimum number of
parareal iterations (i.e., $\pararealItConverge=0$). In contrast, the
\BELabel-\BELabel\ \reviewerB{and \CNLabel-\CNLabel\ methods converge} in the \textit{worst-case} number of
iterations (i.e., $\pararealItConverge=\nptcrs-1$) for $\nptcrs\leq 6$ in both
cases; this occurs because these cases correspond to relatively large coarse
time steps $\tsscrs$, which degrades the accuracy of the backward-Euler
\reviewerB{and Crank--Nicolson schemes}. The number of parareal iterations needed for convergence in the
\BELabel-\BELabel\ \reviewerB{and \CNLabel-\CNLabel\ cases} decreases as the number of coarse time instances
$\nptcrs$ increases; this can be attributed to the decreasing coarse time
step $\tsscrs$, which improves the accuracy of the \reviewerB{time
integrators}.

Figures \ref{fig:repro_speedupVsCPU_ideal_plot_1}--\ref{fig:repro_speedupVsCPU_ideal_plot_2} report the
theoretical speedups of \reviewerB{these} methods under these ideal conditions. Here, the
reported values correspond to Eq.~\eqref{eq:finePararealSpeedup} for the
\BELabel-\BELabel\ \reviewerB{and \CNLabel-\CNLabel\ methods} and Eq.~\eqref{eq:speedupGlobal} for the
\globalForeLabel-\localForeLabel\ method. As expected, the proposed technique
yields near-ideal theoretical speedups, while the typical
\reviewerB{approaches produce}
modest speedups due to \reviewerB{their} slow convergence on this problem. Note that
increasing the memory degrades speedup in this case, as all values for the
memory ensure an exact initial solution in the ideal case; thus, employing a
small memory does not degrade convergence here.

Finally, Figures \ref{fig:repro_timeParallelError_plot_1}--\ref{fig:repro_timeParallelError_plot_2} report parareal convergence for
\reviewerB{these} methods \ourReReading{for $\nptcrs = 10$}. As expected, the proposed \globalForeLabel-\localForeLabel\
method produces a (near) zero error after initialization; on the other hand, the
typical \BELabel-\BELabel\ \reviewerB{and \CNLabel-\CNLabel\ methods exhibit} relatively slow convergence.

\subsection{Predictive case}
\label{ssec:predictive_case}

We now return to the original problem setup with $\ntrain = 4$ training points
and $\nonline = 2$ online points.
\reviewerB{Here, the `ideal case' Assumptions \ref{ass:subspace}--\ref{ass:isomorphic} no
longer hold. To
assess the accuracy of the coarse propagator in this predictive scenario, Figure \ref{fig:proj_error}
reports 
the relative projection error
\begin{gather} 
	\projerrorj(\param)\defeq 
\|
(\identityArg{\nptfn} -
\timebasisj[\timebasisj]^T)\unrollfunc{\redstateEntry{j}(\cdot,\param)}
\| / \|
\unrollfunc{\redstateEntry{j}(\cdot,\param)}\|
\end{gather} 
$j\innat{\nstate}$,
which measures the ability of the
temporal bases $\timebasisj$ to capture the time
evolution of the reduced states. 
Note that this is a global variant of the
quantity that appears in the coarse-propagator error bound in Theorem
\ref{lem:coarseError} and measures the extent to which Assumption \ref{ass:subspace} is
violated. Further, note that $\projerrorj=0$, $j\innat{\nstate}$ for the ideal case.
This figure also reports the relative magnitude of each reduced state
\begin{gather}
\relmagj(\param)\defeq
\|
\unrollfunc{\redstateEntry{j}(\cdot,\param)}\|/
\|
\unrollfunc{\redstate(\cdot,\param)}\|, j\innat{\nstate}.
\end{gather}
\begin{figure}[htbp!] 
	\centering
	\subfigure[Projection error $\projerrorj$,
	$j\innat{\nstate}$]{
	\includegraphics[width=0.48\textwidth]{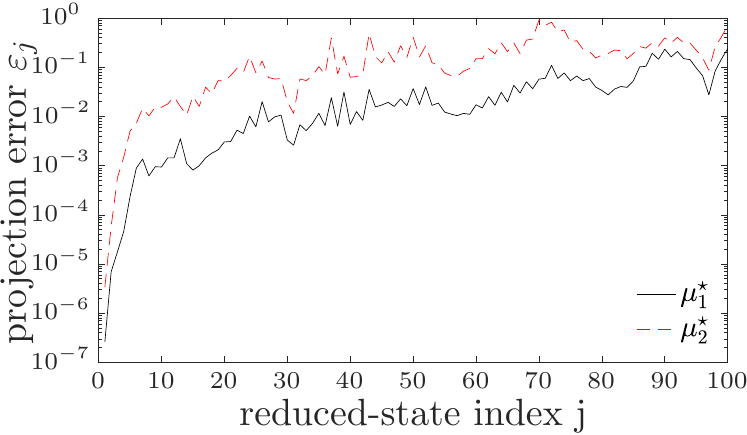}}
	\subfigure[Solution magnitude $\relmagj$,
	$j\innat{\nstate}$]{
	\includegraphics[width=0.48\textwidth]{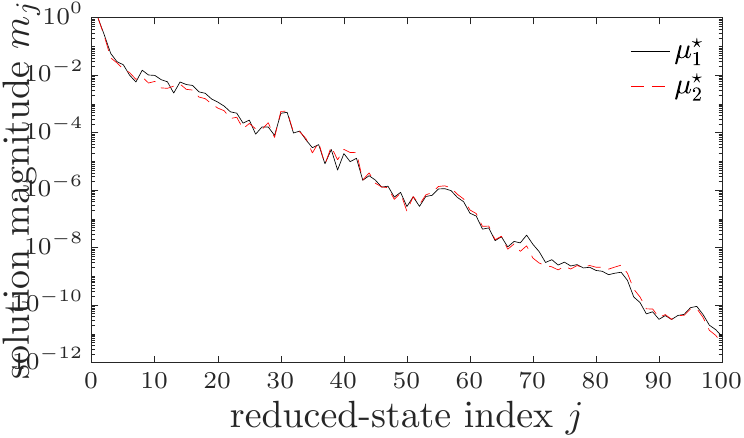}}
 \caption{\textit{Predictive case}. Projection error and solution magnitude for online points $\paramOnlinei{1}$ and $\paramOnlinei{2}$.}
\label{fig:proj_error}
\end{figure} 
Figure \ref{fig:proj_error} shows that the temporal bases are more accurate
(i.e., yield smaller projection errors)
for online point $\paramOnlinei{1}$ than for 
$\paramOnlinei{2}$; this suggests that the method should perform
better (i.e., converge in fewer parareal iterations) for the first online
point. Thus, we can interpret $\paramOnlinei{1}$ and $\paramOnlinei{2}$ as
providing increasingly difficult scenarios for the proposed method in which
the time-evolution bases are increasingly inaccurate. In addition, the figure shows an inverse relationship between the
projection error and the solution magnitude. This is intuitive: the
time-evolution bases are able to accurately capture the time-evolution of the
dominant (low-index) reduced states, while the `noisy' (high-index) reduced
states yield large projection errors.
 Section \ref{sec:param_study} explores this effect further.}

\begin{figure}[htbp!] 
\centering 
\subfigure[$\paramOnlinei{1} = (1.6603, 0.0229)$]{
\includegraphics[width=0.48\textwidth]{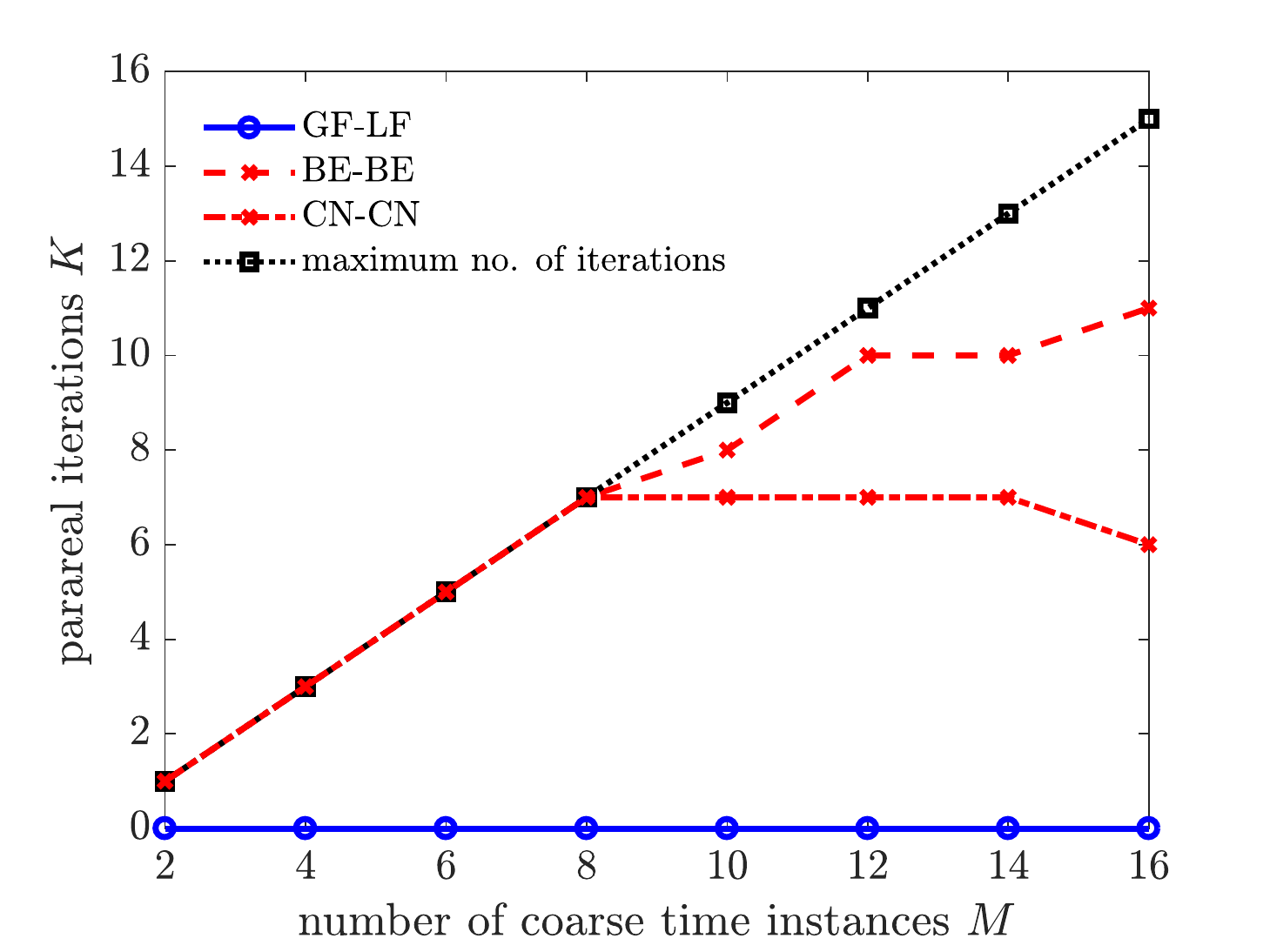} 
\label{fig:iterationVsCPU_plot_1}
}
\subfigure[$\paramOnlinei{2} = (1.5025, 0.0201)$]{
\includegraphics[width=0.48\textwidth]{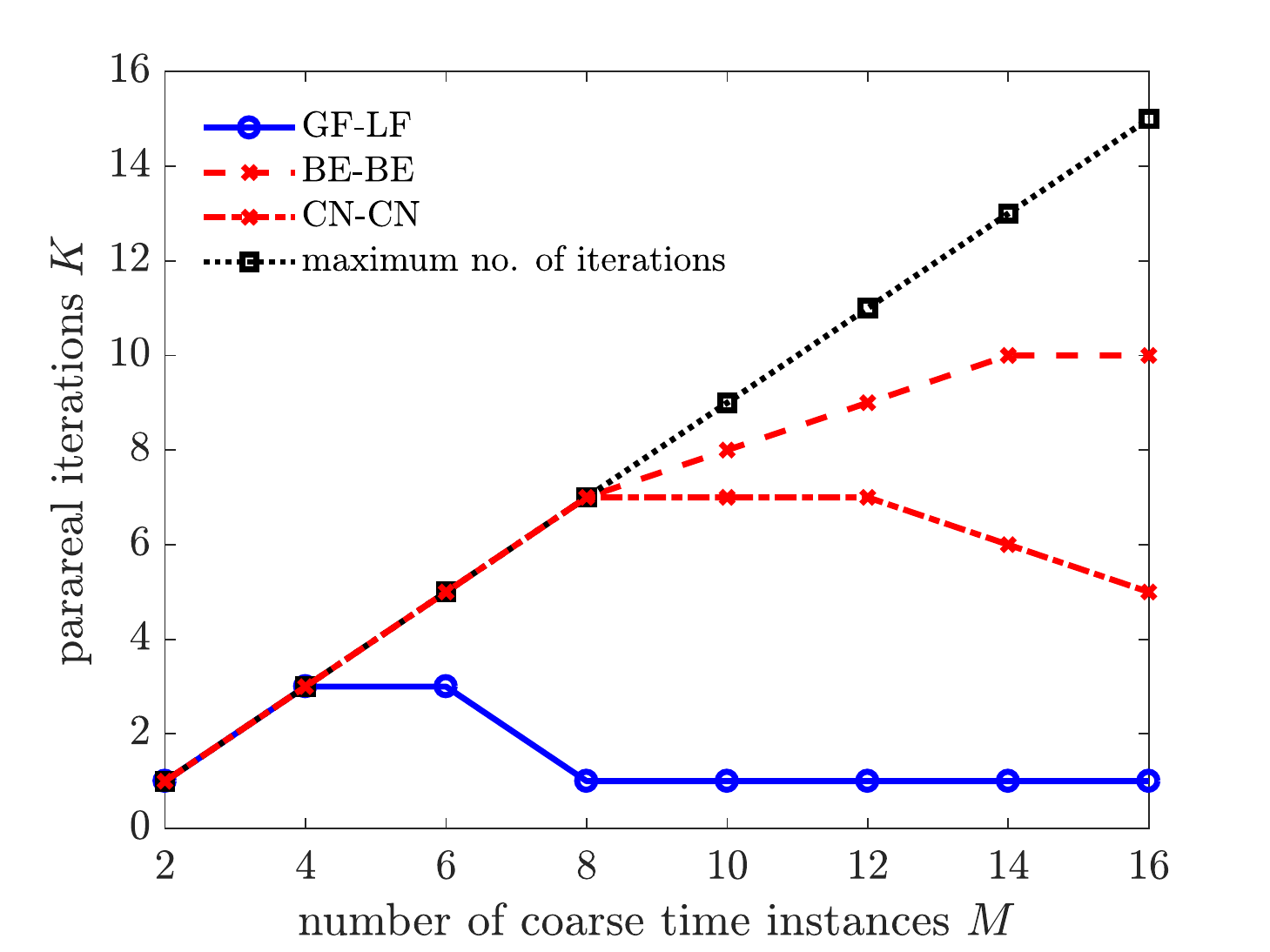} 
\label{fig:iterationVsCPU_plot_2}
}
\centering 
\subfigure[$\paramOnlinei{1} = (1.6603, 0.0229)$]{
\includegraphics[width=0.48\textwidth]{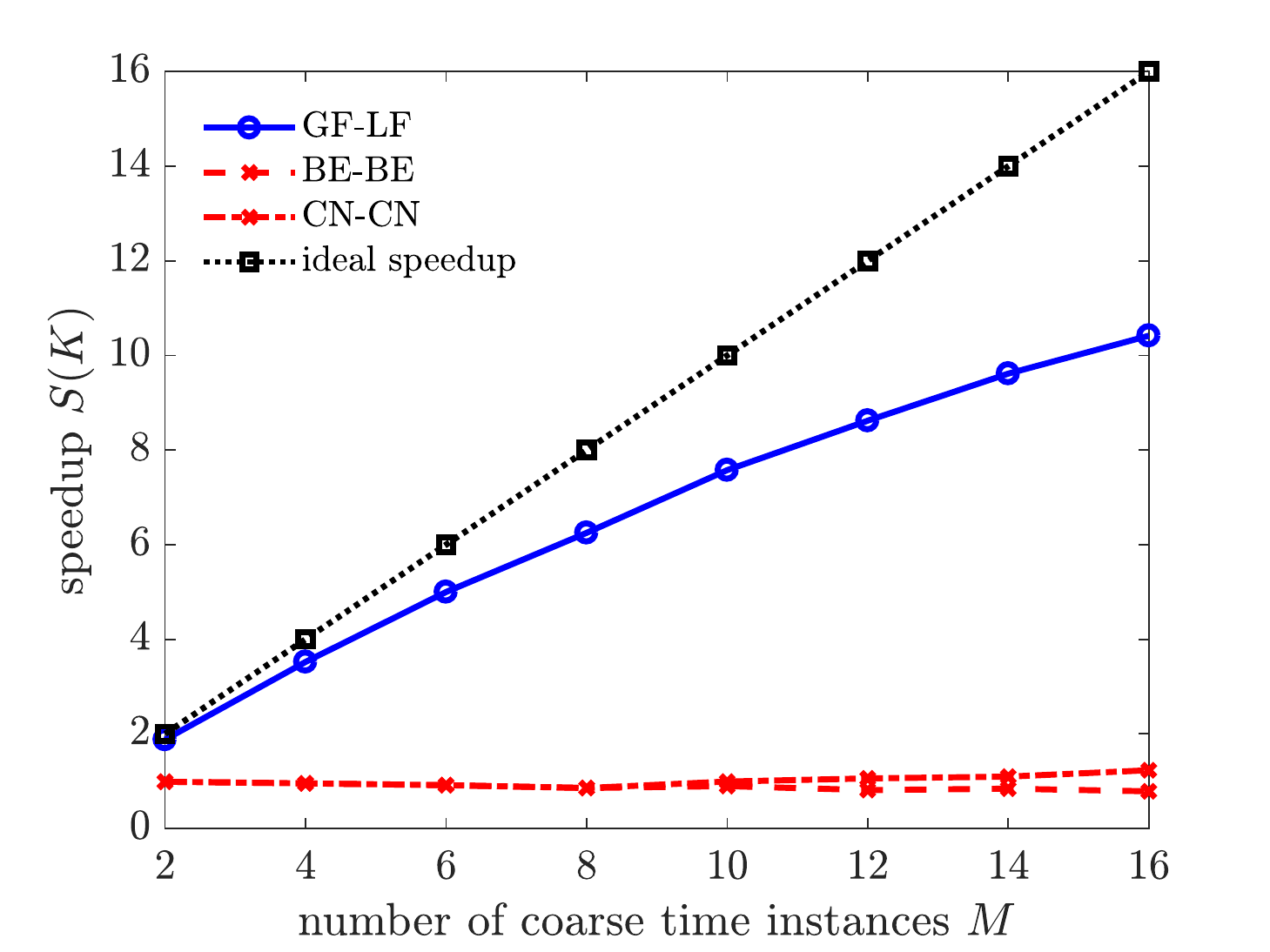} 
\label{fig:speedupVsCPU_ideal_plot_1}
}
\subfigure[$\paramOnlinei{2} = (1.5025, 0.0201)$]{
\includegraphics[width=0.48\textwidth]{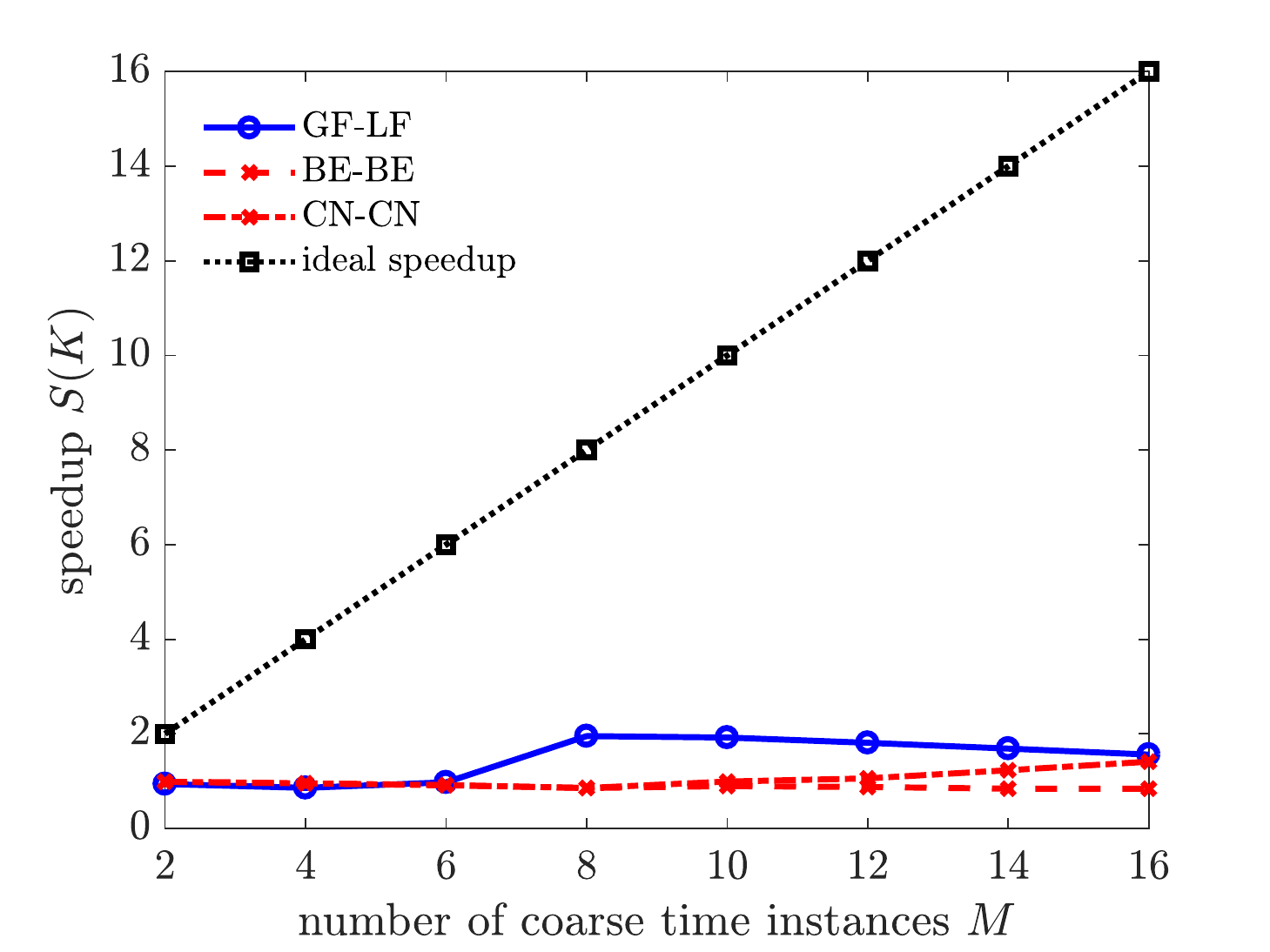} 
\label{fig:speedupVsCPU_ideal_plot_2}
}
\centering 
\subfigure[$\paramOnlinei{1} = (1.6603, 0.0229)$]{
\includegraphics[width=0.48\textwidth]{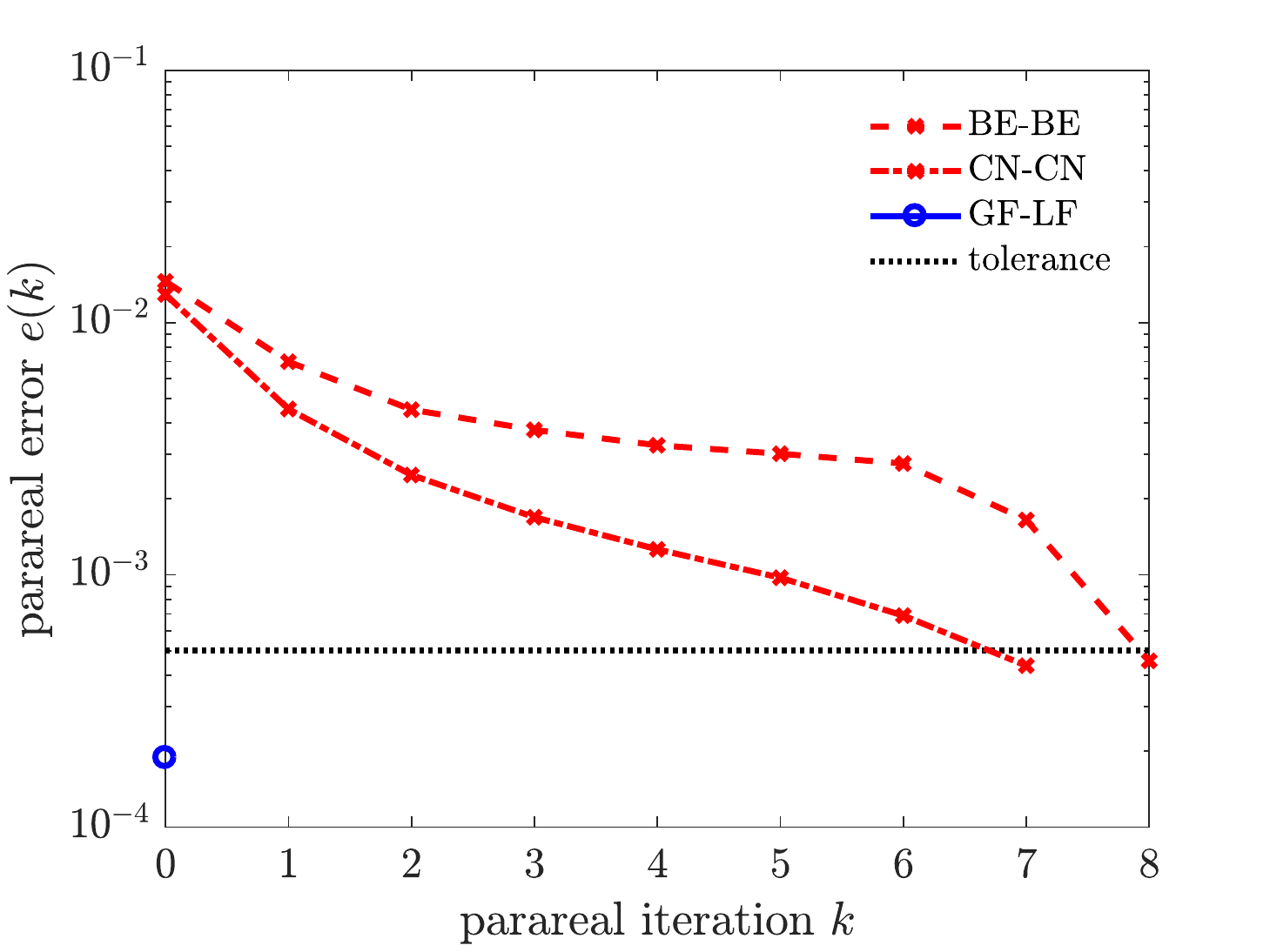} 
\label{fig:timeParallelError_plot_1}
}
\subfigure[$\paramOnlinei{2} = (1.5025, 0.0201)$]{
\includegraphics[width=0.48\textwidth]{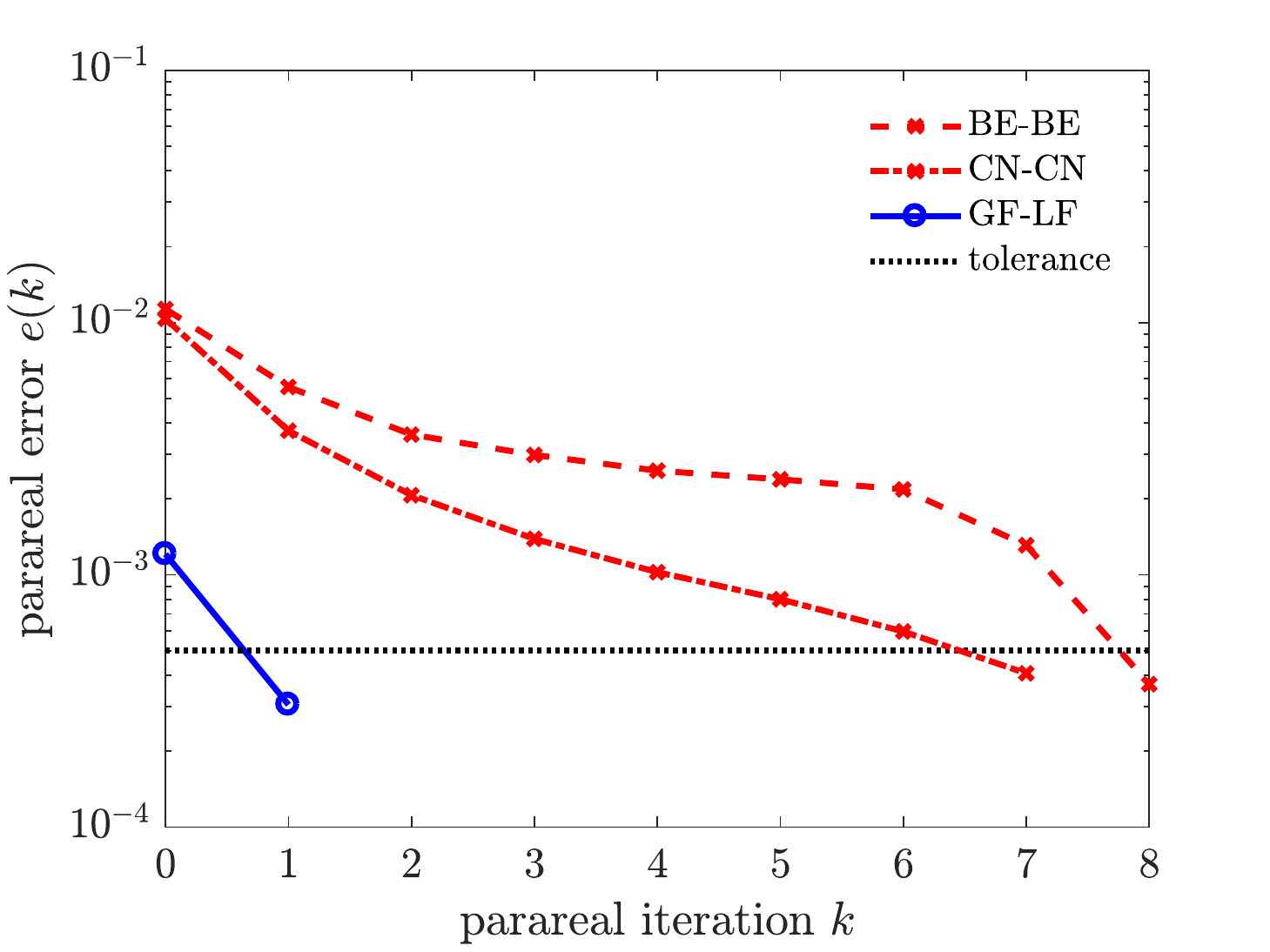} 
\label{fig:timeParallelError_plot_2}
}
\caption{\textit{Predictive case}. Number of parareal iterations $\pararealItConverge$ required for
convergence, theoretical speedups computed via
Eq.~\eqref{eq:finePararealSpeedup} for \BELabel-\BELabel\ and
Eq.~\eqref{eq:speedupGlobal} for \globalForeLabel-\localForeLabel, and
convergence plots. \ourReReading{Figures (e) and (f) associate with $\nptcrs =
10$}. 
}
\label{fig:predictive_case}
\end{figure}

Figures \ref{fig:iterationVsCPU_plot_1}--\ref{fig:iterationVsCPU_plot_2}
report the dependence of the number of parareal iterations on the number of coarse
time instances $\nptcrs$ \ourReReading{for this case}. 
Similar to the ideal case, the proposed
\globalForeLabel-\localForeLabel\ method converges in 
\ourReReading{considerably fewer}
iterations
than the \BELabel-\BELabel\ \reviewerB{and \CNLabel-\CNLabel\ methods}; in fact it converges in the minimum number
of iterations $\pararealItConverge=0$ for $\paramOnlinei{1}$.
\reviewerB{Also, 
	the proposed \globalForeLabel-\localForeLabel\ method exhibits better performance for 
	$\paramOnlinei{1}$ than $\paramOnlinei{2}$ as was suggested by the
projection errors in Figure \ref{fig:proj_error}.}
	As before, the
	\BELabel-\BELabel\ \reviewerB{and \CNLabel-\CNLabel\ methods converge} in the \textit{worst-case} number of
iterations (i.e., $\pararealItConverge=\nptcrs-1$) for $\nptcrs\leq 8$
\reviewerB{for both online points}.  \reviewerB{However, for $\nptcrs\geq 9$,  the \CNLabel-\CNLabel\
	method converges in fewer iterations than the \BELabel-\BELabel\
method, likely due to its higher-order accuracy}.
Figures \ref{fig:speedupVsCPU_ideal_plot_1}--\ref{fig:speedupVsCPU_ideal_plot_2} report the
theoretical speedups of both methods under these ideal conditions. Again, the
proposed technique yields 
better speedups compared with the typical
\reviewerB{methods}, \ourReReading{which is apparent for $\paramOnlinei{1}$ in particular}.

Finally, Figures \ref{fig:timeParallelError_plot_1}--\ref{fig:timeParallelError_plot_2} report parareal convergence for
both methods \ourReReading{for $\nptcrs = 10$}. The proposed \globalForeLabel-\localForeLabel\
method produces a small error after initialization; for $\paramOnlinei{i}$,
the error is smaller than the specified threshold for convergence. In contrast,
the \BELabel-\BELabel\ \reviewerB{and \CNLabel-\CNLabel\ methods exhibit} 
relatively slow convergence \reviewerB{with \CNLabel-\CNLabel\ converging
faster, likely due to its higher-order accuracy}.

These promising results suggest that the proposed
\globalForeLabel-\localForeLabel\ method can deliver significant performance
improvements over standard parareal techniques, even when ideal conditions do
not hold. \reviewerB{We note that numerical results obtained for $\nstate =
	50$ (i.e., a less accurate reduced-order model) reproduce exactly the
	results reported in Figure \ref{fig:predictive_case}, which correspond to
	$\nstate = 100$. This reflects the fact that the proposed
	method's performance is not directly tied to the accuracy of the
reduced-order model; rather, it depends on the ability of the time-evolution
bases to capture the time evolution of the reduced states as discussed above.}

\subsection{Parameter study}\label{sec:param_study}

We now assess the dependence of the proposed \globalForeLabel-\localForeLabel\
method on its parameters, namely the number of restricted variables
$\nrestrict$ and the memory $\memory$.

\begin{figure}[htbp!] 
\subfigure[$\paramOnlinei{1} = (1.6603, 0.0229)$]{
\includegraphics[width=0.45\textwidth]{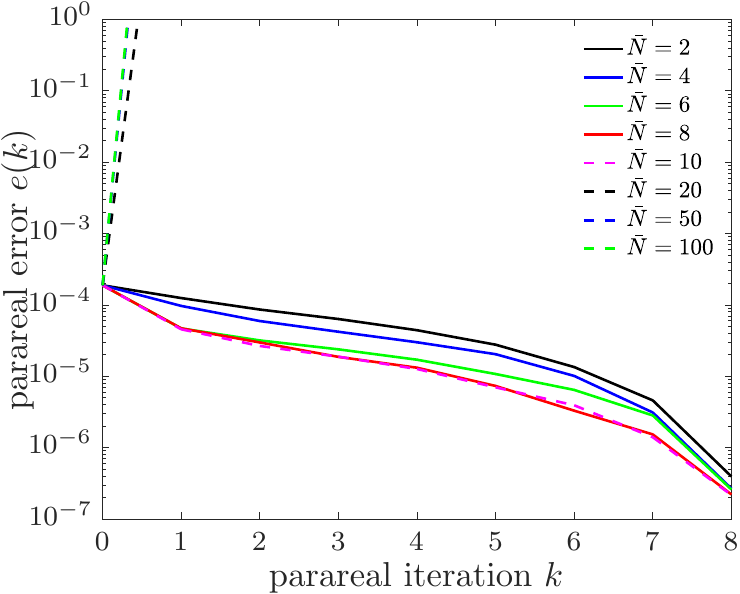} 
\label{fig:nY_Test_plot_1}
}
\subfigure[$\paramOnlinei{2} = (1.5025, 0.0201)$]{
\includegraphics[width=0.45\textwidth]{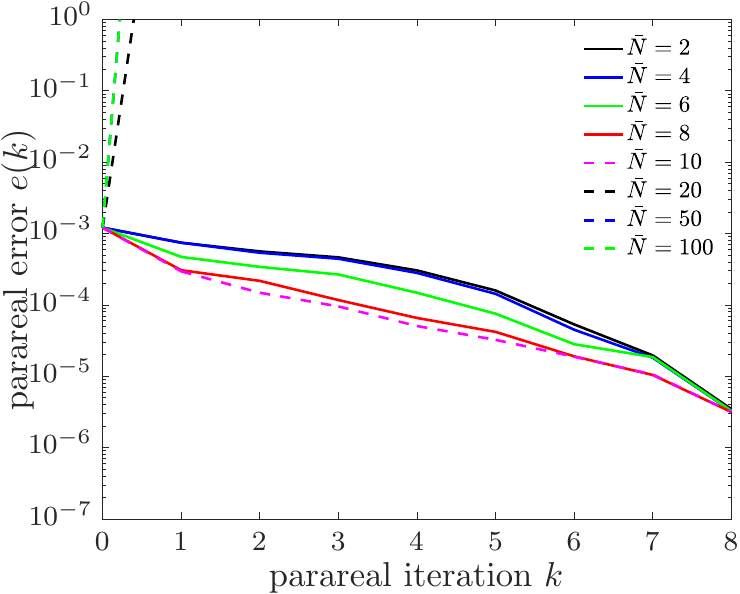} 
\label{fig:nY_Test_plot_2}
}
\caption{\textit{Parameter study}. Convergence plots for the \globalForeLabel-\localForeLabel\
	method for
a range of values for
$\nrestrict$ and a memory of $\memory=8$. Note that a value of $\nrestrict=8$
yields roughly the best overall performance.
}
\label{fig:nY_Test_plot}
\end{figure}

\reviewerB{We first assess the effect of the number of restricted variables
	$\nrestrict$. Recall from Figure \ref{fig:proj_error} that there is an
	inverse relationship between projection error and solution magnitude.
	In particular, low-index reduced states have large solution magnitudes and yield low
	projection errors; high-index reduced states comprise `noise' that cannot be
accurately forecasted due to their high projection errors. 
 To gain additional insight into this, Figure \ref{fig:time_evolution_basis} plots the global
temporal bases $\timebasisj$ associated with different (restricted) solution
components. 
\begin{figure}[htbp!] 
\centering
\subfigure[$\timebasisjArg{1}$]{
\includegraphics[width=0.45\textwidth]{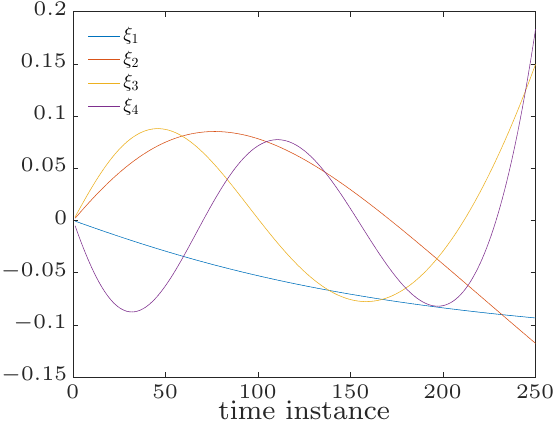} 
\label{fig:time_evolution_basis1}
}
\subfigure[$\timebasisjArg{10}$]{
\includegraphics[width=0.45\textwidth]{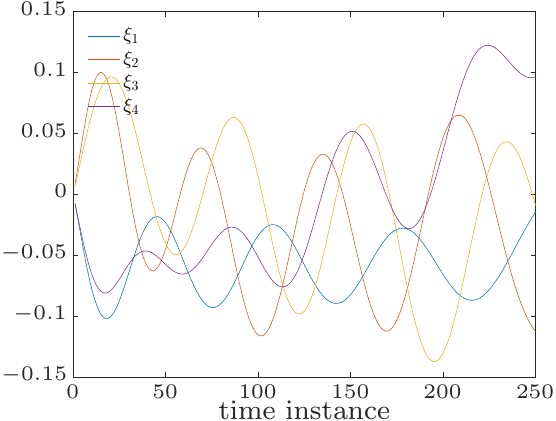} 
\label{fig:time_evolution_basis10}
}
\subfigure[$\timebasisjArg{20}$]{
\includegraphics[width=0.45\textwidth]{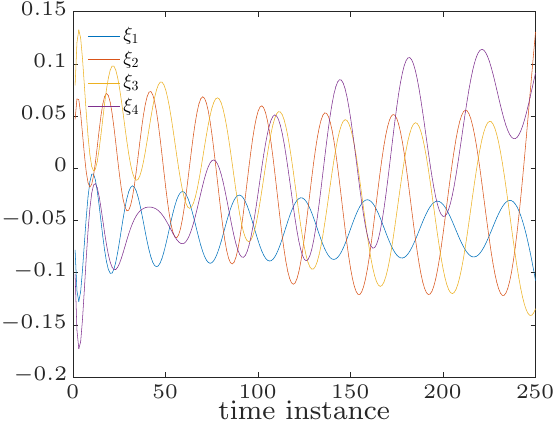} 
\label{fig:time_evolution_basis20}
}
\subfigure[$\timebasisjArg{100}$]{
\includegraphics[width=0.45\textwidth]{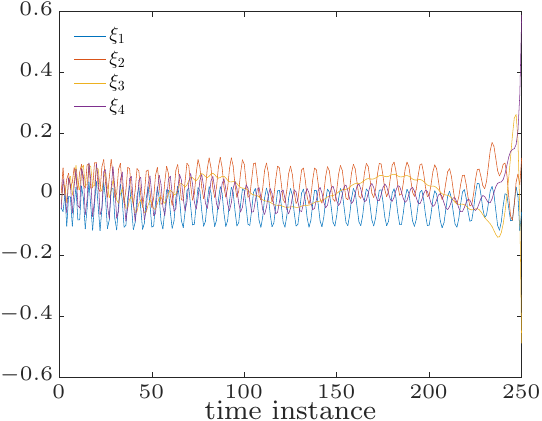} 
\label{fig:time_evolution_basis100}
}
\caption{\textit{Parameter study}. Visualization of global time-evolution bases
$\timebasisj$, $j\in\{1,10,20,100\}$. Time-evolution bases for high-index
POD modes are more highly oscillatory; thus, high-index modes are less
amenable to forecasting.}
\label{fig:time_evolution_basis}
\end{figure} 
Note that the basis vectors are highly
oscillatory for high-index modes, which is consistent with their low relative
magnitudes and interpretation as 
solution
`noise,' as well as their associated large projection errors.
This is consistent
with the discussion in Remark \ref{rem:restrictionTradeoff}: selecting a small value of $\nrestrict$ amounts to forecasting a small
number of solution components, which 
increases the
quantity $\|\restrictOrthArg{\fnprop{\Tmp}{\Tm}{\stateVar}}\|$ appearing
in the bound \eqref{lem:coarseError} for the coarse-propagator error;
alternatively, employing a large value of $\nrestrict$ 
increases the second term in bound \eqref{lem:coarseError} due to the large
projection errors for high-index reduced states. 
Thus, we expect an
intermediate value of $\nrestrict$ to yield the fastest convergence.
}
Figure \ref{fig:nY_Test_plot} reports convergence of the method for $\memory = 8$ for a
range of values for $\nrestrict$. These results show \reviewerB{precisely what
we expect:} the best performance is
obtained (roughly) for an intermediate value of $\nrestrict=8$. 

We next consider the effect of the restricted-state dimension $\nrestrict$
purely when performing global-forecast initialization. 
We find that the
the parareal error after initialization is 
	$\timeParallelError(0) = 1.88\times 10^{-4}$ for $\paramOnlinei{1}$ and
	$\timeParallelError(0) = 1.21\times 10^{-3}$ for
	$\paramOnlinei{\ourReReading{2}}$ and
	for $\nrestrict \in \{6,8,10,15,20,25,30\}$ and a (fixed) memory of $\memory=8$.
		Thus, 
initialization error is insensitive to the parameter $\nrestrict$;
this is an artifact of the intrinsic stability
of the global forecast as discussed in Remark \ref{remark:initialSeedStable}.
Further, it suggests that the first few restricted POD modes
dominate the state information content.

\begin{figure}[htbp!] 
\centering 
\subfigure[$\paramOnlinei{1} = (1.6603, 0.0229)$]{
\includegraphics[width=0.45\textwidth]{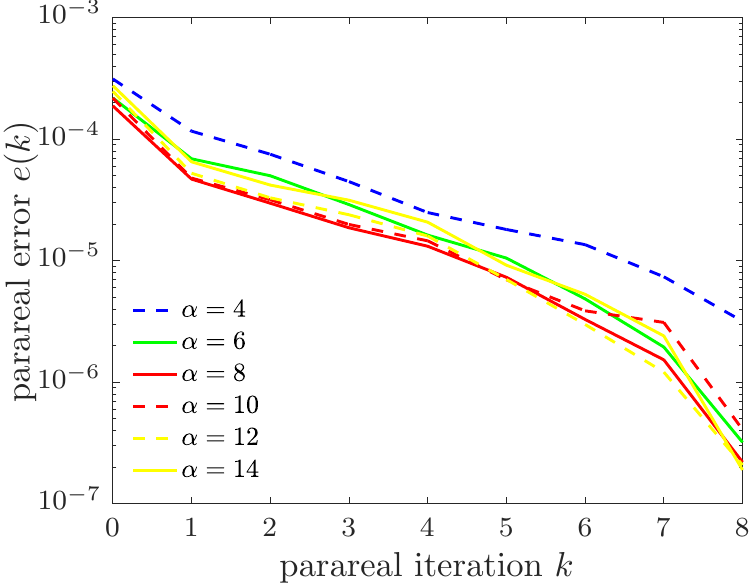} 
\label{fig:memory_Test_plot_1}
}
\subfigure[$\paramOnlinei{2} = (1.5025, 0.0201)$]{
\includegraphics[width=0.45\textwidth]{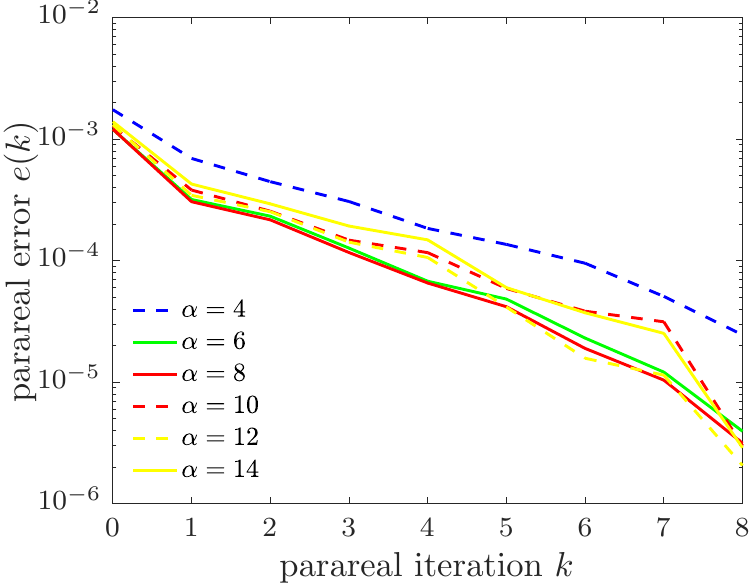} 
\label{fig:memory_Test_plot_2}
}
\subfigure[$\paramOnlinei{1} = (1.6603, 0.0229)$]{
\includegraphics[width=0.45\textwidth]{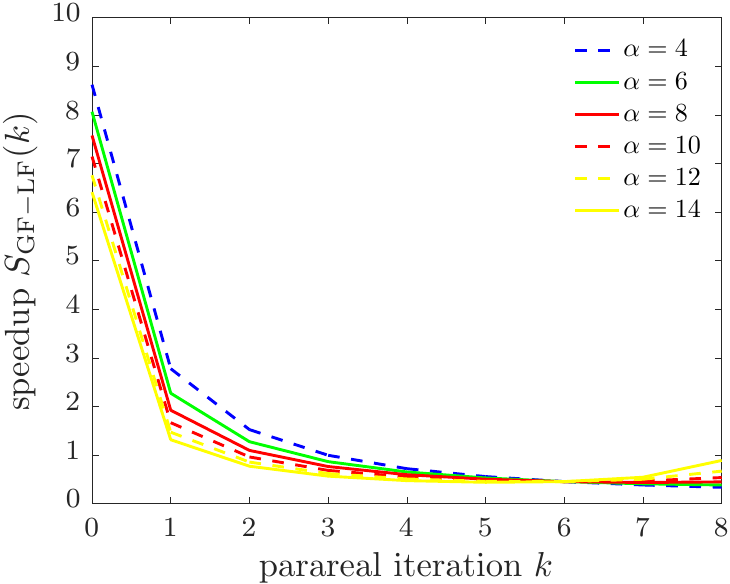} 
\label{fig:memory_Test_theoSpeedup_plot_1}
}
\subfigure[$\paramOnlinei{2} = (1.5025, 0.0201)$]{
\includegraphics[width=0.45\textwidth]{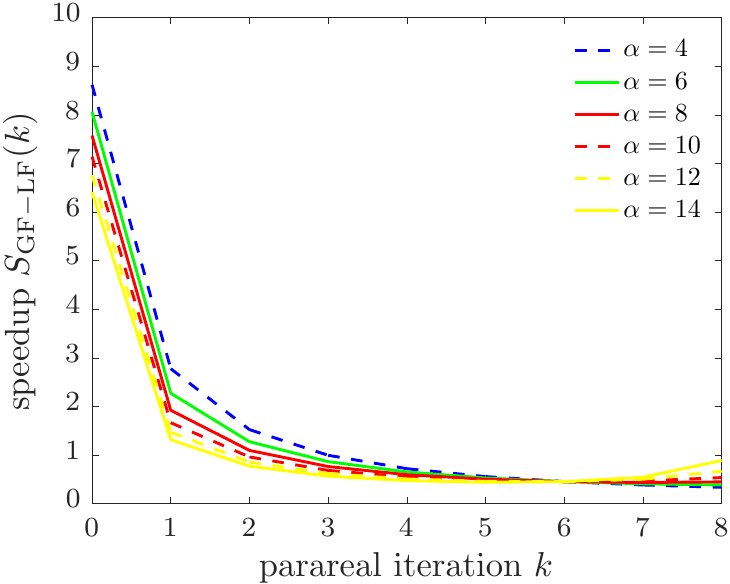} 
\label{fig:memory_Test_theoSpeedup_plot_2}
}
\caption{\textit{Parameter study}. Convergence plots for the \globalForeLabel-\localForeLabel\
method for
a range of values for
$\memory$ and a value of $\nrestrict=8$. Note that a value of $\memory=8$
yields the best overall performance. 
}
\label{fig:memory_Test_plot}
\end{figure}

Next, Figure \ref{fig:memory_Test_plot} reports performance of the method for
a fixed value of $\nrestrict=8$ and a range of values for the memory
$\memory$. First, note that interpolation, which corresponds to $\memory =
\dimBasisjm=4$,
yields the
worst performance in terms of error at a given iteration. This supports the theoretical results discussed in Remark
\ref{rem:interpOversampling}: oversampling (i.e., employing
$\memory>\dimBasisjm$) produces a stabilizing effect. In this case, the value
of the memory leading to best overall performance (in terms of accuracy) is $\memory=8$.
Note that employing the smallest value for the memory yields the best
theoretical speedups if the method were to converge in the same number of
parareal iterations for all values of the memory. This illustrates the tradeoff discussed in Remark
\ref{rem:memoryTradeoff}: increasing the memory $\memory$ reduces the
speedup for a fixed number of iterations needed for convergence; yet, doing so
can also decrease the bound for the error between coarse and fine propagators,
which promotes convergence.


\section{Conclusions}
\label{sec:conclusions} 

This work presented a novel data-driven method for time parallelism. We
applied both local and global forecasting to define initialization
methods, as well as local forecasting to define the coarse propagator. These
methods are data-driven, as they leverage the availability of time-domain data
from which 
low-dimensional time-evolution bases for the state can be extracted; further,
they are well-suited for POD-based reduced-order models, as the required
time-domain data are already available. We performed analysis
demonstrating the method's accuracy, speedup, and stability. Key theoretical
results include:
 \begin{itemize} 
  \item The error between the local-forecast
	coarse propagator and the fine propagator can be bounded by a readily
	interpretable quantity (Theorem \ref{lem:coarseError}),
	\item Ideal conditions exist under which the local-forecast coarse
	propagator is equal to the fine propagator (Theorem \ref{thm:exactCoarse}),
	and
	\item The parareal recurrence is stable with the local-forecast coarse
	propagator (Theorem \ref{thm:finalStability}) with constants that are
	independent of the time discretization (Remark
	\ref{eq:discretizationDependence} and Figure
	\ref{fig:stab_constant_converge_overall}).
\item \reviewerA{Existing convergence results for the parareal recurrence hold
with the proposed coarse propagator, and superlinear convergence can be
obtained under certain conditions (Corollary \ref{cor:convergenceProposed}).}
	 \end{itemize}
	 Key results corroborated by both theoretical analysis and numerical
	 experiments include:
 \begin{itemize} 
	\item Global-forecast initialization is more stable (Remark
	\ref{remark:initialSeedStable}) and produces a more
	accurate solution (Figure \ref{fig:methodComp_plot}) than the local-forecast
	initialization,
	\item Local-forecast coarse propagation is nearly always more accurate than
	backward-Euler coarse propagation, regardless of initialization (Figure
	\ref{fig:methodComp_plot}),
	\item Under ideal conditions, the proposed method converges after
	parareal initialization, i.e., for $\pararealItConverge=0$ in Algorithm
	\ref{alg:parareal}
	(Theorems
	\ref{thm:idealSpeedup}--\ref{thm:idealSpeedupGlobal} and Figures
	\ref{fig:repro_iterationVsCPU_plot_1}--\ref{fig:repro_iterationVsCPU_plot_2}), and
	can realize 
	near-ideal speedups (Figures
	\ref{fig:theoretical_speedup} and
	\ref{fig:repro_speedupVsCPU_ideal_plot_1}--\ref{fig:repro_speedupVsCPU_ideal_plot_2}), 
  \item Increasing the memory $\memory$ can improve coarse-propagation
	accuracy, but incurs a larger cost (Remark \ref{rem:interpOversampling} and
	Figure \ref{fig:memory_Test_plot}), and
	\item Increasing the number of variables included in the forecast
	$\nrestrict$ has two competing effects: it can improve the forecast
	accuracy, but can incur error if the additional variables are difficult to
	forecast, e.g., associate with high-frequency temporal content (Remark
	\ref{rem:restrictionTradeoff} and \reviewerB{Figures \ref{fig:proj_error} and} \ref{fig:nY_Test_plot}).
	 \end{itemize}
Finally, numerical experiments show that in all (predictive) cases where ideal conditions
do not hold, global-forecast
initialization and local-forecast coarse propagation outperforms
backward-Euler initialization and coarse propagation (Figures
\ref{fig:methodComp_plot} and \ref{fig:predictive_case}).

Future work involves applying the proposed methodology in parallel computing
environments with realistic timings, applying the method to parameterized
full-order ODEs (i.e., not reduced-order models), and assessing the viability
of alternative data sources (including physical experiments) to produce
low-dimensional time-evolution bases.

\section*{Acknowledgments}
Kevin Carlberg acknowledges Julien Cortial for insightful and productive
initial conversations on the subject, as well as Yvon Maday for providing
useful feedback.  Sandia National Laboratories is a multi-program laboratory
managed and operated by Sandia Corporation, a wholly owned subsidiary of
Lockheed Martin Corporation, for the U.S. Department of Energy’s National
Nuclear Security Administration under contract DE-AC04-94AL85000.  The
research of \ourReReading{Andrea Barth, Lukas Brencher,
and Bernard Haasdonk} has received funding from the
German Research Foundation (DFG) as part of the Cluster of Excellence in
Simulation Technology (EXC 310/2) at the University of Stuttgart and it is
gratefully acknowledged.
\ourReReading{We also thank the anonymous reviewers for their valuable feedback.}


\bibliography{references}
\bibliographystyle{siam}

\newpage
\appendix
\section{Proofs}\label{sec:proofs}

\begin{proofofthing}{Theorem}{lem:coarseError}
Under Assumptions \ref{ass:prolongOrtho} and \ref{ass:lipschitzProlong}, we have
\begin{align}\label{eq:FGbound}
\begin{split}
\|\fnprop{\Tmp}{\Tm}{\stateVar} -
&\crspropfore{\Tmp}{\Tm}{\stateVar}
\| = 
\|
\prolongateOrthArg{\restrictOrthArg{\fnprop{\Tmp}{\Tm}{\stateVar}}} + 
\prolongateArg{\restrictArg{\fnprop{\Tmp}{\Tm}{\stateVar}}}
-\crspropfore{\Tmp}{\Tm}{\stateVar}
\|\\
&\leq
\lipschitzProlongOrth\|
\restrictOrthArg{\fnprop{\Tmp}{\Tm}{\stateVar}}\| + 
\underbrace{\|
\prolongateArg{\restrictArg{\fnprop{\Tmp}{\Tm}{\stateVar}}}
-\crspropfore{\Tmp}{\Tm}{\stateVar}
\|}_{\text{(I)}}.
\end{split}
\end{align}
We can bound Term (I) using Eqs.~\eqref{eq:coarsePropDefLocalProp} and \eqref{eq:dataDrivenCoarse1} and~\ref{ass:lipschitzProlong} as follows:
\begin{align}\label{eq:stabilityBoundDifferenceOne}
\text{(I)}
&\leq \lipschitzProlong
\|
\restrictArg{\fnprop{\Tmp}{\Tm}{\stateVar}}
-\left[
	\crspropforesym{1}{\timestepdummy }(\stateVar)\
	\cdots\
	\crspropforesym{\nrestrict}{\timestepdummy }(\stateVar)
	\right]^T
\|
=\lipschitzProlong\|\forecastQuantityVecArg{\timestepdummy }{\stateVar}\|,
\end{align}
where we have defined
$\forecastQuantityVecArg{\timestepdummy }{\stateVar}\defeq[\forecastQuantityArg{1}{\timestepdummy }{\stateVar}\
\cdots\ \forecastQuantityArg{\nrestrict}{\timestepdummy }{\stateVar}]^T$ as the vector of errors in the
local forecast, with
\begin{align}
\label{eq:betaDefOne}
\forecastQuantityArg{j}{\timestepdummy }{\stateVar}&\defeq
\restrictjArg{\fnprop{\Tmp}{\Tm}{\stateVar}} -
\crspropforesym{j}{\timestepdummy }(\stateVar)
\\
\label{eq:deltaDefProjErrorSampled}
 &= 
\unitvec{\nTimestepsIntervalArg{\timestepdummy }}^T
\left(
\identityArg{\nTimestepsIntervalArg{\timestepdummy }} - 
\timebasisjm\left[\sampleMat{
0}{\memory}\timebasisjm\right]^+\sampleMat{
0}{\memory}\right)
\left[
\begin{array}{c}
\restrictjArg{\fnprop{\tptfn{\crstofn{\timestepdummy }+1}}{\Tm}{\stateVar}}
-\restrictjArg{\stateVar}\\
\vdots\\
\restrictjArg{\fnprop{\Tmp}{\Tm}{\stateVar}}
-\restrictjArg{\stateVar}
\end{array}
\right].
\end{align}
Using the norm-equivalence relation $\|{\bf x}\|_2\leq \|{\bf x}\|_1$, we have
from 
\eqref{eq:stabilityBoundDifferenceOne} that
$\text{(I)}\leq\lipschitzProlong\sum_{j=1}^\nrestrict
|\forecastQuantityArg{j}{\timestepdummy }{\stateVar}|
$
with 
\begin{align}\label{eq:boundProjErrorGappy}
|\forecastQuantityArg{j}{\timestepdummy }{\stateVar}| &\leq
\Big\|
\left(
\identityArg{\nTimestepsIntervalArg{\timestepdummy }} - 
\timebasisjm\left[\sampleMat{
0}{\memory}\timebasisjm\right]^+\sampleMat{
0}{\memory}\right)\left[
\begin{array}{c}
\restrictjArg{\fnprop{\tptfn{\crstofn{\timestepdummy }+1}}{\Tm}{\stateVar}}
-\restrictjArg{\stateVar}\\
\vdots\\
\restrictjArg{\fnprop{\Tmp}{\Tm}{\stateVar}}
-\restrictjArg{\stateVar}
\end{array}
\right]
\Big\|\\
\label{eq:boundProjError}
&\leq
\normedQuantityGappyLarger{j}{\timestepdummy }\Big\|
(\identityArg{\nptfnArg{\timestepdummy }} - \timebasisjm[\timebasisjm]^T)\left[
\begin{array}{c}
\restrictjArg{\fnprop{\tptfn{\crstofn{\timestepdummy }+1}}{\Tm}{\stateVar}}
-\restrictjArg{\stateVar}\\
\vdots\\
\restrictjArg{\fnprop{\Tmp}{\Tm}{\stateVar}}
-\restrictjArg{\stateVar}
\end{array}
\right]
\Big\|.
\end{align}
Here, we have used a bound for the gappy POD error \cite[Appendix
D]{carlberg2013gnat} and orthogonality of the time-evolution bases $\timebasisjm$.
Note that \eqref{eq:boundProjErrorGappy} expresses the bound in terms of the
gappy POD approximation error of the time evolution (restricted) state, while \eqref{eq:boundProjError} does so in terms of the
orthogonal projection error onto the 
time-evolution bases.
\end{proofofthing}

\begin{proofofthing}{Lemma}{lem:localSubspace}
If Assumption \ref{ass:subspace} holds,
then
$\unrollfunc{\restrictjArg{\stateNP}}=\timebasisj\redCoordRestrictArgs{j}$ for some
$\redCoordRestrictArgs{j}\in\RR{\dimBasisj}$. Applying the $\timestepdummy $th sampling matrix
and utilizing quantities defined in Algorithm \ref{alg:localBasis} with
Assumption \ref{ass:notruncate} (i.e., $\energyCrit = 1.0$)
yields
\begin{align}\label{eq:projectionhofrj}
\sampleMatLocal{\timestepdummy }\unrollfunc{\restrictjArg{\stateNP}}=\sampleMatLocal{\timestepdummy }\timebasisj\redCoordRestrictArgs{j}=
\vectomat{\umvec}{\dimBasism}\entrytodiag{\singularValue}{\dimBasism}\vectomat{\vmvec}{\dimBasism}^T\redCoordRestrictArgs{j}
= \timebasisjm\redCoordRestrictTwoArgs{j}{\timestepdummy },
\end{align}
where we have defined
$\redCoordRestrictTwoArgs{j}{\timestepdummy }\defeq\entrytodiag{\singularValue}{\dimBasism}\vectomat{\vmvec}{\dimBasism}^T\redCoordRestrictArgs{j}$.
Noting that $\restrictjArg{\stateNP}\in\setOfTimeDepFunctions$ implies $\sampleMatLocal{\timestepdummy }\unrollfunc{\restrictjArg{\stateNP}} =
\unrollmfunc{\restrictjArg{\stateNP}}$, we have the desired
result: $\unrollmfunc{\restrictjArg{\stateNP}}\in\range{\timebasisjm}$,
$j\innat{\ndof}$, $\timestepdummy \innatZero{\nTimeIntervals-1}$.
\end{proofofthing}

\begin{proofofthing}{Theorem}{thm:exactCoarse}
Under 
Assumptions \ref{ass:subspace} and \ref{ass:notruncate}, we have from Lemma \ref{lem:localSubspace} that
$\unrollmfunc{\restrictjArg{\stateNP}}\in\range{\timebasisjm}$,
$j\innat{\ndof}$, $\timestepdummy \innatZero{\nTimeIntervals-1}$.
Then, from
	Eqs.~\eqref{eq:defineGenCoords}, \eqref{eq:localForecast}, and \eqref{eq:coarsePropDefLocalProp} it follows
\begin{align*}
\forecastFunctionjm{\Tmp}{\Tm}{\restrictjArg{\fineFillInLocal{\timestepdummy }{\stateNP(\tptcrs{\timestepdummy })}}}
&=
\forecastFunctionjm{\Tmp}{\Tm}{\restrictjArg{\stateNP}} 
=
\restrictjArg{\stateNP(\Tm)} + 
\unitvec{\crstofn{\timestepdummy +1}-\crstofn{\timestepdummy }}^T\timebasisjm
[\sampleMat{0}{\prevtpts}
\timebasisjm]^+\sampleMat{0}{\prevtpts}
\unrollmfunc{\restrictjArg{\stateNP}}\\
&=
\restrictjArg{\stateNP(\Tm)} + 
\unitvec{\crstofn{\timestepdummy +1}-\crstofn{\timestepdummy }}^T\timebasisjm
\redCoordRestrictTwoArgs{j}{\timestepdummy }\\
&=
\restrictjArg{\stateNP(\Tm)} + 
\unitvec{\crstofn{\timestepdummy +1}-\crstofn{\timestepdummy }}^T\unrollmfunc{\restrictjArg{\stateNP}}=
 \restrictjArg{\stateNP(\Tmp)}
 =
 \restrictjArg{\fnprop{\Tmp}{\Tm}{\stateNP(\Tm)}},
\end{align*}
where we have used 
$\unrollmfunc{\restrictjArg{\stateNP}} = \timebasisjm\redCoordRestrictTwoArgs{j}{\timestepdummy } 
$ from Eq.~\eqref{eq:projectionhofrj} and 
$[\sampleMat{0}{\prevtpts}
\timebasisjm]^+\sampleMat{0}{\prevtpts}
\timebasisjm = \identity$. 
Leveraging Assumption \ref{ass:isomorphic} (i.e.,
$\prolongate\restrict=\identityArg{\ndof}$ with $\ndof=\nrestrict$)
along with the definition of the coarse propagator \eqref{eq:dataDrivenCoarse1} yields
the desired result.
\end{proofofthing}

\begin{proofofthing}{Theorem}{thm:localForecastGen}
Under Assumption \ref{ass:integrateCostDominant}, the wall time incurred by a
serial solution is $\nptfn\costsolve$, where $\costsolve\in\RRplus$ is the
wall time required to compute
$\fnprop{\tptfn{\timestepdummy +1}}{\tptfn{\timestepdummy }}{\stateNP(\tptfn{\timestepdummy })}$ for a given
$\timestepdummy \innatZero{\nptfn-1}$. 
	
	Under Assumptions \ref{ass:localinitialize} and \ref{ass:integrateCostDominant}, the wall time incurred by
	initializing the proposed method in Steps
	\ref{step:initializeOne}--\ref{step:initializeLast} of Algorithm
	\ref{alg:parareal} is composed of (1) the local-forecast initialization in Step
	\ref{step:initialize} of Algorithm \ref{alg:parareal}, which incurs
	performing (in serial) fine propagation $\memory$ times in each \ourReReading{coarse} time interval (wall time of
	$\nptcrs\memory\costsolve$) and (2) the (worst-case) parallel fine
	propagation in Steps
	\ref{step:beginFinePropInit}--\ref{step:initializeLast}
	(wall time of $(\nptfnArg{\timestepdummy }-\memory)\costsolve$); here, we have exploited
	the fact that we can reuse the first $\memory$ fine propagations on each
	time interval, as these were computed during local-forecast initialization.
	
	Because the local forecast is also employed as a coarse propagator,
	Step \ref{step:initial_seed_prop} of Algorithm \ref{alg:parareal} can
	be replaced by simply setting $\coarsesolFOM{\timestepdummy +1}{0} = \apxsolFOM{\timestepdummy +1}{0}$,
	which incurs no cost under Assumption \ref{ass:integrateCostDominant}. Then,
	each subsequent iteration requires (1) serial coarse propagation in Step
	\ref{step:serialcoarse}
	(wall time of $(\nptcrs-\pararealit )\memory\costsolve$), and (2) parallel fine propagation 
        in Step \ref{step:fineProp1} (wall time
	of $(\nptfnArg{\timestepdummy }-\memory)\costsolve$).
	The ratio of these costs yields the theoretical speedup. 
Finally, we note
	that additional speedups may be realizable by pipelining operations, i.e.,
	initiating the fine propagation on a given coarse time interval as soon as
	its initial value is available.
\end{proofofthing}

\begin{proofofthing}{Theorem}{thm:globalForecastGen}
Under Assumption \ref{ass:integrateCostDominant}, the wall time incurred by a
serial solution is (again) $\nptfn\costsolve$. 
	Under Assumptions \ref{ass:globalinitialize} and \ref{ass:integrateCostDominant}, the wall time incurred by
	initializing the proposed method in Steps
	\ref{step:initializeOne}--\ref{step:initializeLast} of Algorithm
	\ref{alg:parareal} is composed of (1) global-forecast initialization in Step
	\ref{step:initialize} of Algorithm \ref{alg:parareal}, which incurs
	performing fine propagation $\memory$ times in only the \textit{first} time interval (wall time of
	$\memory\costsolve$) and (2) the (worst-case) parallel fine
	propagation in Steps
	\ref{step:beginFinePropInit}--\ref{step:initializeLast}
	(wall time of $\nptfnArg{\timestepdummy }\costsolve$); note
	that we can no longer reuse fine propagation from initialization beyond the
	first time interval.
	
	Because the local forecast is employed as a coarse propagator, 
	Step \ref{step:initial_seed_prop} of Algorithm \ref{alg:parareal} incurs
	parallel coarse propagation, which requires performing (in parallel) fine propagation
	$\memory$ times in time intervals $1$ to $\nptcrs-1$ (wall time of
	$\memory\costsolve$). Then,
	each subsequent iteration requires (1) serial coarse propagation in Step
	\ref{step:serialcoarse}
	(wall time of $(\nptcrs-\pararealit )\memory\costsolve$), and (2) parallel fine 
       propagation in Step \ref{step:fineProp1} (wall time
	of $(\nptfnArg{\timestepdummy }-\memory)\costsolve$).
	The ratio of these costs yields the theoretical speedup. 
\end{proofofthing}

\begin{proofofthing}{Theorem}{thm:idealSpeedup}
We proceed by induction. Assume that $\apxsolFOM{\timestepdummy }{0} = \stateNP(\tptcrs{\timestepdummy })$, which holds
for $\timestepdummy =0$ by construction. Then, we have from 
Theorem \ref{thm:exactCoarse} 
under Assumptions \ref{ass:subspace}--\ref{ass:isomorphic}
that
$
\crspropfore{\Tmp}{\Tm}{\apxsolFOM{\timestepdummy }{0}} = \stateNP(\Tmp)
$.
Under Assumption \ref{ass:localinitialize} (i.e., initialization is
	performed via local forecasting), we have
	$
\apxsolFOM{\timestepdummy +1}{0} =
\crspropforesymAll({\apxsolFOM{\timestepdummy }{0}};{\tptcrs{\timestepdummy }},{\tptcrs{\timestepdummy +1}}) =
\stateNP(\Tmp).
	$
	By induction, this yields $\apxsolFOM{\timestepdummy }{0} = \stateNP(\Tm)$,
	$\timestepdummy \innatZero{\nptcrs}$.
	This means that the initialized values computed in Step
	\ref{step:initialize} of Algorithm \ref{alg:parareal} are correct under the
	stated assumptions; as a result, the fine propagation performed in Steps
	\ref{step:beginFinePropInit}--\ref{step:initializeLast} will complete
	computation of the correct solution, the error measure in Step
	\ref{step:termination} will evaluate to zero, and the algorithm will
	terminate with $\pararealItConverge=0$.
Finally, Theorem \ref{thm:localForecastGen} is valid under Assumptions
\ref{ass:localinitialize} and \ref{ass:integrateCostDominant}; thus,
$\speedupLocal{0}$ provides the theoretical speedup in this case.
\end{proofofthing}

\begin{proofofthing}{Theorem}{thm:idealSpeedupGlobal}
As in Lemma \ref{lem:localSubspace},
if Assumption \ref{ass:subspace} holds, then 
$\unrollfunc{\restrictjArg{\stateNP}}=\timebasisj\redCoordRestrictArgs{j}$ for some
$\redCoordRestrictArgs{j}\in\RR{\dimBasisj}$. Then, we have from
Eqs.~\eqref{eq:genSolutionGlobalForecast} and \eqref{eq:globalForecast}
 \begin{align}
 \begin{split}
\forecastFunctionjNo\forecastFunctionjArgs{\tptcrs{\timestepdummy }}{0}{\restrictEntryArg{j}{\fineFillInOne{\initstateNP}}}
=\forecastFunctionjNo\forecastFunctionjArgs{\tptcrs{\timestepdummy }}{0}{\restrictjArg{\stateNP}}
&=\restrictjArg{\stateNP(0)} +
\unitvec{\crstofn{\timestepdummy }}^T\timebasisj
[\sampleMat{0}{\prevtpts}
\timebasisj]^+\sampleMat{0}{\prevtpts}
\timebasisj\redCoordRestrictArgs{j}
=\restrictjArg{\stateNP(0)} +
\unitvec{\crstofn{\timestepdummy }}^T\timebasisj
\redCoordRestrictArgs{j}\\
&=\restrictjArg{\stateNP(0)} +
\unitvec{\crstofn{\timestepdummy }}^T\unrollfunc{\restrictjArg{\stateNP}} = 
\restrictjArg{\stateNP(\tptcrs{\timestepdummy })},
	\end{split}
  \end{align} 
	where we have used $\unrollfunc{\restrictjArg{\stateNP}} = \timebasisj\redCoordRestrictArgs{j} 
$ and 
$[\sampleMat{0}{\prevtpts}
\timebasisj]^+\sampleMat{0}{\prevtpts}
\timebasisj = \identity$. Under Assumption \ref{ass:globalinitialize} with the definition of the proposed global-forecast
initialization \eqref{eq:globalForecastInitialize}, Assumption
\ref{ass:isomorphic} (i.e., 
$\prolongate\restrict=\identityArg{\ndof}$), and the construction
$\apxsolFOM{0}{0}=\initstateNP=\stateNP(\tptcrs{0})$, we then have
$
\apxsolFOM{\timestepdummy }{0} =\stateNP(\tptcrs{\timestepdummy }),\ \timestepdummy \innatZero{\nptcrs}.
$
As before, this proves the desired result under the stated assumptions:
the
initialized values computed in Step \ref{step:initialize} of Algorithm
\ref{alg:parareal} are correct. Thus, the fine propagation performed in
Steps \ref{step:beginFinePropInit}--\ref{step:initializeLast} will complete
computation of the correct solution, the
error measure in Step \ref{step:termination} will evaluate to zero, and the algorithm will
terminate with $\pararealItConverge=0$.
Finally, Theorem \ref{thm:globalForecastGen} is valid under
Assumptions \ref{ass:globalinitialize} and \ref{ass:integrateCostDominant};
thus, $\speedupGlobal{0}$ provides the theoretical speedup in this case.
\end{proofofthing}

\begin{proofofthing}{Lemma}{lem:pararealStability}
Defining $\overallCoarse\defeq \constCoarse(1 +
\stabilityCoarse\tsscrs)$ and $\overallFineCoarse\defeq\constFineCoarse(1 +
\stabilityFineCoarse\tsscrs)$, we have from parareal recurrence
\eqref{eq:parareal} and bounds \eqref{eq:coarseBoundGen} and
\eqref{eq:fineCoarseBoundGen} that
\begin{align}\label{eq:pararealBound}
\begin{split}
\|\apxsolFOM{\tixcrs +1}{\itvar +1}\|\leq&
\overallCoarse\|\apxsolFOM{\tixcrs}{\itvar +1}\|
+ \overallFineCoarse\|\apxsolFOM{\tixcrs}{\itvar}\|,\quad
\timestepdummy \innatZero{\nptcrs-1},\ \pararealit \innatZero{\timestepdummy },
\end{split}
\end{align}
which can be written equivalently as 
\begin{align}\label{eq:pararealBoundOffset}
\begin{split}
\|\apxsolFOM{\tixcrs +1}{\itvar}\|\leq&
\overallCoarse\|\apxsolFOM{\tixcrs}{\itvar }\|
+ \overallFineCoarse\|\apxsolFOM{\tixcrs}{\itvar-1}\|,\quad
\timestepdummy \innatZero{\nptcrs-1},\ \pararealit \innat{\timestepdummy +1}.
\end{split}
\end{align}
\noindent We prove by induction over $\timestepdummy $ that
\begin{align}
\label{eq:ineqone}\|\apxsolFOM{\timestepdummy }{\itvar}\|\leq&
\overallCoarse^\timestepdummy \|\initstateNP\|+\sum_{j=1}^\timestepdummy \overallCoarse^{\timestepdummy -j}\overallFineCoarse\|\apxsolFOM{j-1}{\itvar-1}\|,\quad
\timestepdummy \innat{\nptcrs},\ \pararealit \innat{\timestepdummy }
\end{align}
Applying inequality \eqref{eq:pararealBound} with $\timestepdummy =\pararealit =0$ yields
$
 \|\apxsolFOM{1}{1}\|\leq \overallCoarse \|\apxsolFOM{0}{1}\| + \overallFineCoarse \|\apxsolFOM{0}{0}\| = (A + B) \|\initstateNP\|,
$
where we have used $\|\apxsolFOM{0}{\itvar}\| = \|\initstateNP\|$, 
$\itvar\innatZero{\pararealItConverge}$ by construction. Thus, inequality \eqref{eq:ineqone}
holds for $\timestepdummy =1$.

Now assume that inequality~\eqref{eq:ineqone}
holds for some $\timestepdummy \innat{\nptcrs-1}$ and all $\pararealit \innat{\timestepdummy }$; we will show that
inequality \eqref{eq:ineqone} is then satisfied for
$\timestepdummy +1\in\{2,\ldots,\nptcrs\}$. Applying inequality
\eqref{eq:pararealBoundOffset} with $\timestepdummy \innat{\nptcrs-1}$, $\pararealit \innat{\timestepdummy +1}$
yields
\begin{align}\label{eq:ineqthree}
 \|\apxsolFOM{\timestepdummy +1}{\itvar}\|&\leq \overallCoarse \|\apxsolFOM{\timestepdummy }{\itvar}\| + \overallFineCoarse \|\apxsolFOM{\timestepdummy }{\itvar-1}\|
 \leq \overallCoarse \Big(\overallCoarse^\timestepdummy \|\initstateNP\|+\sum_{j=1}^\timestepdummy \overallCoarse^{\timestepdummy -j}\overallFineCoarse\|\apxsolFOM{j-1}{\itvar-1}\|\Big) + \overallFineCoarse \|\apxsolFOM{\timestepdummy }{\itvar-1}\|\\
 &\leq
 \overallCoarse^{\timestepdummy +1}\|\initstateNP\|+\sum_{j=1}^{\timestepdummy +1}\overallCoarse^{\timestepdummy +1-j}\overallFineCoarse\|\apxsolFOM{j-1}{\itvar-1}\|,\quad
 \timestepdummy \innat{\nptcrs-1},\ \pararealit \innat{\timestepdummy +1},
\end{align}
which is Eq.~\eqref{eq:ineqone} for $\timestepdummy +1\in\{2,\ldots,\nptcrs\}$.

We now prove that
\begin{align}
\label{eq:ineqtwo}
\|\apxsolFOM{\timestepdummy }{\itvar}\|\leq
(\overallCoarse + \overallFineCoarse)^\timestepdummy  \|\initstateNP\|,\quad
\timestepdummy \innat{\nptcrs},\ \pararealit  = \timestepdummy .
\end{align}
Applying inequality \eqref{eq:pararealBound} with $\pararealit =\timestepdummy $ yields
$
 \|\apxsolFOM{\timestepdummy }{\timestepdummy }\|\leq  \overallCoarse \|\apxsolFOM{\timestepdummy -1}{\timestepdummy }\| +
 \overallFineCoarse \|\apxsolFOM{\timestepdummy -1}{\timestepdummy -1}\|,\quad \timestepdummy \innat{\nptcrs}.
 $
From the finite-termination property in Eq.~\eqref{eq:finiteTermination}, we
have that $\apxsolFOM{\timestepdummy -1}{\timestepdummy } = \apxsolFOM{\timestepdummy -1}{\timestepdummy -1}=\stateNP(\tptcrs{\timestepdummy -1})$, 
thus this inequality becomes
$ \|\apxsolFOM{\timestepdummy }{\timestepdummy }\|\leq  \overallCoarse \|\apxsolFOM{\timestepdummy -1}{\timestepdummy -1}\| +
 \overallFineCoarse \|\apxsolFOM{\timestepdummy -1}{\timestepdummy -1}\|
 \leq (A+B) \|\apxsolFOM{\timestepdummy -1}{\timestepdummy -1}\|
$
from which 
the desired result directly follows.

Substituting the definitions of $\overallCoarse$ and $\overallFineCoarse$ in
inequalities \eqref{eq:ineqone}--\eqref{eq:ineqtwo}, 
employing the generalized binomial formula, and applying the
inequality $(1+x)^n\leq \exp(nx)$ yields the desired result.
\end{proofofthing}

\begin{proofofthing}{Lemma}{lem:coarseStable}
From Eqs.~\eqref{eq:localForecast} and \eqref{eq:coarsePropDefLocalProp}, Assumptions \ref{ass:stable}
and \ref{ass:lipschitzRestrict},
and the norm-equivalence relation $\|{\bf
x}\|_2\leq\sqrt{n}\|{\bf x}\|_\infty$ for ${\bf x}\in\RR{n}$, 
\ourReReading{Eq.~\eqref{eq:coarsePropDefLocalPropAlg} yields}
\begin{align}\label{eq:Gcrsbound}
|\crspropforesym{j}{\timestepdummy }(\stateVar)|
&\leq
|(1-\sum_{i=1}^\memory\forecastScalar{i}{j})\restrictjArg{\stateVar}|
 + |\sum_{i=1}^\memory\forecastScalar{i}{j}\restrictjArg{\fnprop{\tptfn{\crstofn{\timestepdummy }+i}}{\Tm}{\stateVar}}|\\
&=|(1-\sum_{i=1}^\memory\forecastScalar{i}{j})\restrictjArg{\stateVar}|
 + |
\unitvec{\nTimestepsIntervalArg{\timestepdummy }}^T\timebasisjm\left[\sampleMat{
0}{\memory}\timebasisjm\right]^+[\restrictjArg{\fnprop{\tptfn{\crstofn{\timestepdummy }+1}}{\Tm}{\stateVar}} \ \cdots\ 
\restrictjArg{\fnprop{\tptfn{\crstofn{\timestepdummy }+\memory}}{\Tm}{\stateVar}}
]^T|\\
&\leq\lipschitzRestrictEntry{j}\left[|1-\sum_{i=1}^\memory\forecastScalar{i}{j}|
\|\stateVar\|+
\|\unitvec{\nTimestepsIntervalArg{\timestepdummy }}^T\timebasisjm\left[\sampleMat{
0}{\memory}\timebasisjm\right]^+
\|\left(\sum_{\pararealit =1}^{\memory}
\|\fnprop{\tptfn{\crstofn{\timestepdummy }+\pararealit }}{\Tm}{\stateVar}\|
^2\right)^{1/2}\right]\\
\label{eq:keystabilityieq}
&\leq\lipschitzRestrictEntry{j}\left[\normedQuantityGappySmall{j}{\timestepdummy }
+
\normedQuantityGappy{j}{\timestepdummy }\sqrt{\memory}(1 +
\stabilityConstant\memory\tssfn)
 \right]\|\stateVar\| = 
\lipschitzRestrictEntry{j}\left(\normedQuantityGappySmall{j}{\timestepdummy }
+
\normedQuantityGappy{j}{\timestepdummy }\sqrt{\memory}\right)\left[1 +
\frac{\normedQuantityGappy{j}{\timestepdummy }\memory^{3/2}\stabilityConstant\tssfn}{\normedQuantityGappySmall{j}{\timestepdummy }
+
\normedQuantityGappy{j}{\timestepdummy }\sqrt{\memory}}
 \right]\|\stateVar\|
\\
& = 
\lipschitzRestrictEntry{j}\left(\normedQuantityGappySmall{j}{\timestepdummy }
+
\normedQuantityGappy{j}{\timestepdummy }\sqrt{\memory}\right)\left[1 +
\frac{(\memory/\nfinepercoarse)\normedQuantityGappy{j}{\timestepdummy }\sqrt{\memory}}{\normedQuantityGappySmall{j}{\timestepdummy }
+
\normedQuantityGappy{j}{\timestepdummy }\sqrt{\memory}}\stabilityConstant\tsscrs
 \right]\|\stateVar\|.
  \end{align} 
	From Eq.~\eqref{eq:dataDrivenCoarse1} and Assumption
\ref{ass:lipschitzProlong} follows then
$\|\crspropforesymAll(\stateVar;\tptcrs{\tixcrs},\tptcrs{\tixcrs +
	1})\|\leq\sqrt{\nrestrict}\lipschitzProlong\max_{j\in\nat{\nrestrict}}
|\crspropforesym{j}{\timestepdummy }(\stateVar)|,
$
which produces the desired result.
\end{proofofthing}

\begin{proofofthing}{Lemma}{lem:fineCoarseStable}
Under the stated assumptions, we have
\begin{align}\label{eq:FGboundAgain}
\begin{split}
\|\fnprop{\Tmp}{\Tm}{\stateVar} -
\crspropfore{\Tmp}{\Tm}{\stateVar}
\| =&
\|
\prolongateOrthArg{\restrictOrthArg{\fnprop{\Tmp}{\Tm}{\stateVar}}} + 
\prolongateArg{\restrictArg{\fnprop{\Tmp}{\Tm}{\stateVar}}}
-\crspropfore{\Tmp}{\Tm}{\stateVar}
\|\\
\leq&
\lipschitzOrth(1 + \stabilityConstant\reviewerB{\tsscrs})\|\stateVar\|
+
\lipschitzProlong\|\forecastQuantityVecArg{\timestepdummy }{\stateVar}\|,
\end{split}
\end{align}
where we have used Eq.~\eqref{eq:stabilityBoundDifferenceOne}. Then from
Eqs.~\eqref{eq:betaDefOne} and
\eqref{eq:keystabilityieq}, we have
\ourReReading{
\begin{align}
|\forecastQuantityArg{j}{\timestepdummy }{\stateVar}|\leq&|\restrictjArg{\fnprop{\Tmp}{\Tm}{\stateVar}}| +
|\crspropforesym{j}{\timestepdummy }(\stateVar)|\\
&
\leq
\lipschitzRestrictEntry{j}(1 + \stabilityConstant\tsscrs) \|\stateVar\|
+\lipschitzRestrictEntry{j}\left[\normedQuantityGappySmall{j}{\timestepdummy }
+
\normedQuantityGappy{j}{\timestepdummy }\sqrt{\memory}(1 +
\stabilityConstant\memory\tssfn)
 \right]\|\stateVar\|\\
&=
\lipschitzRestrictEntry{j}
(
1 + \normedQuantityGappySmall{j}{\timestepdummy } + 
\normedQuantityGappy{j}{\timestepdummy }\sqrt{\memory}
+(1
+\normedQuantityGappy{j}{\timestepdummy }\sqrt{\memory}
(\memory/\nfinepercoarse))
\stabilityConstant\tsscrs
)\|\stateVar\|\\
\label{eq:lastInequalityDeltaStability}&=
\lipschitzRestrictEntry{j}
(1 + \normedQuantityGappySmall{j}{\timestepdummy } + 
\normedQuantityGappy{j}{\timestepdummy }\sqrt{\memory})
(
1
+\frac{
\normedQuantityGappy{j}{\timestepdummy }\sqrt{\memory}
(\memory/\nfinepercoarse)+1}{\normedQuantityGappySmall{j}{\timestepdummy } + 
\normedQuantityGappy{j}{\timestepdummy }\sqrt{\memory}+1 }
\stabilityConstant\tsscrs
)\|\stateVar\|.
\end{align}}
Using the norm-equivalence relation $\|{\bf x}\|_2\leq \sqrt{n}\|{\bf
x}\|_\infty$, we have
$
\|\forecastQuantityVecArg{\timestepdummy }{\stateVar}\|\leq\sqrt{\nrestrict}\max_{j\in\nat{\nrestrict}}|\forecastQuantityArg{j}{\timestepdummy }{\stateVar}|.
$
Combining inequalities \eqref{eq:FGboundAgain},
\eqref{eq:lastInequalityDeltaStability}, and the above inequality yields
\ourReReading{
$$
\|\fnprop{\Tmp}{\Tm}{\stateVar} -
\crspropfore{\Tmp}{\Tm}{\stateVar}
\|\leq 
\lipschitzOrth(1 + \stabilityConstant\reviewerB{\tsscrs})\|\stateVar\|
+
\sqrt{\nrestrict}\lipschitzProlong
\max_{j\innat{\nrestrict}}\lipschitzRestrictEntry{j}(\constCoarseLFnj+1)\left(1 +
\stabilityCoarseLFnjMod\stabilityConstant\tsscrs
 \right)\|\stateVar\|,
$$
where
$\stabilityCoarseLFnjMod\defeq\frac{(\memory/\nfinepercoarse)\normedQuantityGappy{j}{\timestepdummy }\sqrt{\memory} +1}{\normedQuantityGappySmall{j}{\timestepdummy }
+
\normedQuantityGappy{j}{\timestepdummy }\sqrt{\memory}+1}
$.
 Noting that $\stabilityCoarseLFnjMod\leq 1$ because
 $\memory/\nfinepercoarse \leq 1$, this quantity can be bounded from above by the quantity that is the desired result.}
\end{proofofthing}

\begin{proofofthing}{Theorem}{thm:finalStability}
Under the stated assumptions, the results of Lemmas \ref{lem:coarseStable} and
\ref{lem:fineCoarseStable} hold.  This implies that the conditions of Lemma 1
for $\crspropsym\leftarrow\crspropforesymAll$
hold with the specified values of $\constCoarse$, $\stabilityCoarse$, $\constFineCoarse$, and $\stabilityFineCoarse$. 
\end{proofofthing}

\reviewerA{
\begin{proofofthing}{Corollary}{cor:convergenceProposed}
Under Assumptions \ref{ass:linearscalar} and \ref{ass:samebasis}, the coarse propagator is the same on
each coarse time interval, i.e., 
$
\crspropforesym{}{i} = \crspropforesym{}{j}  = \crspropforesymAll
$, and
automatically satisfies Assumption
\ref{ass:coarsePropLinear},
as the definition of the coarse propagator \eqref{eq:dataDrivenCoarse2}
simplifies to
$
\crspropforesymAll:(\stateVarScalar;\tptcrs{\timestepdummy },\tptcrs{\timestepdummy  +
	1})\mapsto
	\linearCrspropscalarLF\stateVarScalar
	$. Thus, inequality~\eqref{eq:corollaryOursConvergenceFirst} results from
	applying inequality \eqref{eq:convergeGeneral} with
	$\linearCrspropscalar\leftarrow\linearCrspropscalarLF$. Next, Assumption
	\ref{ass:coarsePropStable} is automatically satisfied under Assumption
	\ref{ass:coarseContract} because $\linearCrspropscalarLF=1+\sum_{i=1}^\memory\forecastScalarScalar{i}
	[(\linearFnpropscalar)^i-1]$. Therefore,
	inequality~\eqref{eq:corollaryOursConvergenceSecond} results from applying
	inequality \eqref{eq:convergeSuperlinear} with
	$\linearCrspropscalar\leftarrow\linearCrspropscalarLF$.
\end{proofofthing}
}

\begin{proofofthing}{Theorem}{thm:finePropGen}
Under Assumption \ref{ass:integrateCostDominant}, the wall time incurred by a
serial solution is $\nptfn\costsolve$, where $\costsolve\in\RRplus$ is the
wall time required to compute
$\fnprop{\tptfn{\timestepdummy +1}}{\tptfn{\timestepdummy }}{\stateNP(\tptfn{\timestepdummy })}$ for a given
$\timestepdummy \innatZero{\nptfn-1}$. Further, the wall time incurred by initializing the
proposed method in Steps
	\ref{step:initializeOne}--\ref{step:initializeLast} of Algorithm
	\ref{alg:parareal}
	is composed of (1) local-forecast initialization in Step
	\ref{step:initialize} of Algorithm \ref{alg:parareal}, which incurs
	$\nptcrs$ applications of the time integrator (wall time of
	$\nptcrs\costsolve$) and (2) the (worst-case) parallel fine
	propagation in Steps
	\ref{step:beginFinePropInit}--\ref{step:initializeLast}
	(wall time of $\costsolve\nptfnArg{\timestepdummy }$).
	
	Because the coarse propagator was also employed for initialization,
	Step \ref{step:initial_seed_prop} of Algorithm \ref{alg:parareal} can be
	replaced by simply setting $\coarsesolFOM{\timestepdummy +1}{0} = \apxsolFOM{\timestepdummy +1}{0}$,
	which incurs no cost under Assumption \ref{ass:integrateCostDominant}. Then,
	each subsequent iteration requires (1) serial coarse propagation in 
	Step \ref{step:serialcoarse} 
	(wall time of $(\nptcrs-\pararealit )\costsolve$), and (2) parallel fine propagation in 
        \ourReReading{Steps
	\ref{step:fineProp1} and \ref{step:fineProp2}} (wall time
	of $\costsolve\nptfnArg{\timestepdummy }$).
	The ratio of these costs yields the theoretical speedup. 
As before, we note that additional speedups may be realizable by pipelining
operations.
\end{proofofthing}

\section{Ideal conditions with Newton-solver initial
guesses}\label{sec:idealSpeedupsNewtonSolver}
We now consider the case of a nonlinear dynamical system wherein the 
forecasts are also employed to generate initial guesses for 
the Newton solver
as proposed in Ref.\
\cite{Carlberg_carlberg2015decreasing}. We therefore
introduce the following assumptions:
\begin{Assumption}[resume=assumption]
\item\label{ass:nonlinear} The velocity $\ivpfuncNP$ is nonlinear.
\item \label{ass:linearmultistep}
	The fine propagator corresponds to an implicit \reviewerA{single-step} scheme (i.e.,
$\fnpropsym$ is such that
$\resLM{i}{\stateNP(\tptfn{i}))}=0$).\footnote{\reviewerA{The algebraic residual for the
	backward-Euler method, for example, 
is $\resLMone{i}:\unknown\mapsto
\unknown - \tssfnArg{i}\fArgNoParam{\unknown}{\tptfn{i}} -
\stateNP(\tptfn{i-1})$ .}}
\item \label{ass:initialGuessLocal}
The Newton-solver initial
guesses are provided by the local forecast, i.e., Newton's method for solving
$\resLM{i}{\unknownTime}=0$ with $\tptfn{i}\in\timeInstanceSetInterval{\timestepdummy }$
at parareal iteration $\pararealit $
employs 
$$\unknownTimeInit{i} =
\prolongateArg{[\forecastFunctionmFour{\tptfn{i}}{\Tm}{\restrictEntryArg{1}{\fineFillInLocal{\timestepdummy }{\apxsolFOM{\timestepdummy }{\pararealit }}}}{1}\
\cdots\
\forecastFunctionmFour{\tptfn{i}}{\Tm}{\restrictEntryArg{\ndof}{\fineFillInLocal{\timestepdummy }{\apxsolFOM{\timestepdummy }{\pararealit }}}}{\ndof}]^T}$$
as an initial guess.
\item \label{ass:initialGuessGlobal} The Newton-solver initial
guesses are provided by the global forecast, i.e., Newton's method for solving
$\resLM{i}{\unknownTime}=0$ with $\tptfn{i}\in\timeInstanceSet$
for all parareal iterations
employs 
as an initial guess
$$\unknownTimeInit{i} =
\prolongateArg{[\forecastFunctionFour{\tptfn{i}}{0}{\restrictEntryArg{1}{\fineFillInOne{\initstateNP}}}{1}\
\cdots\\; 
\forecastFunctionFour{\tptfn{i}}{0}{\restrictEntryArg{\ndof}{\fineFillInOne{\initstateNP}}}{\ndof}]^T}.$$
\end{Assumption}
\begin{corollary}[\textit{Ideal-conditions speedup}: local-forecast initialization with Newton-solver initial guesses]\label{thm:idealSpeedup2}
If Assumptions \ref{ass:subspace}, \ref{ass:notruncate}, \ref{ass:isomorphic}, \ref{ass:localinitialize}, 
\ref{ass:nonlinear}, 
\ref{ass:linearmultistep}, and \ref{ass:initialGuessLocal} hold,
 then
	the method converges after parareal initialization (i.e.,
	$\pararealItConverge=0$ in Algorithm \ref{alg:parareal}), and only $\memory$
	nonlinear systems of algebraic equations are solved in each time interval,
	with the remaining time steps requiring only a single residual evaluation.
	Further, if Assumption \ref{ass:integrateCostDominant} holds, then
the method 
realizes a speedup of
$
\nptfn/\left((\nptcrs-1)\memory +
(\nptfnArg{\timestepdummy }-\memory)\rescost\right)
$
relative to the sequential algorithm without forecasting. Here, we denote
by $\rescost\in(0,1)$ the ratio of the computational cost of computing the
discrete residual $\resLM{j}{\unknownTime}$ relative to that of solving a
system of nonlinear algebraic equations
$\resLM{j}{\unknownTime}=0$.
\end{corollary}

\begin{figure}[htbp] 
\centering 
\subfigure[Local-forecast initialization and Newton initial guesses]{
\includegraphics[width=0.45\textwidth]{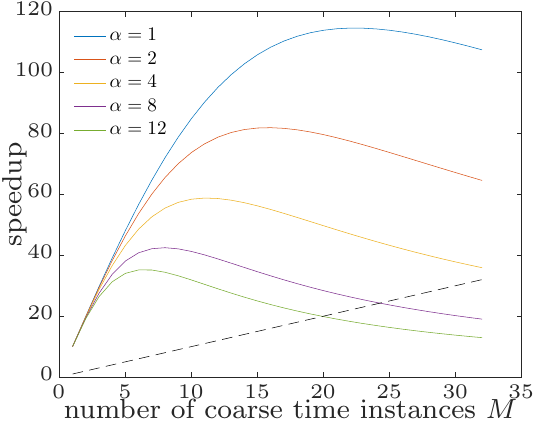} 
\label{fig:theoretical_speedup_local_newton}
}
\subfigure[Global-forecast initialization and Newton initial guesses]{
\includegraphics[width=0.45\textwidth]{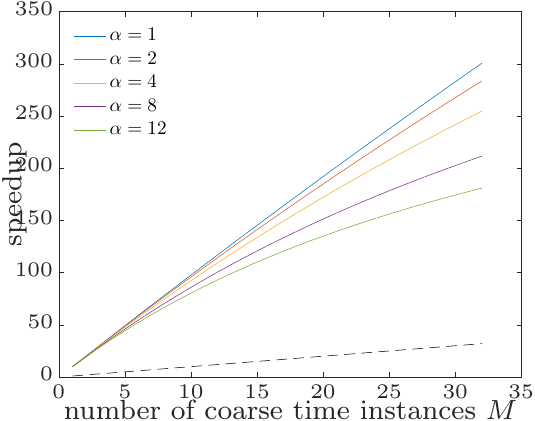} 
\label{fig:theoretical_speedup_global_newton}
}
\caption{Ideal-conditions speedup with Newton-solver initial guesses
(nonlinear dynamical systems with implicit time integration). Plot
corresponds to $\nptfn=5000$ fine time instances, setting the number of processors
equal to the number of coarse time instances $\nptcrs$. Note that the proposed method can realize super-ideal
theoretical speedups by effectively leveraging time-evolution data.}
\label{fig:theoretical_speedup_newton}
\end{figure}

Figure \ref{fig:theoretical_speedup_local_newton}
illustrates this ideal-conditions speedup.
\begin{corollary}[\textit{Ideal-conditions speedup}: global-forecast
initialization with Newton-solver initial
guesses]\label{thm:idealSpeedupGlobal2}
If Assumptions
\ref{ass:subspace},
\ref{ass:isomorphic}, 
\ref{ass:globalinitialize}, 
  \ref{ass:nonlinear}, \ref{ass:linearmultistep}, and
\ref{ass:initialGuessGlobal} hold,
 then the method converges after parareal initialization (i.e.,
	$\pararealItConverge=0$ in Algorithm \ref{alg:parareal}), only $\memory$ nonlinear systems of algebraic equations are solved
in the first time interval, and no algebraic equations are solved in the
remaining time intervals.  All remaining time steps require only a single
residual evaluation.  Further, if Assumption \ref{ass:integrateCostDominant}
holds, then the method 
realizes a theoretical speedup of
$
\nptfn/(\memory + \rescost\nptfnArg{\timestepdummy })
$
relative to the sequential algorithm without forecasting.
\end{corollary}

Figure \ref{fig:theoretical_speedup_global_newton} visualizes 
this
theoretical speedup in the case of global-forecast initialization. By
comparing this figure with
Figure \ref{fig:theoretical_speedup_global}, it is clear that employing the
forecasting approach for both initialization and initial guesses for the
Newton solver can yield \textit{super-ideal} speedups, which highlights the
potential of the proposed approach to realize near-real-time computations.

\begin{proofofthing}{Corollary}{thm:idealSpeedup2}
Under 
Assumptions \ref{ass:subspace} and \ref{ass:notruncate}, we have from Lemma \ref{lem:localSubspace} that
$\unrollmfunc{\restrictjArg{\stateNP}}\in\range{\timebasisjm}$,
$j\innat{\ndof}$, $\timestepdummy \innatZero{\nTimeIntervals-1}$.
Under Assumptions \ref{ass:subspace}, \ref{ass:notruncate},
\ref{ass:isomorphic}, and \ref{ass:localinitialize}, we have from Theorem
\ref{thm:idealSpeedup} that convergence occurs in one parareal iteration
(i.e., $\pararealItConverge=0$ in Algorithm \ref{alg:parareal}).

Then, the local forecast for the $j$th element of the restricted state to time $\tptfn{i}\in\timeInstanceSetInterval{\timestepdummy }$ during Step \ref{step:initialize} in
Algorithm \ref{alg:parareal} satisfies
\begin{align}\label{eq:localForeInitialGuess}
\begin{split}
\forecastFunctionjm{\tptfn{i}}{\Tm}{\restrictjArg{\fineFillInLocal{\timestepdummy }{\apxsolFOM{\timestepdummy }{0}}}}&=
\forecastFunctionjm{\tptfn{i}}{\Tm}{\restrictjArg{\fineFillInLocal{\timestepdummy }{\stateNP(\tptcrs{\timestepdummy })}}}
=\forecastFunctionjm{\tptfn{i}}{\Tm}{\restrictjArg{\stateNP}}\\
&=\restrictjArg{\stateNP(\Tm)} + 
\unitvec{i-\crstofn{\timestepdummy }}^T\timebasisjm
[\sampleMat{0}{\prevtpts}
\timebasisjm]^+\sampleMat{0}{\prevtpts}
\unrollmfunc{\restrictjArg{\stateNP}}\\
&=
\restrictjArg{\stateNP(\Tm)} + 
\unitvec{i-\crstofn{\timestepdummy }}^T\timebasisjm
\redCoordRestrictTwoArgs{j}{\timestepdummy }
=
\restrictjArg{\stateNP(\Tm)} + 
\unitvec{i-\crstofn{\timestepdummy }}^T\unrollmfunc{\restrictjArg{\stateNP}}=
 \restrictjArg{\stateNP(\tptfn{i})}.
%
\end{split}
\end{align}
Under Assumptions \ref{ass:nonlinear} and \ref{ass:linearmultistep}, the residual arising at each time step
is nonlinear and satisfies $\resLM{i}{\stateNP(\tptfn{i})} = 0$,
$i\innat{\nptfn}$.
Under Assumption \ref{ass:initialGuessLocal}, \ourReReading{a} local forecast \eqref{eq:localForeInitialGuess} is employed
as an initial guess for Newton's method for solving
$\resLM{i}{\unknownTime}=0$ with $\tptfn{i}\in\timeInstanceSetInterval{\timestepdummy }$ at
parareal iteration $\pararealit $; this initial guess can be expressed as
\begin{align*}
\unknownTimeInit{i} &=
\prolongateArg{[\forecastFunctionmFour{\tptfn{i}}{\Tm}{\restrictjArg{\fineFillInLocal{\timestepdummy }{\apxsolFOM{\timestepdummy }{0}}}}{1}\
\cdots\
\forecastFunctionmFour{\tptfn{i}}{\Tm}{\restrictjArg{\fineFillInLocal{\timestepdummy }{\apxsolFOM{\timestepdummy }{0}}}}{\ndof}]^T}\\
&=\prolongateArg{[\restrictEntryArg{1}{\stateNP(\tptfn{i})}\
\cdots\
\restrictEntryArg{\ndof}{\stateNP(\tptfn{i})}]^T}=
\stateNP(\tptfn{i}).
\end{align*}
Because $\resLM{i}{\unknownTimeInit{i}} = \resLM{i}{\stateNP(\tptfn{i})}=  0$
under the stated assumptions, the initial residual
is zero such that Newton's method
terminates after simply computing the initial residual. No Newton iterations
are required.

	Under Assumption \ref{ass:integrateCostDominant}, the wall time incurred by the serial solution is
	$\nptfn\costsolve$, where $\costsolve$ in this case corresponds to solving a
	system of nonlinear algebraic equations. The wall time incurred by the
	proposed approach is composed of (1) the serial coarse propagation in \ourReReading{Step}
	\ref{step:initialize} of Algorithm \ref{alg:parareal}, which entails solving
	$\memory$ systems of nonlinear algebraic equations in each time interval
	(i.e., $(\nptcrs-1)\memory\costsolve$) and (2) the (worst-case) parallel
	fine propagation in \ourReReading{Steps}
	\ref{step:beginFinePropInit}--\ref{step:initializeLast}, which no longer requires
	solving systems of nonlinear algebraic equations; it entails computing only
	a single residual for each remaining time instance (i.e.,
	$(\nptfnArg{\timestepdummy }-\memory)\rescost\costsolve$). The ratio of these
	costs yields the theoretical speedup. Again, additional speedups may be
	realizable by pipelining operations.
\end{proofofthing}

\begin{proofofthing}{Corollary}{thm:idealSpeedupGlobal2}
Under Assumptions \ref{ass:subspace}, \ref{ass:isomorphic}, and
\ref{ass:globalinitialize}, we have from
Theorem \ref{thm:idealSpeedupGlobal} that $\unrollfunc{\restrictjArg{\stateNP}}=\timebasisj\redCoordRestrictArgs{j}$ for some
$\redCoordRestrictArgs{j}\in\RR{\dimBasisj}$ (see Eq.~\eqref{eq:projectionhofrj}) and convergence in one parareal
iteration (i.e., $\pararealItConverge=0$ in Algorithm \ref{alg:parareal}).
Then, the global forecast for the $j$th element of the restricted state to
time $\tptfn{i}\in\timeInstanceSet$ during Step \ref{step:initialize} in
Algorithm \ref{alg:parareal} satisfies
(from Eqs.~\eqref{eq:genSolutionGlobalForecast} and \eqref{eq:globalForecast})
 \begin{align} \label{eq:globalForecastProof}
 \begin{split}
\forecastFunctionjNo\forecastFunctionjArgs{\tptfn{i}}{0}{\restrictEntryArg{j}{\fineFillInOne{\initstateNP}}}
=\forecastFunctionjNo\forecastFunctionjArgs{\tptfn{i}}{0}{\restrictjArg{\stateNP}}
&=\restrictjArg{\stateNP(0)} +
\unitvec{i}^T\timebasisj
[\sampleMat{0}{\prevtpts}
\timebasisj]^+\sampleMat{0}{\prevtpts}
\timebasisj\redCoordRestrictArgs{j}
=\restrictjArg{\stateNP(0)} +
\unitvec{i}^T\timebasisj
\redCoordRestrictArgs{j}\\
&=\restrictjArg{\stateNP(0)} +
\unitvec{i}^T\unrollfunc{\restrictjArg{\stateNP}} = 
\restrictjArg{\stateNP(\tptfn{i})}.
	\end{split}
  \end{align} 
Under Assumptions \ref{ass:nonlinear} and \ref{ass:linearmultistep}, the residual arising at each time step
is nonlinear and satisfies $\resLM{i}{\stateNP(\tptfn{i})} = 0$,
$i\innat{\nptfn}$.
Under Assumption \ref{ass:initialGuessGlobal}, global forecast
\eqref{eq:globalForecastProof}
is employed
as an initial guess for Newton's method for solving
$\resLM{i}{\unknownTime}=0$ with $\tptfn{i}\in\timeInstanceSetInterval{\timestepdummy }$ at
all parareal iterations. This initial guess can be written as
\begin{align*}
\unknownTimeInit{i} =
\prolongateArg{[\forecastFunctionFour{\tptfn{i}}{0}{\restrictEntryArg{1}{\fineFillInOne{\initstateNP}}}{1}\
\cdots\
\forecastFunctionFour{\tptfn{i}}{0}{\restrictEntryArg{\ndof}{\fineFillInOne{\initstateNP}}}{\ndof}]^T}
=\prolongateArg{[\restrictEntryArg{1}{\stateNP(\tptfn{i})}\
\cdots\
\restrictEntryArg{\ndof}{\stateNP(\tptfn{i})}]^T}
=\stateNP(\tptfn{i}).
\end{align*}
Because $\resLM{i}{\unknownTimeInit{i}}=\resLM{i}{\stateNP(\tptfn{i})} = 0$
under the stated assumptions, the initial residual is zero so that Newton's
method terminates after simply computing the initial residual without any
required Newton iterations.

Under Assumption \ref{ass:integrateCostDominant}, the wall time incurred by the serial solution is
$\nptfn\costsolve$, where $\costsolve$ in this case corresponds to solving a
	system of nonlinear algebraic equations. The wall time incurred by the proposed forecasting
approach is composed of (1) solving $\memory$ systems of nonlinear algebraic
equations in the first
coarse time interval for initialization in Step \ref{step:initialize} of
Algorithm \ref{alg:parareal} (i.e., $\memory\costsolve$) and (2) the
worst-case parallel fine propagation in Steps
\ref{step:beginFinePropInit}--\ref{step:initializeLast}, which no longer requires
solving any linear algebraic systems of equations; it
entails computing only a single residual for all remaining time instances
(i.e., $\nptfnArg{\timestepdummy }\rescost\costsolve$). The ratio of these costs
yields the theoretical speedup.
\end{proofofthing}
Thus, the method realizes super-ideal speedups when the
	local-forecast is applied to Newton-solver initial guesses in the nonlinear
	case (Corollaries \ref{thm:idealSpeedup2}--\ref{thm:idealSpeedupGlobal2} and
	Figure \ref{fig:theoretical_speedup_newton})
\end{document}